\documentclass[a4paper,11pt]{article}
\usepackage{latexsym}
\usepackage{amsmath}
\usepackage{amsthm}

\usepackage[dvipdfmx]{graphicx}
\usepackage{txfonts}
\usepackage{bm}
\usepackage{color}
\usepackage{epic,eepic}

\usepackage{geometry}
\numberwithin{equation}{section}

\newtheorem{theorem}{Theorem}[section]
\newtheorem{proposition}[theorem]{Proposition}
\newtheorem{corollary}[theorem]{Corollary}
\newtheorem{lemma}[theorem]{Lemma}
\newtheorem{claim}[theorem]{Claim}
\newtheorem{remark}{Remark}[section]
\newtheorem{example}{Example}[section]
\newtheorem{algorithm}[theorem]{Algorithm}

\newcommand{\OMIT}[1]{{\bf [OMIT:} #1 \ {\bf --- end OMIT] }}  
   \renewcommand{\OMIT}[1]{}            

\newcommand{\RR}{{\bf R}}
\newcommand{\ZZ}{{\bf Z}}

\newcommand{\finbox}{\hspace*{\fill}$\rule{0.17cm}{0.17cm}$}
\newcommand{\finboxHere}{\ $\rule{0.17cm}{0.17cm}$}

\newcommand{\smallbox}{\scalebox{0.6}{\mbox{$\square$}}}
\newcommand{\sq}{\sp{\smallbox}}

\newcommand{\llceil}{\left\lceil} 
\newcommand{\rrceil}{\right\rceil} 

\newcommand{\llfloor}{\left\lfloor} 
\newcommand{\rrfloor}{\right\rfloor}

\newcommand{\ol}{\overline}

\newcommand{\ora}{\stackrel{\rightarrow}}

\newcommand{\odotZ}{\overset{....}} 

\newcommand{\odottZ}{\overset{....}}  

\newcommand{\Proof}{\noindent {\bf Proof.  }}
\begin{document}

\title{Discrete Decreasing Minimization, Part I:
\break Base-polyhedra
with Applications in Network Optimization}

\author{Andr\'as Frank\thanks{MTA-ELTE Egerv\'ary Research Group,
Department of Operations Research, E\"otv\"os University, P\'azm\'any
P. s. 1/c, Budapest, Hungary, H-1117. 
e-mail: {\tt frank\char'100 cs.elte.hu}. 
The research was partially supported by the
National Research, Development and Innovation Fund of Hungary
(FK\_18) -- No. NKFI-128673.
}
\ \ and \ 
{Kazuo Murota\thanks{Department of Economics and Business Administration,
Tokyo Metropolitan University, Tokyo 192-0397, Japan, 
e-mail: {\texttt{} murota\char'100 tmu.ac.jp}. 
The research was supported by CREST, JST, Grant Number
JPMJCR14D2, Japan, and JSPS KAKENHI Grant Number 26280004.  }}}

\date{August 2018 / May 2019 / July 2019}


\maketitle


\newpage

\begin{abstract} Motivated by resource allocation problems, 
Borradaile et al.~(2017) investigated orientations of an undirected graph in
which the sequence of in-degrees of the nodes, when arranged in a
decreasing order, is lexicographically minimal in the sense that the
largest in-degree is as small as possible, within this, the next
largest in-degree is as small as possible, and so on.  
They called such an orientation egalitarian but we prefer to use the term 
{\it decreasingly minimal} ($=$dec-min) to avoid confusion with another
egalitarian-felt orientation where the smallest in-degree is as large
as possible, within this, the next smallest in-degree is as large as
possible, and so on.  Borradaile et al. proved that an orientation is
dec-min if and only if 
there is no dipath for which the in-degree of its last
node is at least two larger than the in-degree of its first node. 
They conjectured that an
analogous statement holds for strongly connected dec-min orientations, as well.  
We prove not only this conjecture but its extension to
$k$-edge-connected orientations, as well, even if additional in-degree
constraints are imposed on the nodes.

Resource allocation was also the motivation behind an earlier
framework by Harvey et al.~(2006) who introduced and investigated
semi-matchings of bipartite graphs.  
As a generalization of their results, 
we characterize degree-constrained subgraphs of a bipartite
graph $G=(S,T;E)$ which have a given number of edges and their
degree-sequence in $S$ is decreasingly minimal.  
We also provide a
solution to a discrete version of Megiddo's \lq lexicographically\rq \
optimal (fractional) network flow problem (1974, 1977).

Furthermore, we exhibit a generalization of a result of Levin and Onn
(2016) on \lq shifted\rq \ matroid optimization, and describe a way of
finding a basis of each of $k$ matroids so that the sum of their
incidence vectors is decreasingly minimal.

Our main goal is to integrate these cases into a single framework.
Namely, we characterize dec-min elements of an M-convex set (which is
nothing but the set of integral points of an integral
base-polyhedron), and prove that the set of dec-min elements is a
special M-convex set arising from a matroid base-polyhedron by
translation.  The topic of our investigations may be interpreted as a
discrete counter-part of the work by Fujishige (1980) on the (unique)
lexicographically optimal base of a base-polyhedron.  On the dual
side, as an extension of a result of Borradaile et al.  (2018) on
density decomposition of networks, we exhibit a canonical chain (and
partition) associated with a base-polyhedron.  
We also show that
dec-min elements of an M-convex set are exactly those which minimize
the square-sum of components, and describe a new min-max formula for
the minimum square-sum.

Our approach gives rise to a strongly polynomial algorithm for
computing a dec-min element, as well as the canonical chain.  
The algorithm relies on a submodular function minimizer oracle in the
general case, which can, however, be replaced by more efficient
classic flow- and matroid algorithms in the relevant special cases.

This paper constitutes the first part of a series of our papers on 
discrete decreasing minimization.
In Part~II, we offer a broader structural view from discrete convex analysis (DCA).  
In particular, min-max formulas will be derived as
special cases of general DCA results.  Furthermore, the relationship
between continuous and discrete problems will also be clarified.
In Part~III
we describe the structure of decreasingly minimal integral feasible flows and
develop a strongly polynomial algorithm for finding such a dec-min flow.
In Part~IV we consider
the discrete decreasing minimization problem over the
intersection of two base-polyhedra, and also over submodular flows.
Finally, Part~V deals with discrete decreasing minimality 
with respect to a weight vector.
\end{abstract}


\newpage

\tableofcontents

\newpage


\section{Introduction}

\subsection{Background problems}

There are three independent sources of the topic we study.

\subsubsection{Orientations of graphs}

Let $G=(V,E)$ be an undirected graph.  
Orienting an edge $e=uv$ means the operation that replaces $e$ 
by one of the two oppositely directed edges 
(sometimes called arcs) $uv$ or $vu$.  
A directed graph arising from $G$ by orienting all of its edges 
is called an orientation of $G$.  
A graph orientation problem consists of finding an orientation
of $G$ meeting some specified properties such as in-degree constraints
(lower and upper bounds) and/or various connectivity prescriptions.
One goal is to characterize undirected graphs for which the requested
orientation exists, and a related other one is to design an algorithm
for finding the orientation.  The first orientation result is due to
Robbins \cite{Robbins} who proved that exactly the 2-edge-connected
graphs admit strongly connected (strong, for short) orientations.  
The in-degree constrained orientation problem is equivalent to an
in-degree constrained subgraph problem for directed graphs.  
An explicit characterization was formulated by Frank and Gy\'arf\'as
\cite{FrankP2}, but this may be considered as an appropriate
adaptation (or reformulation) of a characterization of
degree-constraints for subgraphs of a digraph 
(\cite{Ford-Fulkerson}, Theorem 11.1).  
The proper novelty of \cite{FrankP2} was a solution to
an amalgam of the two orientation problems above when one is
interested in the existence of a strong orientation which, in
addition, complies with upper and lower bounds on the in-degrees of nodes.

The literature is quite rich in orientation results, for a relatively
wide overview, see the book \cite{Frank-book}.  There are, however,
other type of requirements for an orientation of $G$ where, rather
than having prescribed upper and lower bounds for the in-degrees of
nodes (or, sometimes, beside these bounds), one is interested 
in the global distribution of the in-degrees of nodes.  
That is, the goal is to find orientations 
(with possible connectivity expectations) 
whose in-degree vector (on the node-set) is felt intuitively evenly distributed:  
\lq fair\rq, \lq equitable\rq, \lq egalitarian\rq,  \lq uniform\rq.  
For example, how can one determine 
the minimum value $\beta_{1}$ of the largest
in-degree of a ($k$-edge-connected) orientation?  
Even further, after determining $\beta_{1}$, it might
be interesting to minimize the number of nodes 
with in-degree $\beta_{1}$ among orientations of $G$ 
with largest in-degree $\beta_{1}$.  
Or, a more global equitability feeling is captured if we minimize the sum
of squares of the in-degrees.  
For example, the in-degree sequence $(1,2, 6,8)$ 
(with square-sum 105) is felt less \lq fair\rq\ or \lq egalitarian\rq\ or
\lq evenly distributed\rq\ than $(4,4,4,5)$ \ (with square-sum 73).
For an in-degree distribution,
yet another natural way to measure
its deviation from uniformity
is to consider the sum 
$\sum [ |\varrho(v) - \varrho(u) | : u,v, \in V  ]$,
where $\varrho(v)$ denotes the in-degree of $v$.
One feels an orientation `fair' if this sum is as small as possible.

A formally different definition was recently suggested and
investigated by Borradaile et al.~\cite{BIMOWZ} who called an
orientation of $G=(V,E)$ {\bf egalitarian} if the highest in-degree of
the nodes is as small as possible, and within this, the second highest
(but not necessarily distinct) in-degree is as small as possible, and
within this, the third highest in-degree is as small as possible, and
so on.  In other words, if we rearrange the in-degrees of nodes in a
decreasing order, then the sequence is lexicographically minimal.  
In order to emphasize that the in-degrees are considered in a decreasing
order, we prefer the use of the more expressive term 
{\bf decreasingly minimal} ({\bf dec-min}, for short) 
for such an orientation, rather than egalitarian.

This change of terminology is reasonable since one may also consider
the mirror problem of finding an {\bf increasingly maximal}
 (or {\bf inc-max}, for short) orientation that is an orientation of $G$ 
in which \ the smallest in-degree is as large as possible, within this,
the second smallest in-degree is as large as possible, and so on.
Intuitively, such an orientation may equally 
be felt \lq egalitarian\rq \ in the informal meaning of the word.

Borradaile et al.~\cite{BIMOWZ}, however, proved that 
an orientation of a graph is
decreasingly minimal (egalitarian in their original term) 
if and only if 
there is no \lq small\rq \ improvement, where a small improvement
means the reorientation of a dipath from some node $s$ to another node
$t$ with in-degrees $\varrho(t)\geq \varrho(s)+2$.
This theorem immediately implies that an orientation is decreasingly minimal 
if and only if 
it is increasingly maximal, and therefore we could retain the
original terminology \lq egalitarian orientation\rq \ used in \cite{BIMOWZ}.

However, when orientations are considered with specific requirements
such as strong (or, more generally, $k$-edge-) connectivity and/or
in-degree bounds on the nodes, the possible equivalence of
decreasingly minimal and increasingly maximal orientations had not yet
been investigated.  Actually, Borradaile et al. conjectured that a
strong orientation of a graph is decreasingly minimal (among strong orientations) 
if and only if 
there is no small improvement preserving
strong connectivity, and this, if true, would imply immediately that
decreasing minimality and increasing maximality do coincide for strong
orientations, as well.

We shall prove this conjecture in its extended form concerning
$k$-edge-connected and in-degree constrained orientations.  
This result implies immediately that the notions of decreasing minimality
and increasing maximality coincide for $k$-edge-connected and
in-degree constrained orientations, as well.  
We hasten to emphasize that this coincidence is not at all automatic or inevitable.  
For example, although Robbins' theorem on strong orientability of
undirected graphs nicely extends to mixed graphs, as was pointed out
by Boesch and Tindell \cite{BT}, it is not true anymore that a
decreasingly minimal strong orientation of a mixed graph is always
increasingly maximal.  
(For a counterexample, see Section \ref{SCstroricntex}.)
This discrepancy may be explained by the fact that the in-degree
vectors of strong orientations of a graph form an M-convex set while
the set of in-degree vectors of strong orientations of a mixed graph
is not M-convex anymore:  it is the intersection of two M-convex sets.
(The definition of an M-convex set was mentioned in the Abstract and
will be introduced formally in Section \ref{egal2base}.)

Interestingly, it will turn out that for $k$-edge-connected and
in-degree constrained orientations of undirected graphs not only
decreasing minimality and increasing maximality coincide but such
orientations are exactly those minimizing the square-sum of the
in-degrees of the nodes.

\subsubsection{A resource allocation problem and network flows}
\label{SCresource}

Another source of our investigations is due to Harvey et al.~\cite{HLLT} 
who solved the problem of minimizing 
$\sum [d_{F}(s) (d_{F}(s)+1):  s\in S]$ 
over the semi-matchings $F$ of a simple
bipartite graph $G=(S,T;E)$.  
Here a semi-matching is a subset $F$ of
edges for which $d_{F}(t)=1$ holds for every node $t\in T$.  
Harada et al.~\cite{HOSY07} solved the minimum cost version of this problem.
The framework of Harvey et al. was extended by Bokal et al.~\cite{BBJ} to
quasi-matchings, and, even further, to degree-bounded quasi-matchings
by Katreni{\v c} and Semani{\v s}in \cite{Katrenic13}.
It turns out that these problems are strongly related to minimization of a
separable convex function over (integral elements of) a
base-polyhedron which has been investigated in the literature under
the name of \lq resource allocation problems under submodular constraints\rq 
(\cite{Fed-Gro86}, \cite{Hochbaum-Hong}, \cite{Hoc07},
\cite{KSI}, \cite{Ibaraki-Katoh.88}, \cite{Katoh-Ibaraki.98}).  
Ghodsi et al.~\cite{GZSS} considered the problem of finding a semi-matching
$F$ of $G=(S,T;E)$ whose degree-vector restricted to $S$ is
increasingly maximal.  
This problem, under the name \lq constrained max-min fairness\rq \ originated 
from a framework to model a fair
sharing problem for datacenter jobs.  
Here $T$ corresponds to the set of available computers while $S$ to the set of users.  
An edge $st$ belongs to $E$ if user $s$ can run her program on computer $t$.
Ghodsi et al. also consider the fractional version when, instead of
finding a semi-matching of $G$, one is to find a real vector
$x:E\rightarrow {\bf R}_{+}$ so that $d_{x}(t)=1$ for every $t\in T$ and
the vector $(d_{x}(s):  s\in S)$ is increasingly maximal. 
 (Here $d_{x}(v):= \sum [ x(uv):  uv\in E]$).  
When $x$ is requested to be $(0,1)$-valued, we are back at the subgraph version.

It should be emphasized that, unlike the well-known situation with
ordinary bipartite matchings or $b$-matchings, in this problem the
optima for the subgraph version and for the fractional version may be
different.  For example, if $T$ consists of a single node $t$,
$S=\{s_{1},s_{2}\}$, and $E=\{s_{1} t, s_{2} t\}$, 
then the original subgraph problem has two inc-max solutions:  
$F_{1}=\{s_{1} t\}$ and $F_{2}=\{s_{2} t\}$,
where the degree-vector in $S$ is $(1,0)$ in the first case, and
$(0,1)$ in the second case. 
 (Note that both $F_{1}$ and $F_{2}$ are also dec-min in $S$.)  
On the other hand, in the fractional version there
is a unique inc-max fractional solution: 
 $x(s_{1} t)=1/2$ and $x(s_{2} t)=1/2$
 (and this happens to be the unique dec-min solution).
Here the \lq fractional degree-vector\rq\ of $x$ in $S$ is $(1/2,1/2)$.  
Obviously, the fractional vector $(1/2,1/2)$ is
decreasingly smaller (and increasingly larger) than $(1,0)$.

We shall solve the following generalization of the subgraph problem.
Suppose that we are also given a positive integer $\gamma $, a lower
bound function $f:V\rightarrow {\bf Z}_{+}$ and an upper bound function
$g:V\rightarrow {\bf Z}_{+}$ with $f\leq g$ where $V:=S\cup T$.  
The problem is to find a subgraph $H=(S,T;F)$ of $G$ 
with $\vert F\vert =\gamma$ 
for which $f(v)\leq d_{F}(v)\leq g(v)$ for every node $v\in V$
and the degree-vector of $H$ in $S$ \ ({\bf !}) \ is increasingly
maximal (or decreasingly minimal).  
It will turn out that in this case
a solution is dec-min 
if and only if
it is inc-max.  
We emphasize that
in this problem the roles of $S$ and $T$ are not symmetric since we
require that the restriction on $S$ of the degree-vector of the
degree-constrained subgraph $H$ in $S$ should be decreasingly minimal.
The symmetric version when the degree-vector requested to be dec-min
on the whole node-set $V=S\cup T$ is definitely more difficult, and
will be solved in \cite{Frank-Murota.4}.  
An explanation for this
difference of the two seemingly quite similar problems is that the
first one may be viewed as a problem concerning a single
base-polyhedron while the second one may be viewed as a problem on the
intersection of two base-polyhedra.

We also solve the degree-constrained subgraph problem when the subgraph, 
in addition, is requested to have a maximum number of edges
(not just an arbitrarily prescribed number $\gamma $), and we are
interested in such a subgraph for which its degree vector 
on the whole node-set $S\cup T$ is decreasingly minimal.

\medskip

There is a much earlier, strongly related problem concerning network flows, 
due to Megiddo \cite{Megiddo74}, \cite{Megiddo77}.  
We are given a digraph $D=(V,A)$ with a source-set $S\subset V$ 
and a single sink-node $t\in V-S$. 
 (The more general case, when the sink-set $T\subseteq V-S$ may consist 
of more than one node can easily be reduced to the special case when $T=\{ t \}$). 
Let $g:A\rightarrow {\bf R}_{+}$ be a capacity function.  
By an $St$-flow, or just a flow, we mean a function 
$x:A\rightarrow {\bf R}_{+}$ 
for which the net out-flow
$\delta_{x}(v)-\varrho_{x}(v)=0$ if $v\in V-(S+t)$ and 
$\delta_{x}(v)-\varrho_{x}(v)\geq 0$ if $v\in S$.
 (Here $\varrho_{x}(v):= \sum [x(uv):uv\in A]$ and 
$\delta_{x}(v):= \sum [x(vu):  vu\in A]$.) 
The flow is feasible if $x\leq g$.  
The {\bf flow  amount} of $x$ is defined by $\varrho_{x}(t)-\delta_{x}(t)$.
Megiddo solved the problem of finding a feasible flow of maximum flow amount which is,
with his term, \lq source-optimal\rq \ at $S$.  
Source-optimality is the same as requiring that the net out-flow vector 
on $S$ is increasingly maximal.  
Note that there is a formally different but
technically equivalent version of Megiddo's problem 
when $S=\{s\}$, \  $T=\{t\}$ 
(with no arcs entering $s$) and we are interested in finding
a feasible flow $x$ of maximum amount for which the restriction of $x$
onto the set of arcs leaving $s$ is increasingly maximal. 
It must be emphasized that the flow in Megiddo's problem 
is not requested to be integer-valued.

The integrality property is a fundamental feature of ordinary network flows.  
It states that in case of an integer-valued capacity function
$g$ there always exists a maximum flow which is integer-valued.  
In this light, it is quite surprising that the integer-valued
 (or discrete) version of Megiddo's inc-max problem 
(source-optimal with his term), 
when the capacity function $g$ is integer-valued and the
max flow is required to be integer-valued, has not been investigated
in the literature, and we consider the present work as the first such attempt.

We solve this discrete version of Megiddo's problem in a more general
form concerning base-polyhedra (or M-convex sets) and this approach
gives rise to a strongly polynomial algorithm.

\subsubsection{Matroid bases}

A third source of discrete decreasing minimality problems is due to
Levin and Onn \cite{Levin-Onn} who used the term \lq shifted optimization\rq.  
They considered the following matroid optimization problem.  
For a specified integer $k$, find $k$ bases 
$Z_{1},Z_{2},\dots ,Z_{k}$ of a matroid $M$ on $S$ in such a way that the vector 
$\sum_{i}\chi_{Z_{i}}$ be, in our terms, decreasingly minimal, where $\chi_{Z}$
is the incidence (or characteristic) vector of a subset $Z$.  
They apply the following natural approach to reduce the problem to classic
results of matroid theory.  
Replace first each element $s$ of $S$ by
$k$ copies to be parallel in the resulting matroid $M'$ on the new ground-set 
$S'=S_{1}\cup S_{2}\cup \cdots \cup S_{k}$ 
where $S_{1},\dots ,S_{k}$
are the $k$ copies of $S$.  
Assign then a \lq rapidly increasing\rq \  cost function to the copies.  
(The paper \cite{Levin-Onn} explicitly
describes what rapidly increasing means).  
Then a minimum cost basis of the matroid $M_{0}$ obtained 
by multiplying $M'$ $k$-times will be a solution to the problem.  
(By definition, a basis of $M_{0}$
 is the union of $k$ disjoint bases of $M'$).

Our goal is to provide a solution to a natural generalization of this
problem when $k$ matroids $M_{1},\dots ,M_{k}$ are given on the common
ground-set $S$ and we want to select a basis $Z_{i}$ of each matroid $M_{i}$ 
in such a way that $\sum_{i} \chi_{Z_{i}}$ 
should be decreasingly minimal.  
The approach of Levin and Onn does not seem to work in this more general setting.  
We can even prescribe upper and lower bounds on
the elements $s$ of $S$ to constrain the number of the bases containing $s$.

\subsection{Main goals}

Each of the three problems above may be viewed as a special case of a
single discrete optimization problem:  
characterize decreasingly minimal elements of an M-convex set 
(or, in other words, dec-min integral elements of a base-polyhedron).  
By one of its equivalent definitions, an M-convex set is nothing but
the set of integral elements of an integral base-polyhedron.
The notion was introduced and investigated by
Murota \cite{Murota98a}, \cite{Murota03}.

We characterize dec-min elements of an M-convex set as those admitting
no local improvement, and prove that the set of dec-min elements is itself 
an M-convex set arising by translating a matroid base-polyhedron with an integral vector.  
This result implies that
decreasing minimality and increasing maximality coincide for M-convex
sets.  We shall also show that an element of an M-convex set is
dec-min precisely if it is a square-sum minimizer.  
Using the characterization of dec-min elements, we shall derive a novel min-max
theorem for the minimum square-sum of elements of 
an integral member of a base-polyhedron.  
Furthermore, we describe a strongly polynomial
algorithm for finding a dec-min element.  
The algorithm relies on a
subroutine to minimize a submodular function but in the special cases
mentioned above this general routine can be replaced by known strongly
polynomial network flow and matroid algorithms.

The structural description of the set of dec-min elements of an M-convex set 
(namely, that this set is a matroidal M-convex set) 
makes it possible to solve the algorithmic problem of 
finding a minium cost dec-min element.  
(In the continuous case this problem simply did not exist 
due to the uniqueness of the fractional dec-min element of a base-polyhedron.)  
We shall also describe a polynomial algorithm for
finding a minimum cost (in-degree constrained) dec-min orientation.
Furthermore, we shall outline an algorithm to solve the minimum cost
version of the resource allocation problem of Harvey et al.~\cite{HLLT}
mentioned at the beginning of Section \ref{SCresource}.
Furthermore, as an essential extension of the algorithm of  Harada et al.~\cite{HOSY07}, 
we describe a strongly polynomial algorithm to solve 
the minimum cost version of the decreasingly minimal degree-bounded subgraph problem.

The topic of our investigations may be interpreted as a discrete
counter-part of the work by Fujishige \cite{Fujishige80} from 1980 on
the lexicographically optimal base of a base-polyhedron $B$, where
lexicographically optimal is essentially the same as decreasingly minimal.  
He proved that there is a unique lexicographically optimal member $x_{0}$
of $B$, and $x_{0}$ is the unique minimum norm 
(that is, the minimum square-sum) element of $B$.  
This uniqueness result reflects a
characteristic difference between the behaviour of the fractional and
the discrete versions of decreasing minimization since in the latter
case the set of dec-min elements (of an M-convex set) is typically not
a singleton, and it actually has, as indicated above, 
a matroidal structure.

Fujishige also introduced the concept of principal partitions
concerning the dual structure of the minimum norm point of a base-polyhedron.  
Actually, he introduced a special chain of the
subsets of ground-set $S$ and his principal partition arises by taking
the difference sets of this chain.  
We will prove that there is an
analogous concept in the discrete case, as well.  
As an extension of
the above-mentioned elegant result of Borradaile et al.~\cite{BMW}
concerning graphs, we show that there is a canonical chain describing
the structure of dec-min elements of an M-convex set.  
We will point out in Part~II \cite{Frank-Murota.2}, that Fujishige's principal
partition is a refinement of our canonical partition.  
Our approach gives rise to a combinatorial algorithm to compute the canonical chain.

\paragraph{Outlook}

The present work is the first member of a five-partite series.  
In Part~II \cite{Frank-Murota.2}, we offer a broader structural view from
discrete convex analysis (DCA).  
In particular, min-max formulas will be derived as special cases of general DCA results.  
Furthermore, the relationship between continuous and discrete problems will also be clarified.  
In Part~III \cite{Frank-Murota.3} 
we describe the structure of decreasingly minimal integral feasible flows and
develop a strongly polynomial algorithm for finding such a dec-min flow.
Part~IV \cite{Frank-Murota.4} 
describes a strongly polynomial algorithm for 
the discrete decreasing minimization
problem over the intersection of two M-convex sets, and also over the
integral elements of an integral submodular flow.  
One of the motivations behind these investigations was the observation mentioned
earlier that, for strong orientations of mixed graphs, dec-min
orientations and inc-max orientations do not coincide.  
The reason behind this phenomenon is that the set of in-degree vectors of strong
orientations of a mixed graph is not an M-convex set anymore.  
It is, in fact, the intersection of two M-convex sets, and therefore the
results of Part~IV can be used to solve this special case, as well.
Finally, in Part~V \cite{Frank-Murota.5}
we consider discrete decreasing minimality with respect to a weight vector.

\subsection{Notation}

Throughout the paper, $S$ denotes a finite non-empty ground-set.  
For elements $s,t\in S$, we say that $X\subset S$ is an {\bf $s\ol t$-set}
if $s\in X\subseteq S-t$.  
For a vector $m\in {\bf R}\sp{S}$ (or function $m:S\rightarrow {\bf R}$), 
the restriction of $m$ to $X\subseteq S$ is denoted by $m\vert X$.  
We also use the notation
$\widetilde m(X)=\sum [m(v):  v\in X]$.  
With a small abuse of notation, 
we do not distinguish between a one-element set $\{s\}$
called a singleton and its only element $s$.  
When we work with a chain 
$\cal C$ of non-empty sets $C_{1} \subset C_{2} \subset \cdots \subset C_{q}$, 
we sometimes use $C_{0}$ to denote the empty set without assuming
that $C_{0}$ is a member of ${\cal C}$.

Two subsets $X$ and $Y$ of $S$ are {\bf intersecting} if 
$X \cap Y\not =\emptyset$ and {\bf properly intersecting} if none of 
$X-Y$, $Y-X$ and $X\cap Y$ is empty.  
If, in addition, $S-(X\cup Y)\not =\emptyset $, the two sets are {\bf crossing}.

We assume that the occurring graphs or digraphs have no loops but
parallel edges are allowed.  
For a digraph $D=(V,A)$, the {\bf in-degree} of a node $v$ is 
the number of arcs of $D$ with head $v$.
The in-degree $\varrho_{D}(Z)=\varrho(Z)$ of a subset $Z\subseteq V$
denotes the number of edges ($=$ arcs) entering $Z$, where an arc $uv$
is said to enter $Z$ if its head $v$ is in $Z$ while its tail $u$ is in $V-Z$.

The {\bf out-degree} $\delta_{D}(Z)=\delta(Z)$ 
is the number of arcs leaving $Z$, that is $\delta(Z)=\varrho(V-Z)$.  
The number of edges of a directed or undirected graph $H$
induced by $Z\subseteq V$ is denoted by $i(Z)=i_H(Z)$.  
In an undirected graph $G=(V,E)$, the {\bf degree} $d(Z)=d_{G}(Z)$ of a subset
$Z\subseteq V$ denotes the number of edges connecting $Z$ and $V-Z$
while $e(Z)=e_{G}(Z)$ denotes the number of edges with one or two end-nodes in $Z$.  
Clearly,  $e(Z)=d(Z)+i(Z)$.

The characteristic (or incidence) vector of $Z$ is denoted by $\chi_{Z}$, 
that is,  $\chi_{Z}(v)=1$ if $v\in Z$ and $\chi_{Z}(v)=0$ otherwise.

For a polyhedron $B$, $\odotZ{B}$ 
 (pronounce:  dotted $B$)
denotes the set of integral members
(elements, vectors, points) of $B$.

For a set-function $h$, we allow it to have value $+\infty $ or $-\infty
$. Unless otherwise stated, $h(\emptyset )=0$ is assume throughout.
Where $h(S)$ is finite, the {\bf complementary function} $\ol h$ is
defined by $\ol h(X) = h(S) - h(S-X)$.  Observe that the complementary
function of $\ol h$ is $h$ itself.

Let $b$ be a set-function for which $b(X)=+\infty $ is allowed but
$b(X)=-\infty $ is not.  The submodular inequality for subsets
$X,Y\subseteq S$ is defined by 
\[
b(X) + b(Y) \geq b(X\cap Y) + b(X\cup Y).
\] 
We say that $b$ is submodular if the submodular inequality holds
for every pair of subsets $X, Y\subseteq S$ with finite $b$-values.
When the submodular inequality is required only for intersecting
(crossing) pairs of subsets, we say that $b$ is {\bf intersecting}
({\bf crossing}) {\bf submodular}.  
When we say that a function $b$ is
submodular, this formally always means that $b$ is fully submodular,
but to avoid any misunderstanding we sometimes emphasize this by
writing (fully) submodular, in particular in an environment where
intersecting or crossing submodular functions also show up.

A set-function $p$ is (fully, intersecting, crossing) supermodular if
$-p$ is (fully, intersecting, crossing) submodular.  
When a submodular function $b$ and a supermodular function $p$ meet 
the {\bf cross-inequality}
\[
b(X)-p(Y) \geq b(X-Y) - p(Y-X)
\]
 for every pair $X,Y\subseteq S$, we say that $(p,b)$ 
is a {\bf paramodular pair} (a {\bf strong pair}, for short).  
If $b$ is intersecting submodular, $p$ is intersecting
supermodular, and the cross-inequality holds for intersecting pairs of
sets, we speak of an {\bf intersecting paramodular pair} (or for
short, a {\bf weak pair}).

For functions $f:S\rightarrow {\bf Z}\cup \{-\infty \}$ and
$g:S\rightarrow {\bf Z}\cup \{+\infty \}$ 
with $f\leq g$,  the polyhedron 
$T(f,g)=\{x\in {\bf R}\sp{S}:  f\leq x\leq g\}$
is called a {\bf box}.  
If $g(s)\leq f(s)+1$ holds for every $s\in S$, we speak of a {\bf small box}.  
For example, the $(0,1)$-box is small.

\paragraph{Acknowledgement} 
 The authors are grateful to the six authors of
the paper by Borradaile et al.~\cite{BIMOWZ} because that work
triggered the present research (and this is so even if we realized
later that there had been several related works).  
We thank S. Fujishige and S. Iwata for discussion about the history of convex
minimization over base-polyhedra.  We also thank A. J\"uttner 
and T. Maehara for illuminating the essence of the Newton--Dinkelbach algorithm.  
J. Tapolcai kindly draw our attention to engineering
applications in resource allocation.  
Z. Kir\'aly played a similar role by finding
an article which pointed to a work of Levin and Onn on decreasingly
minimal optimization in matroid theory.  
We are also grateful to M. Kov\'acs for drawing our attention to some
important papers in the literature concerning fair resource allocation problems.  
Special thanks are due to T. Migler for her continuous availability to answer our questions
concerning the paper \cite{BIMOWZ} and the work by Borradaile, Migler, and Wilfong \cite{BMW}, 
which paper was also a prime driving force in our investigations.
This research was supported through the program ``Research in Pairs''
by the Mathematisches Forschungsinstitut Oberwolfach in 2019.
The two weeks we could spend at Oberwolfach provided 
an exceptional opportunity to conduct particularly intensive research.




\section{Base-polyhedra} 
\label{egal2base}

Let $S$ be a finite non-empty ground-set and let $b$ be a (fully)
submodular integer-valued set-function on $S$ 
for which $b(\emptyset )=0$ and $b(S)$ is finite.  
A (possibly unbounded) {\bf base-polyhedron} $B=B(b)$ 
in ${\bf R}\sp{S}$ is defined by 
\[
  B=\{x\in {\bf R}\sp{S}:  \widetilde x(S)=b(S), \ \widetilde x(Z)\leq b(Z)
  \ \mbox{ for every } \  Z\subset S\}.
\]

A special base-polyhedron is the one of matroids.  
Given a matroid $M$, Edmonds proved that the polytope 
(that is, the convex hull) of the incidence 
(or characteristic) vectors of the bases of $M$ 
is the base-polyhedron $B(r)$ defined by the rank function $r$ of $M$, 
that is, 
$B(r)=\{x\in {\bf R}\sp{S}:  \ \widetilde x(S)=r(S)$ 
and
$\widetilde x(Z)\leq r(Z)$ for every subset $Z\subset S\}$. 
 (Note that there is no need to require explicitly the non-negativity of $x$,
since this follows from the monotonicity of the rank function $r$:
$x(s)=\widetilde x(S) - \widetilde x(S-s) \geq r(S) - r(S-s)\geq 0$).
It can be proved that a kind of converse also holds, namely, every
(integral) base-polyhedron in the unit $(0,1)$-cube is a matroid base-polyhedron.

For a weak pair $(p,b)$, the polyhedron 
$Q=Q(p,b):= \{x:  p(Z) \leq \widetilde x(Z) \leq b(Z)$
for every $Z\subseteq S\}$ is called a
generalized polymatroid (g-polymatroid, for short).  
By convention,
the empty set is also considered a g-polymatroid.  
In the special case, when $p\equiv 0$ and $b$ is monotone non-decreasing, 
we are back at the concept of polymatroids, introduced by Edmonds \cite{Edmonds70}.  
G-polymatroids were introduced by Frank \cite{FrankP6} 
who proved that $Q(p,b)$ is a non-empty integral polyhedron, 
a g-polymatroid uniquely determines its defining strong pair, 
and the intersection of two integral g-polymatroids $Q(p_{1},b_{1})$
and $Q(p_{2},b_{2})$ is an integral polyhedron which is non-empty 
if and only if 
$p_{1}\leq b_{2}$ and $p_{2}\leq b_{1}$.

The book of Frank \cite{Frank-book} includes an overview of basic
properties and constructions of base-polyhedra and g-polymatroids.
For example, the following operations on a g-polymatroid result in a
g-polymatroid:  projection along axes, translation (or shifting) by a
vector, negation (that is, reflection through the origin),
intersection with a box 
$T(f,g):=\{ x:  f\leq x\leq g \}$, 
intersection with a plank 
$\{x :  \alpha \leq \widetilde x(S)\leq \beta \}$, 
and taking a face.

A base-polyhedron is a special g-polymatroid (where $p(S)=b(S)$) and
every g-polymatroid arises from a base-polyhedron by projecting it
along a single axis.  Each of the operations of translation,
intersection with a box, negation, taking a face, when applied to a
base-polyhedron, results in a base-polyhedron.  
The intersection of a g-polymatroid with a hyperplane 
$\{x:  \widetilde x(S)=\gamma \}$ 
is a base-polyhedron.  
Each of the operations above, when applied to an
integral g-polymatroid (base-polyhedron) results in an integral
g-polymatroid (base-polyhedron).  
The (Minkowski) sum of g-polymatroids (base-polyhedra) 
is a g-polymatroid (base-polyhedron).

We call the translation of a matroid base-polyhedron 
a {\bf translated matroid base-polyhedron}.  
It follows that the intersection of a
base-polyhedron with a small box is a translated matroid base-polyhedron.

A base-polyhedron $B(b)$ is never empty, 
and $B(b)$ is known to be an integral polyhedron.  
(A rational polyhedron is {\bf integral} 
if each of its faces contains an integral element.  
In particular, a pointed rational polyhedron is integral 
if all of its vertices are integral.)
By convention, the empty set is also considered a base-polyhedron.
Note that a real-valued submodular function $b$ also defines a
base-polyhedron $B(b)$ but in the present work we are interested only
in integer-valued submodular functions and integral base-polyhedra.

We call the set $\odotZ{B}$ 
of integral elements of an integral base-polyhedron $B$ an {\bf M-convex set}.  
Originally, this basic notion of Discrete Convex Analysis (DCA), 
introduced by Murota \cite{Murota98a} (see, also the book \cite{Murota03}), 
was defined as a set of integral points in 
${\bf R}\sp{S}$ satisfying certain exchange axioms, and it is known 
that the two properties are equivalent (\cite{Murota03}, Theorem 4.15).  
The set of integral elements of a translated matroid
base-polyhedron will be called a {\bf matroidal M-convex set}.

Since in the present work the central notion is that of
base-polyhedra, we define M-convex sets via base-polyhedra.  
The set of integral points of an integral g-polymatroid is called by 
Murota \cite{Murota03} an M$\sp{\natural}$-convex set 
(pronounce M-natural convex).  
Since a base-polyhedron is a special g-polymatroid, an
M-convex set is a special M$\sp{\natural}$-convex set.  
Note that the original definition of M$\sp{\natural}$-convex sets is different 
(and it is a theorem that the two definitions are equivalent).

A non-empty base-polyhedron $B$ can also be defined by a supermodular
function $p$ for which $p(\emptyset )=0$ and $p(S)$ is finite as follows:  
$B=B'(p)=\{x \in {\bf R}\sp{S}:  
     \widetilde x(S)=p(S),  
     \widetilde x(Z)\geq p(Z) \ \mbox{ for every } \  Z\subset S\}$.

For a set $Z\subset S$, $p\vert Z$ \ ($= p-(S-Z)$) denotes the
restriction of $p$ to $Z$ 
(that is, $p\vert Z$ is obtained from $p$ by deleting $S-Z$), 
while $p'=p/Z$ \ ($= p \div (S-Z)$) is the
set-function on $S-Z$ obtained from $p$ by contracting $Z$, which is
defined for $X\subseteq S-Z$ by $p'(X)= p(X\cup Z)-p(Z)$.  
Note that $p/Z$ and $\ol p\vert (S-Z)$ are complementary set-functions.  
It is also known for disjoint subsets $Z_{1}$ and $Z_{2}$ of $S$ that
\begin{equation}
 (p/Z_{1})/Z_{2} = p/(Z_{1}\cup Z_{2}), \label{(Z1Z2)} 
\end{equation}
that is, contracting first $Z_{1}$ and then $Z_{2}$ is the same
as contracting $Z_{1}\cup Z_{2}$.  
(Indeed, this follows from \
$p_{1}'(X \cup Z_{2})-p_{1}'(Z_{2}) 
= p(X \cup Z_{1} \cup Z_{2}) - p(Z_{1}) - [p(Z_{1} \cup Z_{2})- p(Z_{1})]
= p(X \cup Z_{1} \cup Z_{2})-p(Z_{1} \cup Z_{2})$.) \ 
Furthermore, $(p/Z_{1})\vert Z_{2}= (p\vert (Z_{1} \cup Z_{2}))/Z_{1}$.

It is known that $B$ uniquely determines both $p$ and $b$, namely,
$b(Z) = \max\{\widetilde x(Z):  x\in B\}$ 
and 
$p(Z) = \min\{\widetilde x(Z):  x\in B\}$.  
The functions $p$ and $b$ are complementary functions, 
that is, $b(X)= p(S) - p(S-X)$ or $p(X)= b(S) - b(S-X)$
(where $b(S)=p(S)$).

Let $\{S_{1},\dots ,S_{q}\}$ be a partition of $S$ and let $p_{i}$ be a
supermodular function on $S_{i}$.  
Let $p$ denote the supermodular function on $S$ defined by 
$p(X):  = \sum [p_{i}(S_{i}\cap X):  i=1,\dots ,q]$ for $X\subseteq S$.  
The base-polyhedron $B'(p)$ is called the
{\bf direct sum} of the $q$ base-polyhedra $B'(p_{i})$.  
Obviously, a vector $x\in {\bf R}\sp{S}$ is in $B'(p)$ 
if and only if 
each $x_{i}$ is in $B'(p_{i})$ \ $(i=1,\dots ,q)$,
where $x_{i}$ denotes the restriction
$x\vert S_{i}$ of $x$ to $S_{i}$.

Let $Z$ be a subset of $S$ for which $p(Z)$ is finite.  
The {\bf restriction} of a base-polyhedron $B'(p)$ to $Z$ is the
base-polyhedron $B'(p\vert Z)$.

It is known that a face $F$ of a non-empty base-polyhedron $B$ is also
a base-polyhedron.  
The (special) face of $B'(p)$ defined by the single equality 
$\widetilde x(Z)=p(Z)$ is the direct sum of the base polyhedra 
$B'(p\vert Z)$ and $B'(p/Z)$.  
More generally, any face $F$ of $B$ 
can be described with the help of a chain 
$(\emptyset \subset ) \ C_{1} \subset C_{2} \subset \cdots \subset C_{\ell} = S$ 
of subsets by 
$F:= \{z:  z\in B, \ p(C_{i})=\widetilde z(C_{i})$ 
for $i=1,\dots ,\ell\}$.
(In particular, when $\ell=1$, the face $F$ is $B$ itself.)  
Let
$S_{1}:=C_{1}$ and $S_{i}:=C_{i}-C_{i-1}$ 
for $i=2,\dots ,\ell$.  
Then $F$ is the direct sum of the base-polyhedra $B'(p_{i})$,
where $p_{i}$ is a supermodular function on $S_{i}$ defined by 
$p_{i}(X):= p(X\cup C_{i-1})-p(C_{i-1})$ 
for $X\subseteq S_{i}$.  
In other words, $p_{i}$ is a set-function on $S_{i}$ 
obtained from $p$ by deleting $C_{i-1}$ and contracting $S-C_{i}$.  
The unique supermodular function $p_{F}$ 
defining the face $F$ is given by 
$\sum [p_{i}(S_{i}\cap X) :  i=1,\dots ,\ell]$.
The polymatroid greedy algorithm of Edmonds \cite{Edmonds70} along
with the proof of its correctness, when adapted to base-polyhedra,
shows that $F$ is the set of elements $x$ of $B$ minimizing $cx$
whenever $c:S\rightarrow {\bf R} $ is a linear cost function such that
$c(s) = c(t)$ if $s,t\in S_{i}$ for some $i$ and $c(s) > c(t)$ 
if $s\in S_{i}$ and $t\in S_{j}$ for some subscripts $i<j$.

The intersection of an integral base-polyhedron 
$B=B'(p)$ \ ($=B(\ol p)$) 
and an integral box $T(f,g)$ is an integral base-polyhedron.  
The intersection is non-empty 
if and only if
\begin{equation} 
p\leq \widetilde g \quad \hbox{and} \quad \widetilde f \leq \ol p .
\label{(pgfp)} 
\end{equation}
Note that $g$ occurs only in the first inequality and $f$ occurs only
in the second inequality, from which it follows that if $B$ has an
element $x_{1}\geq f$ and $B$ has an element $x_{2}\leq g$, then $B$ has
an element $x$ with $f\leq x\leq g$.  
This phenomenon is often called the linking principle or linking property.

For an element $m$ of a base-polyhedron $B=B(b)$ defined by a (fully)
submodular function $b$, we call a subset $X\subseteq S$ 
{\bf $m$-tight} (with respect to $b$) if $\widetilde m(X)=b(X)$.  
Clearly, the empty set and $S$ are $m$-tight, and $m$-tight sets are closed
under taking union and intersection.  
Therefore, for each subset $Z\subseteq S$, 
there is a unique smallest $m$-tight set $T_{m}(Z;b)$
including $Z$.  When $Z=\{s\}$ is a singleton, we simply write
$T_{m}(s;b)$ to denote the smallest $m$-tight set containing $s$.  
When the submodular function $b$ in this notation is unambiguous from the
context, we abbreviate $T_{m}(Z;b)$ to $T_{m}(Z)$.

Analogously, when $B=B'(p)$ is given by a supermodular function $p$,
we call $X\subseteq S$ {\bf $m$-tight} (with respect to $p$) if
$\widetilde m(X)=p(X)$.  
In this case, we also use the analogous notation 
$T_{m}(Z)=T_{m}(Z;p)$ and $T_{m}(s)=T_{m}(s;p)$.  
Observe that for complementary functions $b$ and $p$, 
$X$ is $m$-tight with respect to
$b$ precisely if $S-X$ is $m$-tight with respect to $p$.

In applications it is important that weaker set-functions may also
define base-polyhedra.  
For example, if $p$ is an integer-valued crossing supermodular function, 
then $B'(p)$ is still an integral base-polyhedron, 
which may, however, be empty.  This result was proved
independently in \cite{FrankJ8} and in \cite{Fujishige84e}. 
To prove theorems on base-polyhedra, it is much easier to work with
base-polyhedra defined by fully sub- or supermodular functions.  
On the other hand, in applications, base-polyhedra are often defined with
a crossing sub- or supermodular (or even weaker) function.  
For example, the in-degree vectors of the $k$-edge-connected orientations
of a $2k$-edge-connected graph are the integral elements of a
base-polyhedron defined by a crossing supermodular function, 
as was pointed out in an even more general setting \cite{FrankJ4}.  
It is exactly the combination of these double features
that makes it possible to prove a conjecture of
Borradaile et al.~\cite{BIMOWZ} even in an extended form.  
Details will be discussed in Sections \ref{SCori2} and \ref{SCori3}.




\section{Decreasingly minimal elements of M-convex sets}

\subsection{Decreasing minimality}

For a vector $x$, let $x {\downarrow}$
denote the vector obtained from $x$ by rearranging
its components in a decreasing order.  For example,
We call two vectors 
$x$ and $y$ (of same dimension) 
{\bf value-equivalent} if 
$x{\downarrow}= y{\downarrow}$.  
For example, $(2,5,5,1,4)$ and
$(1,4,5,2,5)$ are value-equivalent while the vectors $(3,5,5,3,4)$ \
and \ $(3,4,5,4,4)$ \ are not.

A vector $x$ is {\bf decreasingly smaller} than vector $y$, in
notation $x <_{\rm dec} y$ \ if $x{\downarrow}$ \ is lexicographically
smaller than \ $y{\downarrow}$ in the sense that they are not
value-equivalent and $x{\downarrow}(j) < y{\downarrow}(j)$ for the
smallest subscript $j$ for which $x{\downarrow}(j)$ and
$y{\downarrow}(j)$ differ.  
For example, $x = (2,5,5,1,4)$ is
decreasingly smaller than $y =(1,5,5,5,1)$ \ since 
$x{\downarrow}= (5,5,4,2,1)$ 
is \ lexicographically smaller than \
$y{\downarrow}=(5,5,5,1,1)$.
We write $x\leq_{\rm dec} y$ to mean that $x$ is decreasingly smaller than or
value-equivalent to $y$.

For a set \ $Q$ \ of vectors,
$x\in Q$ is {\bf globally decreasingly minimal} or simply 
{\bf decreasingly minimal} ({\bf dec-min}, for short) 
if \ $x \leq_{\rm dec} y$ \ for every \ $y\in Q$.  
Note that the dec-min elements of $Q$ are value-equivalent.  
Therefore an element $m$ of $Q$ is dec-min if
its largest component is as small as possible, within this, its second
largest component (with the same or smaller value than the largest
one) is as small as possible, and so on.  
An element $x$ of $Q$ is
said to be a {\bf max-minimized} element (a {\bf max-minimizer}, for short) 
if its largest component is as small as possible.  
A max-minimizer element $x$ is {\bf pre-decreasingly minimal}
 ({\bf pre-dec-min}, for short) in $Q$ 
if the number of its largest components is as small as possible.  
Obviously, a dec-min element is pre-dec-min, 
and a pre-dec-min element is max-minimized.

In an analogous way, for a vector $x$, we let $x {\uparrow}$
denote the vector obtained from $x$ by rearranging its
components in an increasing order.  
A vector $y$ is {\bf increasingly larger} than vector $x$, 
in notation $y >_{\rm inc} x$, 
if they are not value-equivalent and 
$y{\uparrow}(j)> x{\uparrow}(j)$ 
holds for the smallest subscript $j$ for which 
$y{\uparrow}(j)$ and $x{\uparrow}(j)$ differ.  
We write $y \geq_{\rm inc} x$ 
if either $y >_{\rm inc} x$ or $x$ and $y$ are value-equivalent.  
Furthermore, we call an element $m$ of $Q$ {\bf (globally) increasingly maximal}
 ({\bf inc-max} for short) 
if its smallest component is as large as possible
over the elements of $Q$, within this its second smallest component is
as large as possible, and so on.  
Similarly, we can use the analogous terms {\bf min-maximized} 
and {\bf pre-increasingly maximal} ({\bf pre-inc-max}).

It should be emphasized that a dec-min element of a base-polyhedron
$B$ is not necessarily integer-valued.  
For example, if
$B=\{(x_{1},x_{2}):  \ x_{1}+x_{2}= 1\}$, 
then $x\sp{*}=(1/2,1/2)$ is a dec-min element of $B$.  
In this case, the dec-min members of $\odotZ{B}$ are
$(0,1) $ and $(1,0)$.

Therefore, finding a dec-min element of $B$ and finding a dec-min
element of $\odotZ{B}$ (the set of integral points of $B$) are two
distinct problems, and we shall concentrate only on the second,
discrete problem.  
In what follows, the slightly sloppy term integral
dec-min element of $B$ will always mean a dec-min element of $\odotZ{B}$.  
(The term is sloppy in the sense that an integral dec-min
element of $B$ is not necessarily a dec-min element of $B$).

We call an integral vector $x\in {\bf Z}\sp S$ {\bf uniform} 
if all of its components are the same integer $\ell$, 
and {\bf near-uniform} if
its largest and smallest components differ by at most 1, that is, if
$x(s)\in \{\ell, \ell +1\}$ for some integer $\ell$ for every $s\in S$.  
Note that if $Q$ consists of integral vectors and the
component-sum is the same for each member of $Q$, then any
near-uniform integral member of $Q$ is obviously both decreasingly
minimal and increasingly maximal integral vector.

\subsection{Characterizing dec-min elements} 
\label{SCchardecmin}

Let $B=B(b)=B'(p)$ be a base-polyhedron defined by an integer-valued
submodular function $b$ or a supermodular function $p$ 
(where $b$ and $p$ are complementary set-functions).  
Let $m$ be an integral member of $B$, that is, $m \in \odotZ{B}$.  A set
$X\subseteq S$ is $m$-tight with respect to $b$ precisely if its
complement $S-X$ is $m$-tight with respect to $p$.  
Recall that $T_{m}(s;b)$ denoted the unique smallest $m$-tight set 
(with respect to $b$) containing $s$.  
In other words, $T_{m}(s;b)$ is the intersection
of all $m$-tight sets containing $s$.  
The easy equivalences in the next claim will be used throughout.

\begin{claim}   \label{st}
Let $s$ and $t$ be elements of $S$ and let $m':=m+\chi_{s}-\chi_{t}$.  
The following properties are pairwise equivalent.  
\smallskip

\noindent 
{\rm (A)} \ $m'\in \odotZ{B}$.  
\smallskip

\noindent 
{\rm (P1)} \ There is no $t\ol s$-set which is $m$-tight
with respect to $p$. 
\smallskip

\noindent 
{\rm (P2)} \ $s\in T_{m}(t;p)$.  
\smallskip

\noindent 
{\rm (B1)} \ There is no $s\ol t$-set which is $m$-tight
with respect to $b$.  
\smallskip

\noindent 
{\rm (B2)} \ $t\in T_{m}(s;b)$.  
\finbox 
\end{claim}

A {\bf 1-tightening step} for $m\in \odotZ{B}$ 
is an operation that replaces $m$ by  
$m':=m+\chi_{s}-\chi_{t}$
where $s$ and $t$ are elements of $S$ for which $m(t)\geq m(s)+2$ and
$m'$ belongs to $\odotZ{B}$.  
(See Claim \ref{st} for properties equivalent to $m'\in \odotZ{B}$.)  
Note that $m'$ is both decreasingly
smaller and increasingly larger than $m$.

Since the mean of the components of $m$ does not change 
at a 1-tightening step while the square-sum of the components of $m$
strictly drops, consecutive 1-tightening steps may occur 
only a finite number of times (even if $B$ is unbounded).

A member $m$ of $\odotZ{B}$ 
is {\bf locally decreasingly minimal} in $\odotZ{B}$ 
if there are no two elements $s$ and $t$ of $S$ such that
$m':=m + \chi_{s} - \chi_{t}$ 
is an element of $\odotZ{B}$
and $m'$ is decreasingly smaller than $m$.  
Note that $m'$ is decreasingly smaller than $m$ precisely if
$m(t) \geq  m(s)+2$.
Obviouly, $m$ is locally decreasingly minimal
if and only if there is no 1-tightening step for $m$.
Note that in this case, $m'$ is also increasingly larger than $m$.  
Analogously, $m$ is {\bf locally increasingly maximal}
if there are no two elements $s$ and $t$ of $S$ such that
$m':=m + \chi_{s} - \chi_{t}$ is an element of $\odotZ{B}$ 
and $m'$ is decreasingly larger than $m$.

The equivalence of the properties in the next claim is immediate from the definitions.

\begin{claim}   \label{local} 
For an integral element $m$ of the integral base-polyhedron $B=B(b)=B'(p)$, 
the following conditions are pairwise equivalent.  
\smallskip

\noindent
{\rm (A1)} \ There is no 1-tightening step for $m$.
\smallskip

\noindent 
{\rm (A2)} \ $m$ is locally decreasingly minimal.  
\smallskip

\noindent 
{\rm (A3)} \ $m$ is locally increasingly maximal.  
\smallskip

\noindent 
{\rm (P1)} \ $m(s) \geq m(t)-1$ holds whenever $t\in S$ and
$s\in T_{m}(t;p)$.  
\smallskip

\noindent {\rm (P2)} \ 
Whenever $m(t)\geq m(s)+2$, there is a $t\ol s$-set $X$
which is $m$-tight with respect to $p$.  
\smallskip

\noindent 
{\rm (B1)} \ $m(s) \geq m(t)-1$ holds whenever $s\in S$ and $t\in T_{m}(s;b)$.
\smallskip

\noindent 
{\rm (B2)} \ 
Whenever $m(t)\geq m(s)+2$, there is an $s\ol t$-set $Y$
which is $m$-tight with respect to $b$.  
\finbox 
\end{claim}

For a given vector $m$ in ${\bf R}\sp S$, we call a set
$X\subseteq S$ an {\bf $m$-top set} (or a top-set with respect to $m$)
if $m(u)\geq m(v)$ holds whenever $u\in X$ and $v\in S-X$.  
Both the empty set and the ground-set $S$ are $m$-top sets, 
and $m$-top sets are closed under taking union and intersection.  
If $m(u)>m(v)$ holds
whenever $u\in X$ and $v\in S-X$, we speak of a {\bf strict $m$-top set}.  
Note that the number of strict non-empty $m$-top sets is at
most $n$ for every $m\in \odotZ{B}$ while $m \equiv 0$ exemplifies that
even all of the non-empty subsets of $S$ can be $m$-top sets.

\begin{theorem}  \label{equi.1} 
Let $b$ be an integer-valued submodular function 
and let $p:= \ol b$ be its complementary  (supermodular) function.  
For an integral element $m$ of the integral base-polyhedron $B=B(b)=B'(p)$,
the following four conditions are pairwise equivalent.
\smallskip

\noindent
{\rm (A)} \ There is no 1-tightening step for $m$
(or anyone of the six other equivalent properties holds in 
Claim {\rm \ref{local}}).

\smallskip

\noindent 
{\rm (B)} \ There is a chain $\cal C$ of $m$-top sets 
$(\emptyset \subset ) \ C_{1}\subset C_{2}\subset \cdots \subset C_\ell = S$ 
which are $m$-tight with respect to $p$ 
(or equivalently, whose complements are $m$-tight with respect to $b$) 
such that the restriction $m_{i}=m\vert S_{i}$ of $m$ to $S_{i}$ is
near-uniform for each member $S_{i}$ of the $S$-partition 
$\{S_{1},\dots,S_\ell\}$, where $S_{1}=C_{1}$ 
and $S_{i}:=C_{i}-C_{i-1}$ $(i=2,\dots ,\ell)$.
\smallskip

\noindent 
{\rm (C1)} \ $m$ is (globally) decreasingly minimal in $\odotZ{B}$.
\smallskip

\noindent 
{\rm (C2)} \ $m$ is (globally) increasingly maximal in $\odotZ{B}$.  
\end{theorem}

\Proof \ 
(B)$\rightarrow $(A): \ If $m(t)\geq m(s)+2$, then there is an
$m$-tight set $C_{i}$ containing $t$ and not containing $s$, 
from which Property (A) follows from Claim \ref{local}.

\smallskip

(A)$\rightarrow $(B): \ 
 Let $\cal C$ be a
longest chain consisting of non-empty $m$-tight and $m$-top sets
$C_{1}\subset C_{2}\subset \cdots \subset C_\ell =S$.  
For notational convenience, let 
$C_{0}=\emptyset $ (but $C_{0}$ is not a member of of $\cal C$).  
We claim that $\cal C$ meets the requirement of (B).  
If, indirectly, this is not the case, then there is a subscript 
$i\in \{1,\dots ,\ell\}$ 
for which $m$ is not near-uniform within
$S_{i}:=C_{i}-C_{i-1}$.  
This means that the max $m$-value $\beta_{i}$ in
$S_{i}$ is at least 2 larger than the min $m$-value 
$\alpha_{i}$ in $S_{i}$, that is, 
$\beta_{i}\geq \alpha_{i} +2$.  
Let 
$Z:= \cup [ T_{m}(t;p): t\in S_{i}, \ m(t)=\beta_{i}]$.  
Then $Z$ is $m$-tight.  
Since $C_{i}$ is $m$-tight, $T_{m}(t;p)\subseteq C_{i}$ holds for $t\in S_{i}$ 
and hence $Z\subseteq C_{i}$.  
Furthermore, (A) implies that 
$m(v) \geq \beta_{i} -1$ for every $v\in Z\cap S_{i}$.

Consider the set $C':=C_{i-1}\cup Z$.  
Then $C'$ is $m$-tight, and
$C_{i-1}\subset C'\subset C_{i}$.  
Moreover, we claim that $C'$ is an $m$-top set.  
Indeed, if, indirectly, there is an element $u\in C'$ and an element $v\in S-C'$ 
for which $m(u)<m(v)$, 
then $u\in Z\cap S_{i}$ and 
$v\in C_{i}-Z$ since both $C_{i-1}$ and $C_{i}$ are $m$-top sets.  
But this is impossible since the $m$-value of each element of
$Z$ is $\beta_{i}$ or $\beta_{i} -1$ while the $m$-value of each element
of $C_{i}-Z$ is at most $\beta_{i} -1$.

The existence of $C'$ contradicts the assumption that $\cal C$
was a longest chain of $m$-tight and $m$-top sets, and
therefore $m$ must be near-uniform within each $S_{i}$, that is, 
${\cal C}$ meets indeed the requirements in (B).

\smallskip

(C1)$\rightarrow $(A) and (C2)$\rightarrow $(A):  \ Property (A) must
indeed hold since a 1-tightening step for $m$ results in an element
$m'$ of $\odotZ{B}$ which is both decreasingly smaller and increasingly
larger than $m$.

\smallskip

(B)$\rightarrow $(C1): \ We may assume that the elements of $S$ are
arranged in an $m$-decreasing order 
$s_{1},\dots ,s_n$ 
(that is,
$m(s_{1})\geq m(s_{2})\geq \cdots \geq m(s_n)$) 
in such a way that each
$C_{i}$ in (B) is a starting segment.  
Let $m'$ be an element of
$\odotZ{B}$ which is decreasingly smaller than or value-equivalent to $m$.  
Recall that $m\vert X$ denoted the vector $m$ restricted to a subset $X\subseteq S$.

\begin{lemma}
For each $i=0,1,\dots ,\ell$, vector $m'\vert C_{i}$ is
value-equivalent to vector $m\vert C_{i}$.  
\end{lemma}  

\Proof Induction on $i$.  
For $i=0$, the statement is void so we assume that $1\leq i\leq \ell$.  
By induction, we may assume that the
statement holds for $j\leq i-1$ and we want to prove it for $i$.
Since $m'\vert C_{i-1}$ is value-equivalent to $m\vert C_{i-1}$ and
$C_{i-1}$ is $m$-tight, it follows that $C_{i-1}$ is $m'$-tight, too.

Let $\beta_{i}$ denote the max $m$-value of the elements of
$S_{i}=C_{i}-C_{i-1}$.  By the hypothesis in (B), 
the maximum and the minimum of the $m$-values in $S_{i}$ differ by at most 1. 
Hence we can assume that there are $r_{i}>0$ elements in $S_{i}$ 
with $m$-value $\beta_{i}$ and 
$\vert S_{i}\vert -r_{i}\geq 0$ elements 
with $m$-value $\beta_{i} -1$.

As $m\vert C_{i-1}$ is value-equivalent to $m'\vert C_{i-1}$ and $m'$
was assumed to be decreasingly smaller than or value-equivalent to
$m$, we can conclude that $m'\vert (S-C_{i-1})$ is decreasingly smaller 
than or value-equivalent to $m\vert (S-C_{i-1})$.  
Therefore, $S_{i}$ contains at most $r_{i}$ elements of $m'$-value $\beta_{i}$ and hence
\begin{align*}
p(C_{i}) &\leq \widetilde m'(C_{i}) = \widetilde m'(C_{i-1}) + \widetilde m'(S_{i}) 
\\ &
\leq \widetilde m'(C_{i-1}) + r_{i}\beta_{i} + (\vert S_{i}\vert -r_{i}) (\beta_{i} -1) 
\\ & 
=  \widetilde m(C_{i-1}) + r_{i}\beta_{i}
    + (\vert S_{i}\vert -r_{i}) (\beta_{i} -1) = p(C_{i}),
\end{align*}
from which equality follows everywhere.  
In particular, $S_{i}$ contains exactly $r_{i}$
elements of $m'$-value $\beta_{i}$ and $\vert S_{i}\vert -r_{i}$ elements
of $m'$-value $\beta_{i} -1$, proving the lemma.  
\finbox 
\medskip

By the lemma, $m'$ is value-equivalent to $m$, and hence $m$ is a
decreasingly minimal element of $\odotZ{B}$, that is, (C1) follows.

\smallskip

(B)$\rightarrow $(C2): \ The property in (C1) that $m$ is globally
decreasing minimal in $\odotZ{B}$ is equivalent to the statement that
$-m$ is globally increasing maximal in $-{\odotZ{B}}$, 
that is, (C2) holds with respect to $-m$ and $-{\odotZ{B}}$.  
As we have already proved the implications 
(C2)$\rightarrow$(A)$\rightarrow$(B)$\rightarrow $(C1), 
it follows that (C1) holds for $-m$ and $-{\odotZ{B}}$.  
But (C1) for $-m$ and $-{\odotZ{B}}$ is just the same
as (C2) for $m$ and $\odotZ{B}$.  
\finbox \finboxHere
\medskip

\subsection{Minimizing the sum of the $k$ largest components}

A decreasingly minimal element of $\odotZ{B}$ has the starting property
that its largest component is as small as possible.  
As a natural extension, one may be interested in finding a member of $\odotZ{B}$ 
for which the sum of the $k$ largest components is as small as possible.
We refer to this problem as 
{\bf min $k$-largest-sum}.

\begin{theorem}\label{klargest} 
Let $B$ be an integral base-polyhedron and
$k$ an integer with $1\leq k\leq n$.  Then any dec-min element $m$ of
$\odotZ{B}$ is a solution to Problem min $k$-largest-sum.  
\end{theorem}

\Proof 
Observe first that if $z_{1}$ and $z_{2}$ are dec-min elements of $\odotZ{B}$, 
then it follows from the very definition of decreasing
minimality that the sum of the first $j$ largest components of $z_{1}$
and of $z_{2}$ are the same for each $j=1,\dots ,n$.  
Let $K$ denote the
sum of the first $k$ largest components of any dec-min element, and
assume indirectly that there is a member $y\in \odotZ{B}$ for which the
sum of its first largest components is smaller than $K$.  
Assume that the componentwise square-sum of $y$ is as small as possible.  
By the previous observation, $y$ is not a dec-min element.  
Theorem \ref{equi.1} implies that there are elements $s$ and $t$ of $S$ for
which $y(t)\geq y(s)+2$ and $y':=y - \chi_t + \chi_s$ is in $\odotZ{B}$.  
The sum of the first $k$ largest components of $y'$ is at
most the sum of the first $k$ largest components of $y$, and hence
this sum is also smaller than $K$.  
But this contradicts the choice of
$y$ since the componentwise square-sum of $y'$ is 
strictly smaller than that of $y$.  
\finbox 
\medskip

The theorem implies that a dec-min element $m$ is 
a solution to the min $k$-largest-sum  problem for each $k=1,\dots ,n$.  
In \cite{Frank-Murota.2}, we will point out 
that this property characterizes dec-min elements
and is closely related to
the notion of least majorization investigated in 
\cite{AS18}, \cite{MOA11}, \cite{Tamir95}.

\subsection{An example for the intersection of two base-polyhedra}
\label{egal2.example}

We proved for a single base-polyhedron that an integral element is decreasingly minimal 
if and only if 
it is increasingly maximal
 (and therefore the two properties could jointly be called egalitarian).
The following example shows that the two properties may differ if the
polyhedron is the intersection of two (integral) base-polyhedra.

Let  $S=\{ s_{1},s_{2},s_3,s_4 \}$ 
be the common ground-set of two rank-2
matroids $M_{1}$ and $M_{2}$ which are described by their circuits.  
Both matroids have two 2-element circuits.  
Namely, the circuits of $M_{1}$
are $\{s_{1},s_4\}$ and $\{s_{2},s_3\}$,
while the circuits of $M_{2}$ are $\{s_{1},s_3\}$ and $\{s_{2},s_4\}$.  
Both matroids have four bases:
\begin{eqnarray*} 
{\cal B}_{1} &:=& \{ \{s_{1}, s_{2}\}, \ \{s_3,s_4 \}, \ \{s_{1},s_3\}, \{s_{2},s_4\}\}, 
\\ 
{\cal B}_{2} &:=& \{ \{s_{1}, s_{2}\}, \ \{s_3,s_4\}, \ \{s_{1},s_4\}, \{s_{2},s_3\}\}.  
\end{eqnarray*}
Let $B_{1}$ and $B_{2}$, respectively, denote the base-polyhedra of $M_{1}$ and $M_{2}$.  
(That is, $B_{i}$ is the convex hull of the incidence vectors of 
the four members of ${\cal B}_{i}$.)  The common bases are as follows:
\[
 {\cal B}:= {\cal B}_{1}\cap {\cal B}_{2} = \{ \{s_{1}, s_{2}\}, \  \{s_3,s_4\}\}.
\]
Let $B$ denote the convex hull of the incidence vectors of the two
members of ${\cal B}$.  
That is, $B$ is the line segment connecting
$(1,1,0,0)$ and $(0,0,1,1)$.  
Note that $B=B_{1}\cap B_{2}$ since 
the intersection of two matroid base-polyhedra is an integral polyhedron,
by a theorem of Edmonds \cite{Edmonds70}.
Let $B'_{i}$ denote the base-polyhedron obtained from $B_{i}$ by adding
the vector $(1, -1, 0, 0)$ to $B_{i}$.  Let $B'$ be obtained from $B$ in the same way.  
Then we have:
\begin{eqnarray*} 
\odotZ{B_{1}'} &=& \{(2, 0, 0, 0), \ (1, -1, 1, 1), \ (2, -1, 1, 0), \ (1, 0, 0, 1)\}, 
\\ 
\odotZ{B'_{2}} &=& \{(2, 0, 0, 0), (1, -1, 1, 1), (2, \-1, 0, 1), (1, 0, 1, 0)\}, 
\\ 
\odotZ{B} &=& \{(2, 0, 0, 0), \ (1, -1, 1, 1)\}.  
\end{eqnarray*}
Now $x = (2, 0, 0, 0)$ is an increasingly maximal element of
$\odotZ{B}$ while $y = (1, -1, 1, 1)$ is decreasingly minimal.

Therefore finding a decreasingly minimal integral element and an
increasingly maximal integral element of the intersection $B$ of two
base-polyhedra are two different problems (unlike the analogous
problems for a single base-polyhedron).  
The two problems, however, are equivalent in the sense that 
an element $x$ of $\odotZ{B}$ is decreasingly minimal 
if and only if 
$-x$ is an increasingly maximal
element of $-{\odotZ{B}}$, and if $B$ is the intersection of two
base-polyhedra, then so is $-B$.

Furthermore, we claim that the intersection $B$ has no integral 
least majorized element
(see \cite{AS18}, \cite{MOA11} as well as
Part II \cite{Frank-Murota.2} for the definition).
Indeed, we have only the two possible choices 
$x = (2, 0, 0, 0)$ and $y = (1, -1, 1, 1)$, 
but $x$ is not least majorized
since the largest component of $y$ is smaller than the largest
component of $x$, and $y$ is neither a least majorized element of $\odotZ{B}$ 
since the sum of the 3 largest components of $x$ is smaller
than the sum of the 3 largest components of $y$.

\medskip

In Part IV \cite{Frank-Murota.4}, we shall consider the
even more general problem of finding an increasingly minimal
(integral) submodular flow.  
Note that that submodular flows not only
generalize ordinary flows and circulations but the intersection of two
$g$-polymatroids is also a special submodular flow polyhedron.




\section{Characterizing the set of pre-decreasingly minimal elements}
\label{peak}

We continue to assume that $p$ is an integer-valued (with possible
$-\infty $ values but with finite $p(S)$) 
supermodular function, which implies that $B=B'(p)$ is a non-empty
integral base-polyhedron.  
We have already proved that an integral
element $m$ of $B$ (that is, an element of $\odotZ{B}$) is
decreasingly-minimal ($=$ dec-min) precisely if $m$ is
increasingly-maximal ($=$ inc-max).

One of our main goals is to prove that the set ${\rm dm}(\odotZ{B})$ of all
dec-min elements of $\odotZ{B}$ is an M-convex set, meaning that
there exists an integral base-polyhedron 
$B\sp{\bullet} \subseteq B$ 
such that ${\rm dm}(\odotZ{B})$ 
is the set of integral elements of $B\sp{\bullet}$. 
In addition, we shall show that ${\rm dm}(\odotZ{B})$ is actually a matroidal
M-convex set, that is, $B\sp{\bullet}$ is a special base-polyhedron
which is obtained from a matroid base-polyhedron by translating it
with an integral vector.

The base-polyhedron $B\sp{\bullet}$ will be obtained with the help of a
decomposition of $B$ along a certain \lq canonical\rq\ partition
$\{S_{1},S_{2},\dots ,S_{q}\}$ of $S$ into non-empty sets.  
To this end, we
start by introducing the first member $S_{1}$ of this partition along
with a matroid on $S_{1}$.  
The set $S_{1}$, depending only on $B$, will
be called the peak-set of $S$.

\subsection{Max-minimizers and pre-dec-min elements}
\label{SCmaxminzerpredm}

Recall that an element of $\odotZ{B}$ was called a max-minimizer if 
its largest component was as small as possible, 
while a max-minimizer was called a pre-dec-min element of $\odotZ{B}$ 
if the number of its
maximum components was as small as possible.  
As a dec-min element of
$\odotZ{B}$ is automatically pre-dec-min (in particular, a max-minimizer), 
we start our investigations by studying max-minimizers
and pre-dec-min elements of $\odotZ{B}$.  
For a number 
    $\beta$, we say that a vector is {\bf $\beta$-covered} 
if each of its components is at most $\beta$.
Throughout our discussions, 
\begin{equation} 
\beta_{1}:=\beta(B)
\label{(beta1-def)} 
\end{equation} 
denotes the smallest integer for which
$\odotZ{B}$ has a $\beta_{1}$-covered  element.  
In other words, $\beta_{1}$ is the largest component of a max-minimizer of $\odotZ{B}$.  
Therefore $\beta_{1}$ is the largest component of any 
pre-dec-min (and hence any dec-min) element of $\odotZ{B}$.  
Note that
an element $m$ of $\odotZ{B}$ is $\beta_{1}$-covered 
precisely if $m$ is a max-minimizer.

\begin{theorem}   \label{betamin} 
For the largest component $\beta_{1}$ of a
max-minimizer of $\odotZ{B}$, one has 
\begin{equation} 
\beta_{1} =\max \{ \llceil {p(X)  \over \vert X\vert }\rrceil : 
    \emptyset \not =X\subseteq S\}.
\label{(betamin)} 
\end{equation} 
\end{theorem}

\Proof 
It follows from formula \eqref{(pgfp)} that $B$ has a $\beta$-covered element 
if and only if 
\begin{equation} 
\beta \vert X\vert \geq p(X) \
\hbox{whenever}\ \ X\subseteq S. \label{(beta-covered)} 
\end{equation} 
Moreover,
if $\beta$ is an integer and \eqref{(beta-covered)} holds, then $B$
has an integral $\beta$-covered element.  
As $\beta \vert X\vert \geq
p(X)$ holds for an arbitrary $\beta$ when $X=\emptyset $, it follows
that the smallest integer $\beta$ meeting this \eqref{(beta-covered)}
is indeed 
$\max \{ \llceil {p(X) \over \vert X\vert }\rrceil : \emptyset \not =X\subseteq S\}$.  
\finbox 
\medskip

For a $\beta_{1}$-covered element $m$ of $\odotZ{B}$, let $r_{1}(m)$
denote the number of $\beta_{1}$-valued components of $m$.  Recall that
for an element $s\in S$ we denoted the unique smallest $m$-tight set
containing $s$ by $T_{m}(s)=T_{m}(s;p)$ (that is, $T_{m}(s)$ is the
intersection of all $m$-tight sets containing $s$).  
Furthermore, let
\begin{equation}  \label{(S1m.def)} 
S_{1}(m):=\cup \{T_{m}(t):  m(t)=\beta_{1}\}.  
\end{equation}
Then $S_{1}(m)$ is $m$-tight and $S_{1}(m)$ is actually the unique
smallest $m$-tight set containing all the $\beta_{1}$-valued elements of $m$.

\begin{theorem}   \label{optkrit} 
A $\beta_{1}$-covered element $m$ of $\odotZ{B}$ is pre-dec-min 
if and only if 
$m(s)\geq \beta_{1} -1$ for each $s\in S_{1}(m)$.  
\end{theorem}

\Proof 
Necessity.  
Let $m$ be a pre-dec-min element of $\odotZ{B}$.
For any $\beta_{1}$-valued element $t\in S$ and any element $s\in
T_{m}(t)$, we claim that $m(s)\geq \beta_{1} -1$.  
Indeed, if we had $m(s)\leq \beta_{1} -2$, then the vector $m'$ arising from $m$ by
decreasing $m(t)$ by 1 and increasing $m(s)$ by 1 belongs to $B$
(since $T_{m}(t)$ is the smallest $m$-tight set containing $t$) 
and has one less $\beta_{1}$-valued components than $m$ has, contradicting the
assumption that $m$ is pre-dec-min.

Sufficiency.  
Let $m'$ be an arbitrary $\beta_{1}$-covered integral
element of $B$.  Abbreviate $S_{1}(m)$ by $Z$ and let $h'$ denote the
number of elements $z\in Z$ for which $m'(z)=\beta_{1}$.  
Then
\begin{eqnarray*} 
& &\vert Z\vert (\beta_{1} -1) + r_{1}(m) 
= \widetilde m(Z) = p(Z) \leq \widetilde m'(Z) 
\\ & & 
\ \leq h'\beta_{1} + (\vert Z\vert -h')(\beta_{1} -1) 
       = \vert Z\vert (\beta_{1} -1) + h' 
\\ & &
\ \leq  \vert Z\vert (\beta_{1} -1) + r_{1}(m') ,
\end{eqnarray*}
from which $r_{1}(m)\leq r_{1}(m')$, as required.  
\finbox
\medskip

Define the set-function $h_{1}$ on $S$ as follows.  
\begin{equation} 
h_{1}(X):= p(X) - (\beta_{1} -1)\vert X\vert \ \ \hbox{for }\ \ X\subseteq S.
\label{(h1.def)} 
\end{equation}

\begin{theorem}   \label{betavalmin} 
For the minimum number $r_{1}$ of $\beta_{1}$-valued components 
of a $\beta_{1}$-covered member of $\odotZ{B}$,  one has
\begin{equation} 
r_{1}= \max \{h_{1}(X):  X\subseteq S\}.  
\label{(r1)} 
\end{equation} 
\end{theorem}

\Proof 
Let $m$ be an element of $\odotZ{B}$ for which the maximum of
its components is $\beta_{1}$, and let $X$ be an arbitrary subset of $S$.  
Suppose that $X$ has $\ell$ \ $\beta_{1}$-valued components.  
Then
\begin{equation} 
p(X) \ \leq  \ \widetilde m(X) 
   \ \leq \  \ell \beta_{1} + (\vert X\vert -\ell)(\beta_{1} -1) 
  \ = \  \vert X\vert (\beta_{1} -1) + \ell 
  \ \leq \  \vert X\vert (\beta_{1} -1) + r_{1}(m) ,
\label{(estim)} 
\end{equation}
from which $r_{1}(m) \geq p(X) - (\beta_{1} -1)\vert X\vert =h_{1}(X)$, 
implying that 
\[
 r_{1} = \min \{r_{1}(m) :  m\in \odotZ{B},   
    \mbox{ \ $m$ is  $\beta_{1}$-covered } \} 
  \ \geq \  \max \{h_{1}(X): X\subseteq S\}. 
\]

In order to prove the reverse inequality, we have to find 
a $\beta_{1}$-covered integral element $m$ of $B$ and a subset $X$ of $S$ for
which $r_{1}(m) = h_{1}(X)$, 
which is equivalent to requiring that each of the three
inequalities in \eqref{(estim)} holds with equality.  
That is, the following three {\bf optimality criteria} hold:  
\ (a) \ $X$ is $m$-tight, 
\ (b) \ $X$ contains all $\beta_{1}$-valued components of $m$, and 
\ (c) \ $m(s)\geq \beta_{1} -1$ for each $s\in X$.

Let $m$ be a pre-dec-min element of $B$.  
Then $S_{1}(m)$ is $m$-tight, $S_{1}(m)$ contains all $\beta_{1}$-valued elements and, 
by Theorem \ref{optkrit}, $m(s)\geq \beta_{1} -1$ for all $s\in S_{1}(m)$, 
therefore $m$ and $S_{1}(m)$ satisfy the three optimality criteria.  
\finbox

\medskip

Note that $r_{1}$ is the number of $\beta_{1}$-valued components of any
pre-dec-min element (and in particular, any dec-min element) of $\odotZ{B}$.

\subsection{The peak-set $S_{1}$}

Since the set-function $h_{1}$ introduced in \eqref{(h1.def)} is supermodular, 
the maximizers of $h_{1}$ are closed under taking intersection and union.  
Let $S_{1}$ denote the unique smallest subset of $S$ maximizing $h_{1}$.  
In other words, $S_{1}$ is the intersection of all sets maximizing $h_{1}$.  
We call this set $S_{1}$ the {\bf peak-set} of $B$ (and of $\odotZ{B}$).

\begin{theorem}  \label{smallest} 
For every pre-dec-min 
(and in particular, for every dec-min) 
element $m$ of $\odotZ{B}$, the set $S_{1}(m)$
introduced in \eqref{(S1m.def)} is independent of the choice of $m$
and $S_{1}(m)=S_{1}$, where $S_{1}$ is the peak-set of $B$.  
\end{theorem}

\Proof 
It follows from Theorem \ref{betavalmin} that, given a
pre-dec-min element $m$ of $B$, a subset $X$ is maximizing $h_{1}$ 
precisely if the three optimality criteria mentioned in the proof hold.  
Since $S_{1}(m)$ meets the optimality criteria, 
it follows that $S_{1}\subseteq S_{1}(m)$.  
If, indirectly, there is an element $s\in S_{1}(m) -S_{1}$, 
then $m(s)=\beta_{1} -1$ since $S_{1}$ contains all the
$\beta_{1}$-valued elements.  
By the definition of $S_{1}(m)$, there is a $\beta_{1}$-valued element 
$t\in S_{1}(m)$ for which the smallest $m$-tight set $T_{m}(t)$ contains $s$, 
but this is impossible since $S_{1}$ is an $m$-tight set containing $t$ but not $s$.  
\finbox

\medskip

Since $S_{1}=S_{1}(m)$ is $m$-tight and near-uniform, we obtain that 
\[
\beta_{1}
= \llceil {\widetilde m_{1} (S_{1}) \over \vert S_{1}\vert }\rrceil 
= \llceil {p(S_{1}) \over \vert S_{1}\vert }\rrceil , 
\]
and the definitions of $S_{1}$ and $r_{1}$ imply that 
\begin{equation} 
  r_{1}=p(S_{1}) - (\beta_{1}-1)\vert S_{1}\vert . 
\label{(r11)} 
\end{equation}

\begin{proposition}   \label{PR=S1} 
$S_{1}=\{s\in S:  \mbox{\rm  there is a pre-dec-min element } \
       m \in \odotZ{B} \ \mbox{\rm with } \  m(s)=\beta_{1}\}$.  
For every pre-dec-min element $m$ of $\odotZ{B}$,
$m(s)\geq \beta_{1} -1$ for every $s\in S_{1}$ and  
$m(s)\leq \beta_{1} -1$ for every $s\in S-S_{1}$.
\end{proposition}

\Proof 
If $m(s)=\beta_{1}$ for some pre-dec-min $m$, then $s\in S_{1}(m)=S_{1}$.  
Conversely, let $s\in S_{1}$ and let $m$ be a pre-dec-min element.  
We are done if $m(s)=\beta_{1}$.  
If this is not the case,
then $m(s)=\beta_{1} -1$ by Theorem \ref{optkrit}.  
By the definition of $S_{1}(m)$, 
there is an element $t\in S_{1}(m)$ for which $m(t)=\beta_{1}$ and $s\in T_{m}(t)$.  
But then $m':= m+\chi_{s}-\chi_{t}$
 is in $\odotZ{B}$, $m'(s)=\beta_{1}$ and $m'$ is also
pre-dec-min as it is value-equivalent to $m$.  
\finbox

\subsection{Separating along $S_{1}$}
\label{SCsepS1}

Let $S_{1}$ be the peak-set occurring in Theorem \ref{smallest} and 
let $S_{1}':=S-S_{1}$.  
Let $p_{1}=p\vert S_{1}$ denote the restriction of $p$ to
$S_{1}$, and let $B_{1}\subseteq {\bf R}\sp {S_{1}}$ denote the
base-polyhedron defined by $p_{1}$, that is, $B_{1}:=B'(p_{1})$.  
Suppose that $S_{1}'\not =\emptyset $ and let $p_{1}':=p/S_{1}$, 
that is, $p_{1}'$ is the set-function on $S_{1}'$ 
obtained from $p$ by contracting $S_{1}$ \
($p_{1}'(X)= p(S_{1}\cup X)-p(S_{1})$ for $X\subseteq S_{1}'$).

Consider the face $F$ of $B$ determined by $S_{1}$, that is, 
$F$ is the direct sum of the base-polyhedra $B_{1}=B'(p_{1})$ and $B_{1}'=B'(p_{1}')$.
Then the dec-min elements of $\odotZ{B}_{1}$ are exactly the integral
elements of the intersection of $B_{1}$ and the box given by 
$\{x: \beta_{1} -1\leq x(s)\leq \beta_{1} \ \mbox{ for every } s\}$.  
 Hence the dec-min elements of
$\odotZ{B}_{1}$ are near-uniform.

\begin{theorem}   \label{felbont} 
An integral vector $m=(m_{1},m_{1}')$ is a dec-min element of $\odotZ{B}$ 
if and only if 
$m_{1}$ is a dec-min element of $\odotZ{B}_{1}$ 
and $m_{1}'$ is a dec-min element of $\odotZ{B_{1}'}$.  
\end{theorem}

\Proof 
Suppose first that $m$ is a dec-min element of $\odotZ{B}$.
Then $S_{1}=S_{1}(m)$ by Theorem \ref{smallest} and $m$ is a max-minimizer,
implying that every component of $m$ in $S_{1}$ is of value $\beta_{1} -1$
or value $\beta_{1}$, and $m$ has exactly $r_{1}$ components of value $\beta_{1}$.  
Therefore each of the components of $m_{1}$ is $\beta_{1} -1$ or $\beta_{1}$, 
that is, $m_{1}$ is near-uniform.
Since $m_{1}$ is obviously in $\odotZ{B}_{1}$,
$m_{1}$ is indeed dec-min in $\odotZ{B}_{1}$.

Since $\widetilde m(S_{1}) = p(S_{1})$, for a set $X\subseteq S_{1}'$, we have 
\[
\widetilde m_{1}'(X) = \widetilde m(X) 
  = \widetilde m(S_{1}\cup X)- \widetilde m(S_{1}) 
  = \widetilde m(S_{1}\cup X) - p(S_{1}) 
  \geq p(S_{1}\cup X) - p(S_{1}) = p_{1}'(X).
\]
Furthermore 
\[
\widetilde m_{1}'(S_{1}') 
= \widetilde m(S_{1}') 
= \widetilde m(S_{1}\cup S_{1}') - \widetilde m(S_{1}) 
= p(S_{1}\cup S_{1}') - p(S_{1}) = p_{1}'(S_{1}'),
\] 
that is, $m_{1}'$ is in $\odotZ{B_{1}'}$.  If,
indirectly, $m_{1}'$ is not dec-min, then, by applying Theorem
\ref{equi.1} to $S_{1}'$, $m_{1}'$, and $p_{1}'$, we obtain that there are
elements $t$ and $s$ of $S_{1}'$ for which $m_{1}'(t) \geq m_{1}'(s)+2$ and
\ $(*)$ \ no $t\ol s$-set exists which is $m_{1}'$-tight with respect to $p_{1}'$.  
On the other hand, $m$ is a dec-min element of $\odotZ{B}$ for which 
\[
m(t)=m_{1}'(t) \geq m_{1}'(s)+2 = m(s)+2,
\] 
and hence there must
be a $t\ol s$-set $Y$ which is $m$-tight with respect to $p$.

Since $S_{1}$ is $m$-tight with respect to $p$, 
the set $S_{1}\cup Y$ is also $m$-tight with respect to $p$.  
Let $X:=S_{1}'\cap Y$.  
Then
\[
\widetilde m(X) + \widetilde m(S_{1}) 
  = \widetilde m(S_{1}\cup Y) 
  = p(S_{1}\cup Y) = p (S_{1}\cup X),
\]
and hence
\[
\widetilde m_{1}'(X)  = \widetilde m(X) 
= p(S_{1}\cup X) - \widetilde m(S_{1}) = p(S_{1}\cup X) - p(S_{1}) = p_{1}'(X),
\]
that is, $X$ is a $t\ol s$-set which is $m_{1}'$-tight with respect to $p_{1}'$, 
in contradiction with statement $(*)$ above that no such set exists.

To see the converse, assume that $m_{1}$ is a dec-min element of
$\odotZ{B}_{1}$ and $m_{1}'$ is a dec-min element of $\odotZ{B_{1}'}$.  
This immediately implies that $m$ is in the face $F$ of $B$ determined by $S_{1}$.  
Suppose, indirectly, that $m$ is not a dec-min element of $\odotZ{B}$.  
By Theorem \ref{equi.1}, there are elements $t$ and $s$
of $S$ for which $m(t) \geq m(s)+2$ and 
$(**)$ no $t\ol s$-set exists which is $m$-tight with respect to $p$.  
If $t\in S_{1}$, then $s$
cannot be in $S_{1}$ since the $m$-value of each element of $S_{1}$ is
$\beta_{1}$ or $\beta_{1} -1$.  But $S_{1}$ is $m_{1}$-tight with respect to
$p$ and hence it is $m$-tight with respect to $p$, contradicting property $(**)$.  
Therefore $t$ must be in $S_{1}'$, 
implying, by Proposition \ref{PR=S1}, that $s$ is also in $S_{1}'$.

Since $m_{1}'$ is a dec-min element of $\odotZ{B_{1}'}$, there must be a
$t\ol s$-set $Y\subset S_{1}'$ which is $m_{1}'$-tight with respect to $p_{1}'$.  
It follows that 
\[
\widetilde m(Y) = \widetilde m_{1}'(Y) 
= p_{1}'(Y) = p(S_{1}\cup Y) - p(S_{1}) 
\leq \widetilde m(S_{1}\cup Y) - \widetilde m(S_{1})=\widetilde m(Y),
\]
from which $\widetilde m(S_{1}\cup Y) = p(S_{1}\cup Y)$, 
contradicting property $(**)$ that no $t\ol s$-set exists which is $m$-tight with respect to $p$.  
\finbox 
\medskip

An important consequence of Theorem \ref{felbont} is that, in order to
find a dec-min element of $\odotZ{B}$, it will suffice to find
separately a dec-min element of $\odotZ{B_{1}}$ (which was shown above to
be a near-uniform vector) and a dec-min element of $\odotZ{B_{1}'}$.  
The algorithmic details will be discussed in Section \ref{SCalgo2}.

\begin{theorem}  \label{m1.ekv} 
Let $S_{1}$ be the peak-set of $\odotZ{B}$. 
For an element $m_{1}$ of $\odotZ{B_{1}}$, 
the following properties are pairwise equivalent.
\smallskip

\noindent 
{\rm (A1)} \ 
$m_{1}$ has $r_{1} \ \ (= p(S_{1}) - (\beta_{1}-1)\vert S_{1}\vert $ \ $>0)$\ \ 
components of value $\beta_{1}$ 
and 
$\vert S_{1}\vert -r_{1}$ \ $(\geq 0)$ \ 
components of value $\beta_{1} -1$.

\smallskip

\noindent 
{\rm (A2)} \ $m_{1}$ is near-uniform. 

\smallskip

\noindent 
{\rm (A3)} \ $m_{1}$ is dec-min in $\odotZ{B_{1}}$.  
\smallskip

\noindent 
{\rm (B1)} \ $m_{1}$ is the restriction of a dec-min element
$m$ of $\odotZ{B}$ to $S_{1}$.  
\smallskip

\noindent 
{\rm (B2)} \ $m_{1}$ is the restriction of a pre-dec-min
element $m$ of $\odotZ{B}$ to $S_{1}$.  
\end{theorem}

\Proof 
The implications (A1)$\rightarrow $(A2)$\rightarrow $(A3) and
(B1)$\rightarrow $(B2) are immediate from the definitions.

(A3)$\rightarrow $(B1):  
Let $m_{1}'$ be an arbitrary dec-min element of $\odotZ{B_{1}'}$.  
By Theorem \ref{felbont}, $m:=(m_{1},m_{1}')$ is a dec-min
element of $\odotZ{B}$ and hence $m_{1}$ is indeed the restriction of a
dec-min element of $\odotZ{B}$ to $S_{1}$.

(B2)$\rightarrow $(A1):
By Theorems \ref{optkrit} and \ref{smallest},
we have $m_{1}(s)\geq \beta_{1} -1$ for each $s\in S_{1}(m)=S_{1}$, that is,
$\beta_{1} -1\leq m_{1}(s)\leq \beta_{1}$.  
By letting $r'$ denote the number of $\beta_{1}$-valued components of $m_{1}$, 
we obtain by \eqref{(r11)} that 
\[
r_{1}+(\beta_{1} -1)\vert S_{1}\vert = p_{1}(S_{1})
    = \widetilde m_{1}(S_{1})  = (\beta_{1} -1)\vert S_{1}\vert +r'
\] 
and hence $r'=r_{1}$.  
\finbox

\medskip 
Theorem \ref{felbont} implies that, in order to characterize
the set of dec-min elements of $\odotZ{B}$, it suffices to characterize
the set of dec-min elements of $\odotZ{B_{1}'}$.

\begin{theorem}   \label{beta2} 
Let $\beta_{2}$ denote the smallest integer for
which $\odottZ{B_{1}'}$ has a $\beta_{2}$-covered element, that is, 
$\beta_{2}= \beta (B_{1}')$.  
Then 
\begin{equation} 
\beta_{2} = \max \{ \llceil {p_{1}'(X) \over \vert X\vert }\rrceil :  
                  \emptyset \not =X\subseteq S-S_{1}\},
\label{(betamin2)} 
\end{equation}
 where $p_{1}'(X)= p(X\cup S_{1}) - p(S_{1})$.
Furthermore, $\beta_{2}$ is the largest component in $S-S_{1}$ of every
dec-min element of $\odotZ{B}$, and $\beta_{2} < \beta_{1}$.  
\end{theorem}  

\Proof 
Formula \eqref{(betamin2)} follows by applying Theorem \ref{betamin} 
to base-polyhedron $B_{1}'$ ($=B'(p_{1}')$) in place of $B$.
By Theorem \ref{felbont}, the largest component in $S-S_{1}$ of any
dec-min element $m$ of $\odotZ{B}$ is $\beta_{2}$.  
By Theorem \ref{smallest}, $S_{1}(m)=S_{1}$, and the definition of $S_{1}(m)$ shows
that $m(s)\leq \beta_{1} -1$ holds for every $s\in S-S_{1}$, from which
$\beta_{2}<\beta_{1}$ follows.  
\finbox

\subsection{The matroid $M_{1}$ on $S_{1}$} \label{M1}

It is known from the theory of base-polyhedra that the intersection of
an integral base-polyhedron with an integral box is an integral base-polyhedron.  
Moreover, if the box in question is small, then the
intersection is actually a translated matroid base-polyhedron 
(meaning that the intersection arises from a matroid base-polyhedron by
translating it with an integral vector).  
This result is a consequence of the theorem that 
($*$) any integral base-polyhedron in the unit
$(0,1)$-cube is the convex hull of (incidence vectors of) the bases of a matroid.

Consider the special small integral box $T_{1}\subseteq {\bf Z}\sp{S_{1}}$ defined by 
\[
T_{1}:=\{x:  \beta_{1} -1\leq x(s)\leq \beta_{1}\}
\]
and its intersection $B_{1}\sp{\bullet} := B_{1}\cap T_{1}$ 
with the base-polyhedron
$B_{1}$ investigated above.  
Therefore $B_{1}\sp{\bullet}$ is a translated matroid
base-polyhedron and Theorem \ref{m1.ekv} implies the following.

\begin{corollary}   \label{trans.mat} 
The dec-min elements of $\odotZ{B_{1}}$ are exactly 
the integral elements of the translated matroid base-polyhedron $B_{1}\sp{\bullet}$.  
\finbox 
\end{corollary}

Our next goal is to reprove Corollary \ref{trans.mat} by concretely
describing the matroid in question 
and not relying on the background theorem ($*$) mentioned above.  
For a dec-min element $m_{1}$ of $\odotZ{B_{1}}$, let 
\[
   L_{1}(m_{1}):= \{s\in S_{1}: m_{1}(s)=\beta_{1}\}.
\] 
We know from Theorem \ref{m1.ekv} that 
$\vert L_{1}(m_{1})\vert =r_{1}$.  
Define a set-system ${\cal B}_{1}$ as follows:
\begin{equation} 
{\cal B}_{1}:=\{ L\subseteq S_{1}:  L=L_{1}(m_{1})
 \ \hbox{ for some dec-min element $m_{1}$ of $\odotZ{B_{1}}\}.$}\ 
\label{(B1def)} 
\end{equation}

We need the following characterization of ${\cal B}_{1}$.

\begin{proposition} \label{PRbaseM1}
 An $r_{1}$-element subset $L$ of $S_{1}$ is in ${\cal B}_{1}$
if and only if 
\begin{equation} 
\vert L\cap X\vert \geq p_{1}'(X):= p_{1}(X) - (\beta_{1} -1)\vert X\vert  
\ \hbox{ whenever } \  X\subseteq S_{1}.  
\label{(p01)} 
\end{equation} 
\end{proposition}

\Proof 
Suppose first that $L\in {\cal B}_{1}$, that is, 
there is a dec-min element $m_{1}$ of $\odotZ{B_{1}}$ for which $L=L_{1}(m_{1})$.  
Then
\[
(\beta_{1} -1)\vert X\vert + \vert X\cap L\vert = \widetilde m_{1}(X) \geq p_{1}(X),
\] 
for every subset $X\subseteq S_{1}$ from which \eqref{(p01)} follows.

To see the converse, let $L\subseteq S_{1}$ be an $r_{1}$-element set meeting \eqref{(p01)}.  
Let
\begin{equation} 
m_{1}(s):= \begin{cases} 
          \beta_{1} & \ \ \hbox{if}\ \ \ s\in L \cr
         \beta_{1} -1 & \ \ \hbox{if}\ \ \ s\in S-L.  
         \end{cases} 
\end{equation}
Then obviously $L=L_{1}(m_{1})$.  Furthermore,
\[
\widetilde m_{1}(S_{1}) 
  = (\beta_{1} -1)\vert S_{1}\vert + \vert L\vert 
  = (\beta_{1} -1)\vert S_{1}\vert + r_{1} =p(S_{1})
\] 
and
\[
  \widetilde m_{1}(X) 
   = (\beta_{1} -1)\vert X\vert + \vert L\cap X\vert
  \geq p_{1}(X) \ \hbox{whenever}\ \ X\subset S_{1},
\] 
showing that $m_{1}\in B_{1}$.  
Since $m_{1}\in T_{1}$, we conclude that $m_{1}$ is a dec-min element of $\odotZ{B_{1}}$.  
\finbox

\begin{theorem}   \label{matroid1} 
The set-system \ ${\cal B}_{1}$ defined in \eqref{(B1def)} 
forms the set of bases of a matroid $M_{1}$ on ground-set $S_{1}$.  
\end{theorem}

\Proof 
The set-system ${\cal B}_{1}$ is clearly non-empty and all of its
members are of cardinality $r_{1}$.
It is widely known \cite{Edmonds70} 
that for an integral submodular function $b$ on a ground-set $S_{1}$ the set-system 
\[
\{ \  L\subseteq S_{1}:  \ \vert L\cap X\vert \leq b(X) \ \hbox{ whenever } \ 
 X\subset S_{1}, \ \vert L\vert = b(S_{1}) \ \},
\] 
if non-empty,
satisfies the matroid basis axioms.  This implies for the supermodular
function $p_{1}'$ that the set-system 
$\{L:  \vert L\cap X\vert \geq p_{1}'(X) \ \hbox{ whenever } \ 
     X\subset S_{1}, \ \vert L\vert = p_{1}'(S_{1}) \}$, 
if non-empty, forms the set of bases of a matroid.  
By applying this fact to the supermodular function $p_{1}'$ defined by 
$p_{1}'(X):= p_{1}(X) - (\beta_{1} -1)\vert X\vert $, 
one obtains that ${\cal B}_{1}$ is non-empty
and forms the set of bases of a matroid.  
\finbox

\medskip

In this way, we proved the following more explicit form of Corollary \ref{trans.mat}.

\begin{corollary}   \label{M1B1} 
Let $\Delta_{1}:S_{1}\rightarrow {\bf Z} $ denote
the integral vector defined by $\Delta_{1}(s):=\beta_{1} -1$  for $s\in S_{1}$.  
A member $m_{1}$ of $\odotZ{B_{1}}$ is decreasingly minimal 
if and only if 
there is a basis $B_{1}$ of $M_{1}$ such that $m_{1}=\chi_{B_{1}} + \Delta_{1}$.
\finbox 
\end{corollary}

\subsection{Value-fixed elements of $S_{1}$}
\label{value-fixed}

We say that an element $s\in S$ is {\bf value-fixed} with respect to $\odotZ{B}$ 
if $m(s)$ is the same for every dec-min element $m$ of $\odotZ{B}$.  
In Section \ref{SCoptdualset}, we will show a description of value-fixed
elements of $\odotZ{B}$.  In the present section, we consider the
value-fixed elements with respect to $B_{1}$, 
that is, $s\in S_{1}$ is value-fixed if 
$m_{1}(s)$ is the same for every dec-min element $m_{1}\in \odotZ{B_{1}}$.  
Recall that $m_{1}\in \odotZ{B_{1}}$ was shown to be dec-min precisely if 
$\beta_{1}-1\leq m_{1}(s)\leq \beta_{1}$ for each $s\in S_{1}$.

A {\bf loop} of a matroid is an element $s\in S_{1}$ not belonging to any basis.  
(Often the singleton $\{s\}$ is called a loop, that is,
$\{s\}$ is a one-element circuit).  
A {\bf co-loop} (or cut-element or isthmus) 
of a matroid is an element $s$ belonging to all bases.

\begin{proposition}  
$M_{1}$ has no loops.  
\end{proposition}  

\Proof
 By Proposition \ref{PR=S1}, for every $s\in S_{1}$ there is a
pre-dec-min element $m$ of $\odotZ{B}$ for which $m(s)=\beta_{1}$.  
Then $m_{1}:=m\vert S_{1}$ is a pre-dec-min element of $\odotZ{B_{1}}$ 
by Theorem \ref{m1.ekv} from which $s_{1}$ belongs to a basis of $M_{1}$ 
by Corollary \ref{M1B1}.  
\finbox

\medskip

The proposition implies that:

\begin{proposition}   \label{fixbetai} 
If $s\in S_{1}$ is value-fixed (with respect to $B_{1}$), 
then $m_{1}(s)=\beta_{1}$ for every dec-min element $m_{1}$ of $\odotZ{B_{1}}$.  
\finbox 
\end{proposition}

By Corollary \ref{M1B1}, an element $s\in S_{1}$ is a co-loop of $M_{1}$
if and only if 
$m_{1}(s)=\beta_{1}$ holds for every dec-min element $m_{1}$ of $\odotZ{B_{1}}$.  
This and Theorem \ref{felbont} imply the following.

\begin{theorem}   \label{co-loop} 
For an element $s\in S_{1}$, the following properties are pairwise equivalent. 
\smallskip

\noindent 
{\rm (A)} \ $s$ is a co-loop of $M_{1}$.  
\smallskip

\noindent 
{\rm (B)} \ $s$ is value-fixed.  
\smallskip

\noindent 
{\rm (C)} \ $m(s)=\beta_{1}$ holds for every dec-min element $m$ of $\odotZ{B}$.  
\finbox 
\end{theorem}

Our next goal is to characterize the set of value-fixed elements of $S_{1}$.  
Consider the family of subsets $S_{1}$ defined by 
\begin{equation} 
{\cal F}_{1}:= \{X\subseteq S_{1}:  \ \beta_{1} \vert X\vert = p_{1}(X)\}.
\label{(scriptF1)} 
\end{equation}
The empty set belongs to ${\cal F}_{1}$ and it is possible that 
${\cal F}_{1}$ has no other members.  
By standard submodularity arguments,
${\cal F}_{1}$ is closed under taking union and intersection.  
Let $F_{1}$ denote the unique largest member of ${\cal F}_{1}$.  
It is possible that $F_{1}=S_{1}$ in which case we call $S_{1}$ {\bf degenerate}.

\begin{theorem}   \label{F1} 
An element $s\in S_{1}$ is value-fixed 
if and only if 
$s\in F_{1}$.  
\end{theorem}  

\Proof 
Let $m_{1}$ be a dec-min member of $\odotZ{B_{1}}$.  Then
\[
 \beta_{1}\vert F_{1}\vert \geq \widetilde m_{1}(F_{1}) \geq p_{1}(F_{1}) =
\beta_{1}\vert F_{1}\vert 
\]
 and hence we must have $\beta_{1}=m_{1}(s)$ for
every $s\in F_{1}$, that is, the elements of $F_{1}$ are indeed value-fixed.

Conversely, let $s$ be value-fixed, that is, $m_{1}(s)=\beta_{1}$ for
each dec-min element $m_{1}$ of $\odotZ{B_{1}}$.  
Let $m_{1}$ be a dec-min member of $\odotZ{B_{1}}$.  
Let $Z$ denote the unique smallest set containing $s$
for which $\widetilde m_{1}(Z)= p_{1}(Z)$.  
(That is, $Z=T_{m_{1}}(s;p_{1}).)$
We claim that $m_{1}(t)=\beta_{1}$ for every element $t\in Z$.  
For if $m_{1}(t)=\beta_{1} -1$ for some $t$, then 
$m'_{1}:=m_{1}- \chi_{s}+\chi_{t}$ 
would also be a dec-min member of $\odotZ{B_{1}}$, 
contradicting the assumption that $s$ is value-fixed.
Therefore $ p_{1}(Z) = \widetilde m_{1}(Z) = \beta_{1}\vert Z\vert $ from
which the definition of $F_{1}$ implies 
that $Z\subseteq F_{1}$ and hence $s\in F_{1}$.  
\finbox



\section{The set of dec-min elements of an M-convex set}
\label{canonical}

Let $B=B'(p)$ denote again an integral base-polyhedron defined by the
(integer-valued) supermodular function $p$.  
As in the previous section, 
$\odotZ{B}$ continues to denote the M-convex set consisting of
the integral vectors (points, elements) of $B$.  
Our present goal is
to provide a complete description of the set of decreasingly-minimal
(= egalitarian) elements of $\odotZ{B}$
by identifying a partition of the ground-set, 
to be named the canonical partition,
inherent in this problem.
As a consequence, we show that the set of
dec-min elements has a matroidal structure and this feature makes it
possible to solve the minimum cost dec-min problem.

\subsection{Canonical partition and canonical chain}

In Section \ref{peak} we introduced the integer $\beta_{1}$ as the
minimum of the largest component of the elements of $\odotZ{B}$ as well
as the notion of peak-set $S_{1}$ of $S$.  
We considered the face of $B$ defined by $S_{1}$ 
that was the direct sum of base-polyhedra
$B_{1}=B'(p_{1})$ and $B_{1}'=B'(p_{1}')$, 
where $p_{1}$ denoted the restriction
of $p$ to $S_{1}$ while $p_{1}'$ 
arose from $p$ by contracting $S_{1}$ 
(that is, $p_{1}'(X)=p(S_{1}\cup X)-p(S_{1})$).

A consequence of Theorem \ref{felbont} is that, in order to
characterize the set of dec-min elements of $\odotZ{B}$, 
it suffices to characterize separately the dec-min elements 
of $\odotZ{B_{1}}$ and the dec-min elements of $\odotZ{B_{1}'}$.  
In Theorem \ref{m1.ekv}, we characterized the dec-min elements of $\odotZ{B_{1}}$ 
as those belonging to the small box 
$T_{1}:=\{x\in {\bf R}\sp{S_{1}}:  \beta_{1} -1\leq x(s)\leq \beta_{1}$
 \ for \ $s\in S_{1}\}$.  
We also proved that the set
$\odotZ{B_{1}\sp{\bullet}}$
of dec-min elements of $\odotZ{B_{1}}$ can be described
with the help of matroid $M_{1}$.  
If the peak-set $S_{1}$ happens to be
the whole ground-set $S$, then the characterization of the set of
dec-min elements of $\odotZ{B}$ is complete.  
If $S_{1}\subset S$, then
our remaining task is to characterize the set of dec-min elements of
$\odotZ{B_{1}'}$.  This can be done by repeating iteratively the
separation procedure to the base-polyhedron 
$B_{1}' = B'(p_{1}') \subseteq {\bf R}\sp{S-S_{1}}$ 
described in Section \ref{peak} for $B$.

In this iterative way, we are going to define a partition 
${\cal P}\sp{*}=\{S_{1},S_{2}, \dots , S_{q}\}$ of $S$ 
which determines a chain 
${\cal C}\sp{*}=\{C_{1},C_{2},\dots ,C_{q}\}$ where 
$C_{i}:=S_{1}\cup S_{2}\cup \cdots \cup S_{i}$ (in particular $C_{q}=S)$, 
and the supermodular function 
\[
p_{i}':=p/C_{i} \quad  \mbox{ on set } \  \ol {C_{i}}:=S-C_{i} 
\]
which defines the base-polyhedron 
$B_{i}'=B'(p_{i}')$ in ${\bf R}\sp{\ol C_{i}}$.  
Moreover, we define iteratively a decreasing sequence 
$\beta_{1}>\beta_{2}>\cdots >\beta_{q}$ 
of integers, a small box 
\begin{equation} 
 T_{i}:=\{x\in {\bf R}\sp{S_{i}}: \ \beta_{i}-1 \leq x(s)\leq \beta_{i} \ \hbox{for}\ s\in S_{i}\},
\label{(boxi)} 
\end{equation} 
and the supermodular function $p_{i}$ on $S_{i}$, where 
\begin{equation} 
p_{i}:=p_{i-1}'\vert S_{i}  \quad (=(p/C_{i-1})\vert S_{i}), 
\label{(pidef)} 
\end{equation} 
that is,
\[
 p_i(X)= p(X\cup C_{i-1}) - p(C_{i-1}) \ \hbox{for}\ \ X\subseteq S_i. 
\]
Let $B_{i} := B'(p_{i})\subseteq {\bf R}\sp{S_{i}}$
be the base-polyhedron defined by $p_i$.

In the general step, suppose that the pairwise disjoint non-empty sets
$S_{1}, S_{2},\dots ,S_{j-1}$ 
have already been defined, along with the decreasing sequence 
$\beta_{1}>\beta_{2}>\cdots >\beta_{j-1}$ 
of integers.  
If $S=S_{1}\cup \cdots \cup S_{j-1}$, then by taking $q:=j-1$, 
the iterative procedure terminates.  
So suppose that this is not the case, that is, $C_{j-1}\subset S$.  
We assume that $p_{j-1}$
on $S_{j-1}$ has been defined as well as $p_{j-1}'$ on $\ol {C_{j-1}}$.

Let 
\begin{equation} 
\beta_{j}=\max \{ \llceil {p_{j-1}'(X) \over \vert X\vert }\rrceil :  
    \emptyset \not =X\subseteq \ol {C_{j-1}}\},
\label{(betaminj} 
\end{equation}
that is,
\begin{equation} 
\beta_{j}=\max \{ \llceil {p(X\cup C_{j-1}) - p(C_{j-1}) \over \vert X\vert }\rrceil :  
    \emptyset \not =X\subseteq \ol {C_{j-1}}\}.
\label{(betaj} 
\end{equation}
Note that, by the iterative feature of these definitions,
Theorem \ref{beta2} implies that 
\[
\beta_{j}<\beta_{j-1}.  
\]
Furthermore, let $h_{j}$ be a set-function on $\ol {C_{j-1}}$ defined as follows:
\begin{equation} 
h_{j}(X):= p_{j-1}'(X) - (\beta_{j}-1)\vert X\vert 
 \ \ \hbox{for }\ \ X\subseteq \ol {C_{j-1}}, 
\label{(hj.def)}
\end{equation}
and let $S_{j}\subseteq \ol {C_{j-1}}$ be the peak-set of $\ol
{C_{j-1}}$ assigned to $B_{j-1}':= B'(p_{j-1}')$, 
that is, $S_{j}$ is the smallest subset of $\ol C_{j-1}$ maximizing $h_{j}$.  
Finally, let $p_{j}:=p_{j-1}'\vert S_{j}$ and let $p_{j}':= p_{j-1}'/S_{j}$.  
Observe by \eqref{(Z1Z2)} that $p_{j}'= p/C_{j}$.  
Therefore $p_{j}$ is a set-function on $S_{j}$ while $p_{j}'$ is defined on $\ol {C_{j}}$.

We shall refer to the partition ${\cal P}\sp{*}$ and the chain 
${\cal C}\sp{*}$ defined above as the {\bf canonical partition} and 
{\bf canonical chain} of $S$, respectively, assigned to $B$, 
while the sequence $\beta_{1}>\cdots >\beta_{q}$ 
will be called the {\bf essential value-sequence} of $\odotZ{B}$.  
Let $B\sp{\oplus}$ denote the face of $B$ 
defined by the canonical chain ${\cal C}\sp{*}$, 
that is, $B\sp{\oplus}$ is the direct sum of the $q$ base-polyhedra 
$B'(p_{i}) \ (i=1,\dots ,q)$.  
Finally, let $T\sp{*}$ be the direct sum of the small boxes $T_{i}$ $(i=1,\dots ,q)$, 
that is, $T\sp{*}$ is the integral box
defined by the essential value-sequence as follows:
\begin{equation} 
T\sp{*}:= \{x\in {\bf R}\sp{S}:  \ \beta_{i}-1\leq x(s)\leq \beta_{i} 
  \ \ \hbox{whenever}\ s\in S_{i} \ (i=1,\dots ,q)\}, 
\label{(box)} 
\end{equation}
and let 
\[
B\sp{\bullet}:= B\sp{\oplus}\cap T\sp{*}.
\]
is always an integral base-polyhedron and hence $B\sp{\bullet}$ 
is an integral base-polyhedron.  
Furthermore, $B\sp{\bullet}$ is the direct sum of the $q$ base-polyhedra 
$B_{i}\cap T_{i}$ ($i=1,\dots ,q$), where $B_{i}=B'(p_{i})$, 
implying that a vector $m$ is in $\odotZ{B\sp{\bullet}}$ 
if and only if each $m_{i}$ is in $\odotZ{B_{i}}\cap T_{i}$, where $m_{i}=m\vert S_{i}$.

\begin{theorem}  \label{main.1} 
Let $B=B'(p)$ be an integral base-polyhedron on ground-set $S$.  
The set of decreasingly-minimal elements of $\odotZ{B}$ 
is (the {\rm M}-convex set) $\odotZ{B\sp{\bullet}}$.
Equivalently, an element $m\in \odotZ{B}$ is decreasingly-minimal 
if and only if 
its restriction $m_{i}:=m\vert S_{i}$ to $S_{i}$ belongs to
$B_{i}\cap T_{i}$ for each $i=1,\dots ,q$, 
where $\{S_{1},\dots ,S_{q}\}$ is the canonical partition of $S$ belonging to $B$, 
\ $T_{i}$ is the small box defined in \eqref{(boxi)}, 
and $B_{i}$ is the base-polyhedron
$B'(p_{i})$ belonging to the supermodular set-function $p_{i}$ defined in \eqref{(pidef)}.  
\end{theorem}

\Proof 
We use induction on $q$.  
Suppose first that $q=1$, that is,  $S_{1}=S$ and $B_{1}=B$.  
If $m$ is a dec-min element of $B$, 
then the equivalence of Properties (A1) and (A3) in Theorem \ref{m1.ekv}
implies that $m$ is in $\odotZ{B\sp{\bullet}}$.  
If, conversely, $m\in \odotZ{B\sp{\bullet}}$, 
then $m$ is near-uniform and, by the equivalence of
Properties (A1) and (A3) in Theorem \ref{m1.ekv} again, $m$ is dec-min.

Suppose now that $q\geq 2$ and consider the base-polyhedron 
$B_{1}' = B'(p_{1}')$ 
appearing in Theorem \ref{felbont}.  
The iterative definition of the canonical partition ${\cal P}\sp{*}$ 
implies that the canonical partition of $S-S_{1}$ assigned to 
$B_{1}'$ is $\{S_{2},\dots ,S_{q}\}$ 
and the essential value-sequence belonging to $B_{1}'$ is 
$\beta_{2}>\beta_3>\cdots >\beta_{q}$.  
Also, the canonical chain 
${\cal C}':= \{C_{2}',\dots ,C_{q}'\}$
of $B_{1}'$ consists of the sets 
$C_{i}'=S_{2}\cup \cdots \cup S_{i} = C_{i}-S_{1}$ \ $(i=2,\dots ,q)$.

By applying the inductive hypothesis to $B_{1}'$, we obtain that an
integral element $m_{1}'$ of $B_{1}'$ is dec-min 
if and only if
$m_{1}'$ is in the face of $B_{1}'$ defined by chain ${\cal C}'$ and $m_{1}'$ 
belongs to the box 
$T':= \{x\in {\bf R}\sp{S-S_{1}}:  \ \beta_{i}-1\leq x(s)\leq
   \beta_{i} \ \ \hbox{whenever}\ s\in S_{i} \ (i=2,\dots ,q)\}$.  
By applying Theorem \ref{felbont}, we are done in this case as well.
\finbox 
\medskip

\begin{corollary} \label{COchardecmin}
Let $B=B'(p)$ be an integral base-polyhedron on ground-set $S$.  
Let $\{C_{1},\dots ,C_{q}\}$ be the canonical chain,
$\{S_{1},\dots ,S_{q}\}$ the canonical partition of $S$,
and $\beta_{1}> \beta_{2}> \dots > \beta_{q}$ 
the essential value-sequence belonging to $\odotZ{B}$.  
Then an element $m\in \odotZ{B}$ is decreasingly-minimal 
if and only if 
each $C_{i}$ is $m$-tight 
(that is, $\widetilde m(C_{i})=p(C_{i})$) 
and $\beta_{i} -1 \leq m(s) \leq \beta_{i}$ holds for each 
$s\in S_{i}$ \ $(i=1,\dots ,q)$.  
\finbox 
\end{corollary}

\subsection{Obtaining the canonical chain and value-sequence from a dec-min element}

The main goal of this section is to show that 
the canonical chain and value-sequence can be rather easily obtained 
from an arbitrary dec-min element of $\odotZ{B}$.
This approach will be crucial in developing a polynomial algorithm 
in Section~\ref{SCcompcanchain}
for computing the essential value-sequence along with  the canonical chain and partition.

Let $m$ be an element of $\odotZ{B}$.  
We called a set $X\subseteq S$ $m$-tight if $\widetilde m(X)=p(X)$.  
Recall from Section \ref{egal2base} that, 
for a subset $Z\subseteq S$, $T_{m}(Z)=T_{m}(Z;p)$
denoted the unique smallest $m$-tight set including $Z$, that is,
$T_{m}(Z)$ is the intersection of all the $m$-tight sets including $Z$.
Obviously,
\begin{equation} 
T_{m}(Z)= \cup (T_{m}(z) :  z\in Z). 
 \label{(Tm)} 
\end{equation}

Let $m$ be an arbitrary dec-min element of $\odotZ{B}$.  
We proved that $m$ is in the face $B\sp{\oplus}$ of $B$ defined 
by the canonical chain 
${\cal C}\sp{*} =\{C_{1},\dots ,C_{q}\}$
belonging to $B$.  
Therefore each $C_{i}$ is $m$-tight with respect to $p$.  
Furthermore $m_{i}:=m\vert S_{i}$ 
belongs to the box $T_{i}$ defined in \eqref{(boxi)}.  
This implies that 
$m(s)\geq \beta_{i}-1$ for every $s\in C_{i}$ and 
$m(s')\leq \beta_{i+1}$ for every $s'\in \ol C_{i}$. 
 (The last inequality holds indeed
since $s'\in \ol C_{i}$ implies that $s'\in S_{j}$ for some $j\geq i+1$
from which $m(s')\leq \beta_{j}\leq \beta_{i+1}$.)  
Since $\beta_{i+1}\leq \beta_{i}-1$, we obtain that each $C_{i}$ is an $m$-top set.

Since $m_{i}$ is near-uniform on $S_{i}$ with values $\beta_{i}$ and
possibly $\beta_{i}-1$, we obtain
\[
 \beta_{i}= \llceil {\widetilde m_{i} (S_{i}) \over \vert S_{i}\vert }\rrceil
= \llceil {p_{i} (S_{i}) \over \vert S_{i}\vert }\rrceil 
= \llceil {p (C_{i})- p(C_{i-1}) \over \vert S_{i}\vert }\rrceil .
\]
Let 
$L_{i}:=\{s\in S-C_{i-1}:  m(s) =\beta_{i}\}$
and let $r_{i}:=\vert L_{i}\vert $. 
Then $p_{i}(S_{i}) = \widetilde m_{i}(S_{i}) = (\beta_{i}-1)\vert
S_{i}\vert + r_{i}$ and hence
\begin{equation} 
r_{i}= p(C_{i}) - p(C_{i-1}) - (\beta_{i}-1)\vert S_{i}\vert .
\label{(ri)} 
\end{equation}

The content of the next lemma is that, once $C_{i-1}$ is given, the
next member $C_{i}$ of the canonical chain (and hence $S_{i}$, as well)
can be expressed with the help of $m$.  
Recall that
$T_{m}(L_{i})=T_{m}(L_{i};p)$ denoted the smallest $m$-tight set including $L_{i}$.

\begin{lemma}   \label{Sim} 
$C_{i}= C_{i-1} \cup T_{m}(L_{i};p)$.  
\end{lemma}

\Proof 
Recall the definition of function $h_{i}$ given in \eqref{(hj.def)}.  
We have 
\begin{equation} 
h_{i}(S_{i}) = r_{i} 
\label{(hiSi)} 
\end{equation}
since 
$h_{i}(S_{i}) 
  =p_{i-1}'(S_{i}) -(\beta_{i}-1)\vert S_{i}\vert 
  = p(S_{i}\cup C_{i-1}) - p(C_{i-1}) - (\beta_{i}-1)\vert S_{i}\vert 
  =\widetilde m(C_{i}) - \widetilde m(C_{i-1}) - (\beta_{i}-1)\vert S_{i}\vert
  = \widetilde m(S_{i}) - (\beta_{i}-1)\vert S_{i}\vert = r_{i}$.

Since $L_{i}\subseteq C_{i}$ and each of $C_{i-1}$, $C_{i}$, and $T_{m}(L_{i})$
are $m$-tight, we have $C_{i-1} \cup T_{m}(L_{i};p)\subseteq C_{i}$.  
For $X':= T_{m}(L_{i}) \cap \ol {C_{i-1}}$ we have
\begin{eqnarray*} h_{i}(X') 
&=& p(C_{i-1}\cup T_{m}(L_{i})) -p(C_{i-1}) - (\beta_{i}-1)\vert X_{i}'\vert 
\\ &=& 
\widetilde m(C_{i-1}\cup T_{m}(L_{i}))
- \widetilde m(C_{i-1}) - (\beta_{i}-1)\vert X_{i}'\vert 
\\ &=&
\widetilde m(X') - (\beta_{i}-1)\vert X_{i}'\vert = \vert L_{i}\vert =r_{i} = h_{i}(S_{i}), 
\end{eqnarray*}
that is, $X'$ is also a maximizer of $h_{i}(X)$.  
Since $S_{i}$ was the smallest maximizer of $h_{i}$, 
we conclude that $C_{i-1} \cup T_{m}(L_{i};p) \supseteq C_{i}$.  
\finbox

\medskip 
The lemma implies that both the essential value-sequence
$\beta_{1}>\cdots >\beta_{q}$ 
and the canonical chain ${\cal C}\sp{*}$
belonging to $\odotZ{B}$ can be directly obtained from $m$.

\begin{corollary}   \label{Cim} 
Let $m$ be an arbitrary dec-min element of $\odotZ{B}$.  
The essential value-sequence and the canonical chain belonging to $\odotZ{B}$ 
can be described as follows.
Value $\beta_{1}$ is the largest $m$-value and $C_{1}$ is the smallest
$m$-tight set containing all $\beta_{1}$-valued elements.  
Moreover, for $i=2,\dots ,q,$ \ $\beta_{i}$ is the largest value of 
$m\vert \ol {C_{i-1}}$ and $C_{i}$ is the smallest $m$-tight set 
(with respect to $p$) 
containing each element of $m$-value at least $\beta_{i}$.
\finbox 
\end{corollary}

A detailed algorithm based on this corollary will be described in
Section \ref{SCcompcanchain}. 
Note that a dec-min element $m$ of $\odotZ{B}$ 
may have more than $q$ distinct values.  
For example, if $q=1$ and $L_{1}\subset C_{1}=S$, then
$m$ has two distinct values, namely $\beta_{1}$ on the elements of
$L_{1}$ and $\beta_{1} -1$ on the elements of $S-L_{1}$,
while its essential value-sequence consists of the single member $\beta_{1}$.

\medskip

\paragraph{A direct proof} 
Corollary \ref{Cim} implies that
the chain of subsets and value-sequence assigned to a dec-min element $m$ of $\odotZ{B}$ 
in the corollary do not depend on the choice of $m$.
Here we describe an alternative, direct proof of this consequence.

\begin{theorem} \label{THnew54} 
Let $m$ be an arbitrary dec-min element of $\odotZ{B}$.  
Let $\beta_{1}$ denote the largest value of $m$ and 
let $C_{1}$ denote the smallest $m$-tight set (with respect to $p$)
containing all $\beta_{1}$-valued elements.  
Moreover, for $i=2,3,\dots,q$, let $\beta_{i}$ denote 
the largest value of $m\vert \overline{C_{i-1}}$
and let $C_{i}$ denote the smallest $m$-tight set 
containing each element of $m$-value at least $\beta_{i}$.  
Then the chain $C_{1}\subset C_{2}\subset \cdots \subset C_{q}$ 
and the sequence 
$\beta_{1}>\beta_{2}>\cdots >\beta_{q}$ do not depend on the choice of $m$.
\end{theorem}

\Proof 
Let $z$ be dec-min element of $\odotZ{B}$.  
We use induction on the number of elements $t$ of $S$ for which $m(t)>z(t)$.  
If no such an element $t$ exists, then $m=z$ and there is nothing to prove.  
So assume that $z\not =m$.

Let $L_{i}:=\{t\in S_{i}:  m(t)=\beta _{i}\}$.  
As $m$ is dec-min, the definition of $C_{i}$ implies that 
$m(s)=\beta_{i}-1$ holds for every element 
$s \in S_{i}-L_{i}$.
Let $t\in L_{i}$ and let $s\in T_{m}(t)-L_{i}$.
Then $m':= m + \chi_{s} - \chi_{t}$ is also a dec-min element of $\odotZ{B}$, 
and we say that $m'$ is obtained from $m$ by an elementary step.  
Observe that $T_{m}(t) = T_{m'}(s)$ and hence the chain and 
the value-sequence assigned to $m'$ is the same as those assigned to $m$.

Let $i$ denote the smallest subscript for which $m\vert S_{i}$ and $z\vert S_{i}$ differ.  
Since $z$ is dec-min, $z(s)\leq \beta _{i}$ holds for every $s\in S_{i}$.  
Let $L'_{i}:=\{t\in S_{i}:  z(t)=\beta _{i}\}$.  
Then
$z(v)\leq \beta _{i}-1$ for every $v\in S_{i}-L'_{i}$, 
and $\vert L'_{i}\vert \leq \vert L_{i}\vert $ as $z$ is dec-min.  
Therefore
\[
\widetilde z(S_{i})
 \ \leq \ 
\beta _{i}\vert L'_{i}\vert + (\beta _{i}-1)(\vert S_{i}-L'_{i}\vert ) 
\ = \ 
(\beta _{i}-1)\vert S_{i}\vert + \vert L'_{i}\vert 
\ \leq \
(\beta _{i}-1)\vert S_{i}\vert + \vert L_{i}\vert .
\]
On the other hand,
\begin{align*}
\widetilde z(S_{i}) & = \widetilde z(C_{i}) -\widetilde z(C_{i-1}) =
\widetilde z(C_{i}) - \widetilde m(C_{i-1}) 
\\
 & \geq p(C_{i}) -
\widetilde m(C_{i-1}) = \widetilde m(C_{i}) - \widetilde m(C_{i-1})=
\widetilde m(S_{i})= (\beta _{i}-1)\vert S_{i}\vert + \vert L_{i}\vert .
\end{align*}
Therefore we have equality throughout, in particular, 
$\widetilde z(C_{i})=p(C_{i})$, 
$\vert L'_{i}\vert = \vert L_{i}\vert $, and $z(v)=\beta_{i}-1$ for every $v\in S_{i}-L'_{i}$.

Let $t\in L_{i}$ be an element for which $m(t)>z(t)$.  
Then $m(t)=\beta_{i}$ and $z(t)=\beta _{i}-1$.  
It follows that $T_{m}(t)$ contains an element $s$ for which $z(s) > m(s)$, 
implying that $m(s)=\beta _{i}-1$ and $z(s)=\beta _{i}$.  
Now $m(t)>m'(t) = z(t)$ holds for the dec-min element 
$m':= m + \chi_{s} - \chi_{t}$ 
obtained from $m$ by an elementary step, and therefore we are done by induction.  
\finbox

\subsection{Matroidal description of the set of dec-min elements}
\label{SCmatrdescdmset}

In Section \ref{M1}, we introduced a matroid $M_{1}$ on $S_{1}$ and proved in 
Corollary \ref{trans.mat}
that the dec-min elements of $\odotZ{B_{1}}$
are exactly the integral elements of the translated base-polyhedron of $M_{1}$, 
where the translation means the addition of the constant vector
$(\beta_{1} -1,\dots ,\beta_{1} -1)$ of dimension $\vert S_{1}\vert $. 
The same notions and results can be applied to each subscript $i=2,\dots,q$.
Furthermore, by formulating Lemma~\ref{matroid1} 
for subscript $i$ in place of 1, we obtain the following.

\begin{proposition} \label{matroidi} 
The set-system \
${\cal B}_{i}:=\{ L\subseteq S_{i}:  L=L_{i}(m_{i})$ 
{\rm for some dec-min element} $m_{i}$ {\rm of} $\odotZ{B_{i}}\}$ \
forms the set of bases of a matroid $M_{i}$ on ground-set $S_{i}$.  
An $r_i$-element subset $L$ of $S_i$ is a basis of $M_i$ if and only if 
\begin{equation} 
\vert L\cap X\vert  \geq    p_{i}'(X)  
  := p_{i}(X)  - (\beta_{i}-1)\vert X\vert 
\label{(p0)} 
\end{equation} 
holds for every $X\subseteq S_{i}$.
\finbox 
\end{proposition}

It follows that a vector $m_{i}$ on $S_{i}$ is a dec-min element of $\odotZ{B_{i}}$ 
if and only if 
$\beta _{i}-1\leq m_{i}(s)\leq \beta _{i}$ 
for each $s\in S_{i}$ and the set 
$L_{i}:=\{s\in S_{i}:  m_{i}(s)=\beta _{i}\}$ is a basis of $M_{i}$.  
Let $M\sp{*}$ denote the direct sum of matroids $M_{1},\dots ,M_{q}$ and
let $\Delta \sp{*}\in {\bf Z}\sp{S}$ denote the translation vector
defined by 
\[
\Delta \sp{*}(s):= \beta_{i}-1 \ \ \hbox{whenever}\ \ s\in S_{i}, \ i=1,\dots ,q.
\]

By integrating these results, we obtain the following characterization.

\begin{theorem}   \label{matroid-eltolt} 
Let $B$ be an integral base-polyhedron.  
An element $m$ of (the M-convex set) $\odotZ{B}$ is decreasingly minimal 
if and only if 
$m$ can be obtained in the form $m=\chi_{L}+ \Delta \sp{*}$
where $L$ is a basis of the matroid $M\sp{*}$.  
The base-polyhedron $B\sp{\bullet}$ 
arises from the base-polyhedron of $M\sp{*}$ 
by adding the translation vector $\Delta\sp{*}$.  
Concisely, the set of dec-min elements of $\odotZ{B}$ is a
matroidal M-convex set.  
\finbox 
\end{theorem}

\paragraph{Cheapest dec-min element} 
An important algorithmic consequence of Theorems \ref{main.1} and \ref{matroid-eltolt} is 
that they help solve the 
{\bf cheapest dec-min element} problem, 
which is as follows.  
Let $c:S\rightarrow {\bf R} $ be a cost function and
consider the problem of computing a dec-min element $m$ of an M-convex
set $\odotZ{B}$ for which $cm$ is as small as possible.

By Theorem \ref{matroid-eltolt} the set 
$\odotZ{B\sp{\bullet}}$ 
of dec-min elements of
$\odotZ{B}$ can be obtained from a matroid $M\sp{*}$ 
by translation.
Namely, there is a vector $\Delta \sp{*} \in {\bf Z}\sp{S}$ 
such that $m$ is in 
$\odotZ{B\sp{\bullet}}$ 
if and only if 
there is a basis $L$ of $M\sp{*}$
for which  $m= \chi_{L} + \Delta \sp{*}$.  
Note that the matroid $M\sp{*}$
arises as the direct sum of matroids $M_{i}$ defined on the members
$S_{i}$ of the canonical partition.  
$M_{1}$ is described in Proposition \ref{PRbaseM1}  
and the other matroids $M_{i}$ may be determined analogously in an iterative way.  
To realize this algorithmically, 
we must have a strongly polynomial algorithm 
to compute the canonical partition 
as well as the essential value-sequence.  
Such an algorithm will be described in Section \ref{SCcompcanchain}.

Therefore, in order to find a minimum $c$-cost dec-min element of $\odotZ{B}$, 
it suffices to find a minimum $c$-cost basis of $M\sp{*}$.
Note that, in applying 
the greedy algorithm to the matroids $M_{i}$
in question, 
we need a rank oracle, which can be realized with the help
of a submodular function minimization oracle 
by relying on the definition of bases in \eqref{(p01)}.

Recall that for integral bounds $f \leq g$,
the intersection $B_{1}$ of a base-polyhedron $B$
and the box $T(f,g)$, if non-empty, 
is itself a base-polyhedron.
Therefore the algorithm above can be applied to the 
M-convex set $\odotZ{B_{1}}$, that is, 
we can compute a cheapest dec-min element of the intersection
$\odotZ{B_{1}} = \odotZ{B} \cap T(f,g)$.




\section{Integral square-sum and difference-sum minimization}
\label{SCsqdiffsum}

For a vector  $z \in {\bf Z}\sp{S}$, 
we can conceive several natural functions
to measure the uniformity of its component values $z(s)$
for $s \in S$.
Here are two examples:
\begin{align}
&\mbox{\bf square-sum}: \quad \ \  
W(z):=\sum [z(s)\sp{2}:s\in S],
\label{squaresumdef}
\\
&\mbox{\bf difference-sum}: \  
\Delta(z):=\sum [ | z(s) - z(t) |:s \not= t, \  s,t \in S].
\label{diffsumdef}
\end{align}
For vectors  $z_{1}$ and $z_{2}$ with
$\widetilde z_{1}(S)=\widetilde z_{2}(S)$,
$z_{1}$ may be felt more uniform than $z_{2}$ if $W(z_{1}) < W(z_{2})$,
and $z_{1}$ may also be felt more uniform 
if $\Delta(z_{1}) < \Delta(z_{2})$.
The first goal of this section is to show,
by establishing a fairly general theorem,
that a dec-min element of an M-convex set $\odotZ{B}$
is simultaneously a minimizer of these two functions.
The second goal of this section is to derive a
min-max formula for the minimum integral square-sum of an element of an M-convex set 
$\odotZ{B}$, along with characterizations of (integral) square-sum minimizers
and dual optimal solutions.

\subsection{Symmetric convex minimization} 
\label{SCnonsepmin}

Let $S$ be a non-empty ground-set of $n$ elements: 
$S = \{ 1,2,\ldots,n \}$.
We say that function $\Phi: \ZZ\sp{S} \to \RR$ is {\bf symmetric} if
\begin{equation} \label{symmtry}
 \Phi(z(1), z(2), \ldots, z(n)) = 
 \Phi(z(\sigma(1)), z(\sigma(2)), \ldots, z(\sigma(n))) 
\end{equation}
for all permutations $\sigma$ of $(1,2,\ldots,n)$.
We call a function $\Phi: \ZZ\sp{S} \to \RR$ {\bf convex} if
\begin{equation} \label{nonstrcnvI}
 \lambda \Phi(x) + (1-\lambda) \Phi(y) \geq  \Phi( \lambda x + (1-\lambda) y)
\end{equation}
whenever
$x, y \in \ZZ\sp{S}$, \ $0 < \lambda < 1$, and  
$\lambda x + (1-\lambda) y$ is an integral vector;
and {\bf strictly convex} if
\begin{eqnarray}  \label{strcnvI}
 \lambda \Phi(x) + (1-\lambda) \Phi(y) >  \Phi( \lambda x + (1-\lambda) y)
\end{eqnarray}
whenever
$x, y \in \ZZ\sp{S}$, \ $0 < \lambda < 1$, 
and  $\lambda x + (1-\lambda)  y$ is an integral vector.

In the special case where $\varphi$ is a function in one variable, 
it can easily be shown that the convexity of $\varphi$
is equivalent to the weaker requirement that the inequality
\begin{equation}  \label{univarconvfndef}
  2\varphi (k) \leq \varphi (k-1) + \varphi (k+1)
\end{equation}
holds for every integer $k$.
It is strictly convex in the sense of \eqref{strcnvI} if and only if 
$2\varphi (k) < \varphi (k-1) + \varphi (k+1)$ 
holds for every integer $k$.
For example, $\varphi(k)=k\sp{2}$ is strictly convex 
while $\varphi (k)=\vert k \vert$ is convex but not strictly.  
Given a function $\varphi $ in one variable, define $\Phi$ by 
\begin{equation}  \label{symsepconvfndef}
\Phi(z) := \sum [\varphi (z(s)):  s\in S]
\end{equation}
for $z \in {\bf Z}\sp{S}$.  
Such a function $\Phi$ is called a {\bf symmetric separable convex function};
note that $\Phi$ is indeed convex in the sense of \eqref{nonstrcnvI}.
When $\varphi $ is strictly convex, $\Phi$ is also called strictly convex.

\begin{example} \rm \label{RMsquaresum}
The square-sum $W(z)$ in \eqref{squaresumdef}
is a symmetric convex function
which is separable and strictly convex.
\finbox
\end{example}

\begin{example} \rm \label{RMdiffsum}
The difference-sum $\Delta(z)$ in \eqref{diffsumdef}
is a symmetric convex function
which is neither separable nor strictly convex.
More generally, for a nonnegative integer $K$, 
the function defined by
\[
\Delta_K(z) := \sum [ ( |z(s) - z(t)|  - K)\sp{+}:s \not= t, \  s,t \in S]
\]
is a symmetric convex function,
where $(x )\sp{+} = \max \{ x, 0 \}$.
\finbox
\end{example}

The following statements show a close relationship 
between decreasing minimality and 
the minimization of symmetric convex $\Phi$ over an M-convex set $\odotZ{B}$.

\begin{proposition} \label{decminIsPhiMin=nonsep} 
Let $B$ be an integral base-polyhedron
and $\Phi$ a symmetric convex function.  
Then each dec-min element of $\odotZ{B}$ 
is a minimizer of $\Phi$ over $\odotZ{B}$.  
\end{proposition} 

\Proof 
Since the dec-min elements of $\odotZ{B}$
are value-equivalent and $\Phi$ is symmetric,
the $\Phi$-value of each dec-min element is the same value $\mu$.
We claim that $\Phi(m)\geq \mu$ for each $m \in \odotZ{B}$.  
Suppose indirectly that there is an element $m$ of 
$\odotZ{B}$ for which $\Phi(m) < \mu$.
Then $m$ is not dec-min in $\odotZ{B}$
and Property (A) in Theorem~\ref{equi.1} implies that 
there is a 1-tightening step for $m$
resulting in decreasingly smaller member of $\odotZ{B}$,
that is,
there exist $s, t \in S$ such that
$m(t) \geq m(s)  + 2$
and 
$m' := m + \chi_{s} - \chi_{t} \in \odotZ{B}$.

Let $\alpha = m(t) - m(s)$, where $\alpha \geq 2$,
and define
$z = m + \alpha (\chi_{s} - \chi_{t})$.
Since $z$ is obtained from $m$ by 
interchanging the components at $s$ and $t$,
$\Phi(m) = \Phi(z)$ by symmetry (\ref{symmtry}).
Note that the vector $z$ may not be a member of $\odotZ{B}$.
For $\lambda = 1 - 1/\alpha$ we have
\begin{equation}  \label{decminnonsep=pr1}
 \lambda m + (1-\lambda)  z = 
\left( 1 - \frac{1}{\alpha} \right) m 
+ \frac{1}{\alpha} \left( m + \alpha (\chi_{s} - \chi_{t}) \right)  
=  m + \chi_{s} - \chi_{t} = m'
\in \odotZ{B} \ (\subseteq \ZZ^{S} ) ,
\end{equation}
from which 
$\lambda \Phi(m) + (1-\lambda) \Phi(z) \geq  \Phi(m')$
by convexity (\ref{nonstrcnvI}).
Since $\Phi(m) = \Phi(z)$,
this implies
$\Phi(m) \geq  \Phi(m')$.
After a finite number of
such 1-tightening steps, we arrive at a dec-min element
$m_{0}$ of $\odotZ{B}$,
for which $\mu = \Phi(m_{0})\leq \Phi(m) < \mu$, a contradiction.
\finbox

\medskip 

Note that if $\Phi $ is convex but not strictly convex,
then $\Phi$ may have minimizers that are not dec-min elements.  
This is exemplified by the identically zero function $\Phi$ for which
every member of $\odotZ{B}$ is a minimizer.  
However, for strictly convex functions we have the following characterization.

\begin{theorem}  \label{decmin=Phi=nonsep} 
Given an integral base-polyhedron $B$ and
a symmetric strictly convex function $\Phi$, 
an element $m$ of $\odotZ{B}$ is a minimizer of $\Phi$ 
if and only if 
$m$ is a dec-min element of $\odotZ{B}$.  
\end{theorem}  

\Proof
 If $m$ is a dec-min element, then $m$ is a $\Phi$-minimizer by
Proposition \ref{decminIsPhiMin=nonsep}.  
To see the converse, let $m$ be a
$\Phi$-minimizer of $\odotZ{B}$.  If, indirectly, $m$ is not a dec-min
element, then Property (A) in Theorem~\ref{equi.1} implies that 
there is a 1-tightening step for $m$,
that is,
there exist $s, t \in S$ such that
$m(t) \geq m(s)  + 2$
and 
$m' := m + \chi_{s} - \chi_{t} \in \odotZ{B}$.
For $\alpha = m(t) - m(s)$,
$\lambda = 1 - 1/\alpha$, and
$z = m + \alpha (\chi_{s} - \chi_{t})$,
we have
\eqref{decminnonsep=pr1},
from which we obtain
$\lambda \Phi(m) + (1-\lambda) \Phi(z) >  \Phi(m')$
by strict convexity (\ref{strcnvI}),
and hence
$\Phi(m) >  \Phi(m')$,
a contradiction to the assumption that
$m$ is a $\Phi$-minimizer.  
\finbox 

\medskip

We obtain the following  as corollaries of this theorem.

\begin{corollary}  \label{COdecmin=separable} 
Let $B$ be an integral base-polyhedron
and $\Phi$ a symmetric separable convex function.  
Then each dec-min element of $\odotZ{B}$ 
is a minimizer of $\Phi$ over $\odotZ{B}$,
and the converse is also true 
if, in addition,  $\Phi$ is strictly  convex.
\finbox
\end{corollary}  

\begin{corollary}  \label{COdecmin=squaresum} 
For an M-convex set $\odotZ{B}$,
an element $m$ of $\odotZ{B}$ is a square-sum minimizer 
if and only if 
$m$ is a dec-min element of $\odotZ{B}$.  
\finbox
\end{corollary}

A dec-min element is also characterized as a difference-sum minimizer.

\begin{theorem}  \label{THdecmin=diffsum} 
For an M-convex set $\odotZ{B}$,
an element $m$ of $\odotZ{B}$ is a difference-sum minimizer 
if and only if $m$ is a dec-min element of $\odotZ{B}$.  
\end{theorem}

\Proof
By Proposition \ref{decminIsPhiMin=nonsep}
every dec-min element is a difference-sum minimizer.
To show the converse, suppose indirectly that there is 
difference-sum minimizer $m$
that is not dec-min in $\odotZ{B}$.
Property (A) in Theorem~\ref{equi.1} implies that 
there is a 1-tightening step for $m$,
that is,
there exist $s, t \in S$ such that
$m(t) \geq m(s)  + 2$
and 
$m' := m + \chi_{s} - \chi_{t} \in \odotZ{B}$.
Here we observe that
$|m'(s) - m'(t)| =  |m(s) - m(t)|  -2$
and
\[
(|m'(v) - m'(s)| + |m'(v) - m'(t)|) - (|m(v) - m(s)| + |m(v) - m(t)|)
=
\begin{cases}
 -2 &  \hbox{if}\ \ m(s) < m(v) < m(t)
\\
  \phantom{-} 0 &  \hbox{otherwise}.
\end{cases} 
\]
This shows $\Delta(m') \leq \Delta(m) -2$, a contradiction.
\finbox

\begin{remark} \rm \label{RMcontvsdiscrete}
We emphasize that there is a fundamental difference between the
problems of finding a minimum square-sum element over a
base-polyhedron $B$ and over the M-convex set $\odotZ{B}$ (the set of
integral elements of $B$).  
In the first case (investigated by Fujishige \cite{Fujishige80,Fujishigebook}), 
there alway exists a single, unique solution, 
while in the second case, the square-sum minimizer elements
of $\odotZ{B}$ have an elegant matroidal structure.  
Namely, 
Corollary \ref{COdecmin=squaresum}
shows that the square-sum minimizers 
are exactly the dec-min elements of $\odotZ{B}$ and hence, 
by Theorem \ref{matroid-eltolt}, the set of square-sum minimizers of an
M-convex set arises from the bases of a matroid 
by translating their incidence vectors with a vector.
\finbox
\end{remark}

\begin{remark} \rm \label{RMorderingsqsum}
Corollary \ref{COdecmin=squaresum} says 
that an element $m$ of an M-convex set
$\odotZ{B}$ is dec-min precisely if $m$ is a square-sum minimizer.
One may feel that it would
have been a more natural approach to derive this equivalence by
showing that $x \leq_{\rm dec} y$ holds precisely if $W(x)\leq W(y)$.
Perhaps surprisingly, however, this equivalence fails to hold, that
is, the square-sum is not order-preserving with respect to the
quasi-order $\leq_{\rm dec}$.  
To see this, consider the following
four vectors in increasing order:
\[
 m_{1}=(2,3,3,1) <_{\rm dec} m_{2}=(3,3,3,0) <_{\rm dec} m_3=(2,2,4,1)
 <_{\rm dec} m_4=(3,2,4,0).
\]
Their square-sums admit a different order:
\[
W(m_{1})= 23, \quad W(m_{2})=27, \quad  W(m_3) = 25, \quad W(m_4) = 29.
\]
It is worth mentioning that the four vectors $m_{i}$ $(i=1,2,3,4)$
arise from the matroid $M_{1}$ given in Section \ref{egal2.example}.
Namely,
if $B_{1}'$ denotes the base-polyhedron
obtained from the base-polyhedron of $M_{1}$ by adding $(2,2,3,0)$, then
$\odotZ{B_{1}'}=\{m_{1},m_{2},m_3,m_4\}$.  
Among these four elements, $m_{1}$
is the unique dec-min element and the unique square-sum minimizer but
the decreasing-order and the square-sum order of the other three
elements are different.

We remark that if $\varphi$
in \eqref{symsepconvfndef}
is not only strictly convex but \lq rapidly\rq \ increasing 
as well, then $x <_{\rm dec} y$
can be proved to be equivalent to $\Phi(x) < \Phi (y)$.  
This intuitive notion of rapid increase will be formalized
in Part II \cite{Frank-Murota.2}.
\finbox
\end{remark}

\begin{remark} \rm \label{RMintersqsum}
For the intersection of two M-convex sets,
dec-min elements and square-sum minimizers may not coincide.  
Here is an example.
Let 
$\odotZ{B_{1}}=\{ (3,3,3,0 ), \ (2,2,4,1), \ (2,3,3,1), \ (3,2,4,0) \}$,
which is obtained by adding $(2,2,3,0)$ to 
the incidence vectors of the four bases of 
matroid $M_{1}$ in Section~\ref{egal2.example}.
Similarly, let
$\odotZ{B_{2}}=\{ (3,3,3,0 ), \ (2,2,4,1), \ (3,2,3,1), \ (2,3,4,0) \}$,
which is obtained by adding the same vector $(2,2,3,0)$ to 
the incidence vectors of the  bases of 
matroid $M_{2}$ in Section \ref{egal2.example}.
Let $\odotZ{B}$ denote the intersection of $\odotZ{B_{1}}$ and $\odotZ{B_{2}}$.
Then $\odotZ{B}= \{ (3,3,3,0), (2,2,4,1) \}$.  
Here $(3,3,3,0)$ is the
unique dec-min element while $(2,2,4,1)$ is the unique square-sum minimizer, 
demonstrating that the two notions of optima may differ 
for the intersection of two M-convex sets.
\finbox
\end{remark}

\begin{remark} \rm \label{RMgroenevelt}
An immediate consequence of 
Corollary \ref{COdecmin=separable}
is that a square-sum
minimizer of $\odotZ{B}$ minimizes an arbitrary symmetric separable
discrete convex function.  
Note, however, that this consequence immediately follows from a much earlier result by 
Groenevelt \cite{Groenevelt91} who characterized 
the elements of $\odotZ{B}$ minimizing a (not-necessarily symmetric) 
discrete convex function $\Phi$ defined for 
$z \in {\bf Z}\sp{S}$ by $\Phi(z) := \sum [\varphi_s(z(s)):  s\in S]$ 
(where $\varphi_s$ is a convex function in one variable for each $s\in S$).
\finbox
\end{remark}

\begin{remark} \rm \label{RMtwosepconv}
For $a,b,c \geq 0$, the function defined by
\[
\Phi(z) = a \sum_{s \in S} |z(s)| 
  + b \sum_{s \neq t} |z(s) - z(t)| + c \sum_{s \neq t} |z(s)+z(t)| 
\] 
is a symmetric convex function.
More generally, a function of the form 
\[
\Phi(z) = \sum_{s \in S} \varphi_{1}(z(s)) 
 + \sum_{s \neq t}  \varphi_{1}(|z(s) - z(t)|) + \sum_{s \neq t} \varphi_{3}(z(s)+z(t))  ,
\] 
where $\varphi_{1}, \varphi_{2}, \varphi_{3}: \ZZ \to \RR$,
are (discrete) convex functions,
is a symmetric convex function which is not separable.
Such a function is an example of the so-called 
2-separable convex functions.
By Theorem \ref{decmin=Phi=nonsep}, 
a dec-min element of $\odotZ{B}$ is a minimizer of function $\Phi$ over $\odotZ{B}$.
The minimization of 2-separable convex functions 
is investigated in depth by Hochbaum and others \cite{AHO04,Hoc02,Hoc07}
using network flow techniques.
\finbox
\end{remark}

\begin{remark} \rm \label{RMmaruyama}
Theorem \ref{decmin=Phi=nonsep} 
is a discrete counterpart of a result of 
Maruyama \cite{Mar78}
for the continuous case. See also Nagano \cite[Corollary 13]{Nag07}.
Symmetric convex function minimization is studied,
mainly for the continuous case,  
in the literature of majorization \cite{AS18}, \cite{MOA11}.
\finbox
\end{remark}

\begin{remark} \rm \label{RMnoduality}
A min-max formula  can be derived for the square-sum
and, more generally, for separable convex functions
from the Fenchel-type duality theorem in DCA.
However, we cannot use the Fenchel-type duality theorem 
to obtain a min-max formula for non-separable symmetric convex functions,
since non-separable symmetric convex functions are not necessarily
{\rm M}$\sp{\natural}$-convex.
\finbox
\end{remark}

\subsection{Min-max theorem for integral square-sum}
\label{SCsquare.minmax}

Recall the notation
$W(z)=\sum [z(s)\sp{2}:s\in S]$
for the square-sum of $z \in {\bf Z}\sp{S}$.  
Given a polyhedron $B$, 
we say that an element $m\in \odotZ{B}$ is a {\bf square-sum minimizer} 
(over $\odotZ{B}$) or that $m$ is an {\bf integral square-sum minimizer} 
of $B$ 
if $W(m)\leq W(z)$ holds for each $z\in \odotZ{B}$.  
The main goal of this section is to derive a
min-max formula for the minimum integral square-sum of an element of an M-convex set 
$\odotZ{B}$, along with a characterization of (integral) square-sum minimizers.

A set-function $p$ on $S$ can be considered as a function defined on $(0,1)$-vectors.  
It is known that $p$ can be extended in a natural
way to every vector $\pi$ in ${\bf R}\sp{S}$, as follows.  
For the sake of this definition, we may assume that the elements of $S$ are
indexed in a decreasing order of the components of $\pi$, that is,
$\pi(s_{1})\geq \cdots \geq \pi(s_{n})$ (where the order of the
components of $\pi$ with the same value is arbitrary).
For $j=1,\dots ,n$,  let $I_{j}:=\{s_{1},\dots ,s_{j}\}$ and let
\begin{equation}
\hat p(\pi):= p(I_{n})\pi(s_{n}) 
 + \sum_{j=1}\sp {n-1} p(I_{j})[\pi(s_{j})-\pi(s_{j+1})].
\end{equation}
Obviously, $p(Z)=\hat p(\chi_Z)$.  
The function $\hat p$ is called the {\bf linear extension} of $p$. 

\begin{remark} \rm \label{RMlovext}
The linear extension was first
considered by Edmonds \cite{Edmonds70} who proved for a polymatroid
$P=P(b)$ defined by a monotone, non-decreasing submodular function $b$ that 
$\max \{\pi x :  x\in \odotZ{P}\} = \hat b(\pi )$ when $\pi$ 
is non-negative.  
The same approach shows for a base-polyhedron $B=B'(p)$
defined by a supermodular function $p$ that 
$\min \{\pi x :  x\in \odotZ{B}\} = \hat p(\pi)$.  
Another basic result is due to Lov\'asz \cite{Lovasz83} 
who proved that $p$ is submodular 
if and only if 
$\hat p$ is convex.  
We do not, however, explicitly need these results, 
and only remark that in the literature the linear extension is often called
Lov\'asz extension.
\finbox
\end{remark}

\medskip

Our approach is as follows.  First, we consider an arbitrary
set-function $p$ on $S$ (supermodular or not) along with the
polyhedron
\[
 B=B'(p):=\{x:  x\in {\bf R}\sp{S}, \ \widetilde x(Z)\geq p(Z)
 \ \mbox{for every $Z\subset S$ and $\widetilde x(S)=p(S)$} \}, 
\]
and develop an easily checkable lower bound for the minimum
square-sum over the integral elements of $B$.  If this lower bound is
attained by an element $m$ of $\odotZ{B}$, then $m$ is certainly a
square-sum minimizer independently of any particular property of $p$.
For general $p$, the lower bound (not surprisingly) is not 
always attainable.  
We shall prove, however, that it is attainable when $p$
is fully supermodular.  That is, we will have a min-max theorem for
the minimum square-sum over an M-convex set $\odotZ{B}$, or in other
words, we will have an easily checkable certificate for an element $m$
of $\odotZ{B}$ to be a minimizer of the square-sum.

\begin{remark} \rm \label{RManalogy}
For readers familiar with K{\H o}nig's matching theorem,
we mention an analogy for our approach.  
Consider the well-known problem of finding the minimum number $\tau $ of 
nodes of an arbitrary graph $G=(V,E)$ hitting all the edges.  
If $M$ is any matching of $G$,
then $\vert M\vert $ is a lower bound for $\tau $. 
This implies that 
if $T$ is a subset of nodes hitting $E$ and $M$ is a matching of $G$
for which $\vert T\vert =\vert M\vert $, 
then $T$ is certainly a minimum cardinality subset of nodes hitting all the edges.  
That is, $M$ is a certificate for the minimality of $T$.  
Obviously, such a certificate does not always exist 
when $G$ is arbitrary (as demonstrated by a triangle), 
but K{\H o}nig' classic theorem asserts
that a bipartite graph always includes a subset $T$ of nodes hitting
$E$ and a matching $M$ for which $\vert T\vert =\vert M\vert $, or in
other words, $\tau $ is equal to $\nu $, the maximum cardinality of a matching.
\finbox
\end{remark}

We shall need the following two claims.

\begin{claim}   \label{hcestim}
 For $m,\pi \in {\bf Z}\sp{S}$, one has 
\begin{equation} 
\sum_{s\in S} \llfloor {\pi(s) \over 2}\rrfloor \llceil {\pi(s) \over 2}\rrceil 
    \geq \sum_{s\in S} m(s)[\pi(s)-m(s)].  
\label{(hcestim.1)}
\end{equation}
Moreover, equality holds 
if and only if 
\begin{equation} 
 m(s)\in \bigg\{ \llfloor {\pi(s) \over 2}\rrfloor , 
                 \llceil {\pi(s) \over 2}\rrceil \bigg\} \ \
   \ \hbox{\rm for every} \ s\in S.
\label{(hcestim.2)} 
\end{equation} 
\end{claim}

\Proof 
The claim follows by observing that 
$\llfloor {a \over 2}\rrfloor \llceil {a \over 2}\rrceil \geq b(a-b)$ 
holds for any pair of integers
$a$ and $b$, where equality holds precisely if 
$b \in \big\{\llfloor {a \over 2}\rrfloor , \llceil {a \over 2}\rrceil \big\}$.  
\finbox

\medskip

Let $p$ be an arbitrary set-function on $S$ with $p(\emptyset )=0$ and
consider an integral element $m$ of the polyhedron $B=B'(p)$.  
Recall that a non-empty subset $X\subseteq S$ was called 
a strict $\pi$-top set if $\pi(u)>\pi(v)$ held whenever $u\in X$ and $v\in S-X$.  
In what follows, for an $m\in \odotZ{B}$, $m$-tightness of a subset
$Z\subseteq S$ means $\widetilde m(Z)=p(Z)$.

\begin{claim}   \label{pimpi} 
For $m\in \odotZ{B}$ and $\pi \in {\bf Z}\sp{S}$, one has 
\begin{equation} 
\hat p(\pi) \leq \sum_{s\in S}m(s)\pi(s).
\label{(pimpi)} 
\end{equation} 
Moreover, equality holds 
if and only if 
each (of the at most $n$) strict $\pi$-top set is $m$-tight.  
\end{claim}

\Proof 
Suppose that the elements of $S$ are indexed in such a way that
$\pi(s_{1})\geq \pi(s_{2})\geq \cdots \geq \pi(s_{n})$.  
For $j=1,\dots ,n$,   let $I_{j}:=\{s_{1},\dots ,s_{j}\}$.   
Then
\begin{eqnarray*} 
\hat p(\pi) &=& p(I_{n})\pi(s_{n}) 
    + \sum_{j=1}\sp{n-1} p(I_{j})[\pi(s_{j})-\pi(s_{j+1})] 
\\ &\leq & 
\widetilde m(I_{n})\pi(s_{n}) 
   + \sum_{j=1}\sp{n-1} \widetilde m(I_{j})[\pi(s_{j})-\pi(s_{j+1})]
\\ &=& 
\sum_{1\leq i\leq j\leq n} m(s_{i})\pi(s_{j}) 
   - \sum_{1\leq i\leq j\leq n-1} m(s_{i})\pi(s_{j+1}) 
\\ &=&
 \sum_{1\leq i\leq j\leq n} m(s_{i})\pi(s_{j})  
  - \sum_{1\leq i< j'\leq n} m(s_{i})\pi(s_{j'}) 
\\ &=& 
\sum_{j=1}\sp {n} m(s_{j})\pi(s_{j}), 
\end{eqnarray*}
from which \eqref{(pimpi)} follows.  
Furthermore, we have equality in \eqref{(pimpi)} 
precisely if $\widetilde m(I_{j})=p(I_{j})$
holds whenever $\pi(s_{j})-\pi(s_{j+1})>0$.  
But this latter condition is equivalent to requiring 
that each strict $\pi$-top set is $m$-tight.  
\finbox 
\medskip

\begin{proposition}   \label{lowerbound}
Let $p$ be an arbitrary set-function on $S$ with $p(\emptyset )=0$ 
and let $m$ be an integral element of the polyhedron $B=B'(p)$.
Then
\begin{equation} 
 \sum_{s\in S} m(s)\sp{2}
 \geq \hat p(\pi) - \sum_{s\in S} \llfloor {\pi(s) \over 2}\rrfloor 
 \llceil {\pi(s) \over 2}\rrceil 
\label{(minleqmax)} 
\end{equation}
whenever $\pi \in {\bf Z}\sp{S}$ is an integral vector.
Furthermore, equality holds for $m$ and $\pi$ 
if and only if 
the following {\bf optimality criteria} hold:
\begin{eqnarray} 
{\rm (O1)} & & \hbox{\eqref{(hcestim.2)} holds:} \ \
  m(s)\in 
  \bigg\{ \llfloor {\pi(s) \over 2}\rrfloor , \llceil {\pi(s) \over 2}\rrceil \bigg\} 
 \ \ \hbox{for every }\  s\in S , 
\label{(optcrit.1)} 
\\ 
{\rm (O2)} & & \hbox{each strict $\pi$-top-set is $m$-tight with respect to $p$}.  
\label{(optcrit.2)}
\end{eqnarray}
\end{proposition}  

\Proof 
Let $\pi \in {\bf Z}\sp{S}$.  By the two preceding claims,
\begin{equation} 
\sum_{s\in S} m(s)\sp{2} 
   = \sum_{s\in S} m(s)\pi(s) 
   - \sum_{s\in S} m(s) [\pi(s)-m(s)] 
   \geq \hat p(\pi) 
   - \sum_{s\in S} \llfloor {\pi(s) \over 2}\rrfloor \llceil {\pi(s) \over 2}\rrceil ,
\label{(westim)} 
\end{equation} 
from which \eqref{(minleqmax)} follows.

The claims also immediately imply that we have equality in
\eqref{(minleqmax)} precisely if the optimality criteria (O1) and (O2) hold.  
\finbox
\medskip

The min-max formula in the next theorem concerning min square-sum over
the integral elements of an integral base-polyhedron can be derived from 
the much more general Fenchel-type duality theorem in DCA 
(Discrete Convex Analysis), due to Murota 
(see \cite{Murota98a} and also Theorem 8.21, page 222, in the book \cite{Murota03}).  
The advantage of the present min-max formula is
that it does not need the notion of conjugate functions which is an
essential part of the general result in \cite{Murota98a}.  
In addition, our proof relies on the characterization of dec-min elements
described in Theorem \ref{equi.1} and does not need the general tools of DCA.

\begin{theorem}   \label{minmax}
Let $B=B'(p)$ be a base-polyhedron defined by
an integer-valued fully supermodular function $p$.  
Then
\begin{equation} 
 \min \{ \sum_{s\in S} m(s)\sp{2}:  m\in \odotZ{B} \} 
 = \max \{\hat p(\pi) - \sum_{s\in S} \llfloor {\pi(s) \over 2}\rrfloor 
 \llceil {\pi(s) \over 2}\rrceil :  \pi \in {\bf Z}\sp{S} \}.  
\label{(minmax)} 
\end{equation}
\end{theorem}

\Proof 
By Proposition \ref{lowerbound}, 
$\min\geq \max$ holds in \eqref{(minmax)} and hence all what 
we have to prove is that there is
an element $m\in \odotZ{B}$ and an integral vector 
$\pi \in {\bf Z}\sp{S}$ meeting the two optimality criteria formulated 
in Proposition \ref{lowerbound}.
Let $m$ be an arbitrary dec-min element of $\odotZ{B}$.  
By Property
(B) of Theorem \ref{equi.1}, there is a chain $(\emptyset \subset ) \
C_{1}\subset C_{2}\subset \cdots \subset C_\ell= S$ of $m$-tight and
$m$-top sets for which the restrictions of $m$ onto the difference sets 
$S_{i}:=C_{i}-C_{i-1}$ ($i=1,\dots ,\ell$) are near-uniform in $S_{i}$
(where $C_{0}:=\emptyset $).  
Note that $\{S_{1},\dots ,S_\ell\}$ is a partition of $S$.

For $i=1,\dots ,\ell$, let $\beta_{i}(m):= \max \{m(s):  s\in S_{i}\}$.
Define $\pi_{m}:S\rightarrow {\bf Z} $ by
\[
 \pi_{m}(s):= 2\beta_{i}(m)-1  \ \hbox{ if } \ s\in S_{i} \ \ (i=1,\dots ,\ell)  . 
\]
We have 
\[
\llfloor {\pi_{m}(s)/ 2}\rrfloor 
  = \beta_{i}(m)-1 \leq m(s)
 \leq \beta_{i}(m)= \llceil {\pi_{m}(s)/2} \rrceil 
\]
for every $s\in S_{i}$, and
hence Optimality criterion (O1) holds for $m$ and $\pi_{m}$.

We claim that each strict $\pi_{m}$-top set $Z$ is a member of chain $\cal C$.  
Indeed, as $\pi_{m}$ is uniform in each $S_{j}$, 
if $Z$ contains an element of $S_{j}$, then $Z$ includes the whole $S_{j}$.
Furthermore, since each member of $\cal C$ is an $m$-top set, we have
$\beta_{1}(m)\geq \beta_{2}(m)\geq \cdots \geq \beta_\ell(m)$, 
and hence if
$Z$ includes $S_{j}$, then it includes each $S_{i}$ with $i<j$.  
Therefore
every strict $\pi_{m}$-top set is indeed a member of the chain,
implying Optimality criterion (O2).  
\finbox

\medskip

It should be noted that the optimal dual solution $\pi_{m}$ obtained in
the proof of the theorem is actually an {\bf odd} vector in the sense
that each of its component is an odd integer.

\begin{corollary}   \label{minmaxb} 
There is an odd dual optimizer $\pi$ in
the min-max formula \eqref{(minmax)}, that is, the min-max formula in
Theorem {\rm \ref{minmax}} can be re-written as follows:
\begin{equation} 
\min \{ \sum_{s\in S} m(s)\sp{2}:  m\in \odotZ{B} \}
 = \max \{\hat p(\pi) - \sum_{s\in S} { {\pi(s)\sp{2} - 1} \over 4} : 
  \ \pi \in {\bf Z}\sp{S}, \ \pi \ \hbox{\rm is odd} \ \}.  
\label{(minmax.odd)} 
\end{equation}
\end{corollary}

We emphasize that for the proof of Theorem \ref{minmax} and Corollary
\ref{minmaxb} we relied only on Theorem~\ref{equi.1} and did not need
the characterization of the set of dec-min elements of $\odotZ{B}$ 
given in Section~\ref{canonical}.

In the proof of Theorem \ref{minmax}, we chose an arbitrary dec-min
element $m$ of $\odotZ{B}$ and an arbitrary chain of $m$-tight and
$m$-top sets 
such that $m$ is near-uniform on each difference set.  
In Section \ref{canonical}, we proved that there is a single canonical
chain ${\cal C}\sp{*}$ which meets these properties for every dec-min
element of $\odotZ{B}$.  
Therefore the dual optimal $\pi\sp{*}$ assigned to ${\cal C}\sp{*}$ 
is also independent of $m$.  
Namely, consider the canonical $S$-partition $\{S_{1},\dots ,S_{q}\}$ and
the essential value-sequence $\beta_{1}>\cdots >\beta_{q}$.  
Define $\pi\sp{*}$ by
\begin{equation} 
\pi\sp{*}(s):= 2\beta_{i}-1 \ \hbox{ if } \ s\in S_{i} \ (i=1,\dots ,q).  
\label{(pi*def)} 
\end{equation}
As we pointed out in the proof of Theorem \ref{minmax}, this $\pi\sp{*}$ 
is also a dual optimum in \eqref{(minmax)}.  
We shall prove in the next section that 
$\pi\sp{*}$ is actually the unique smallest dual
optimum in \eqref{(minmax)}.

\subsection{The set of optimal duals to integral square-sum minimization}
\label{SCoptdualset}

We proved earlier that an element $m\in \odotZ{B}$ is a square-sum
minimizer precisely if it is a dec-min element.  
This and Theorem \ref{main.1} imply that the square-sum minimizers 
of $\odotZ{B}$ are the integral members of a base-polyhedron $B\sp{\bullet}$
obtained by intersecting a particular face of $B$ with a special small box.  
This means that the integral square-sum minimizers form an M-convex set.

One of the equivalent definitions of an L$\sp {\natural}$-convex set $L$ 
(pronounce L-natural convex) in Discrete Convex Analysis is that
$L$ is the set of integer-valued feasible potentials.  
Formally,
$L=\{\pi \in {\bf Z}\sp{S}:  \pi(v)-\pi(u) \leq g(uv)\}$
where $g$ is an integer-valued function on the ordered pairs of elements of $S$.
(By a theorem of Gallai, $L$ is non-empty 
if and only if 
there is no dicircuit of negative total $g$-weight.)

Our next goal is to show that the dual optima in Theorem \ref{minmax}
form an L$\sp {\natural}$-convex set $\Pi$, and we
provide a description of $\Pi$ 
as the integral solution set of feasible potentials in a box.

Recall that the optimality criteria for a dec-min element $m$ of
$\odotZ{B}$ and for an integral vector $\pi$ were 
given by (O1) and (O2) in \eqref{(optcrit.1)}--\eqref{(optcrit.2)}. 
These immediately imply the following.

\begin{proposition}  
For an integral vector $\pi$, the following are equivalent.  
\smallskip

\noindent 
{\rm (A)} \ 
$\pi$ is a dual optimum (that is, $\pi$ belongs to $\Pi$).  
\smallskip

\noindent 
{\rm (B)} \ 
There is a dec-min element $m$ of $\odotZ{B}$
such that $m$ and $\pi$ meet the optimality criteria.
\smallskip

\noindent
{\rm (C)} \ 
For every dec-min $m$ of $\odotZ{B}$, 
$m$ and $\pi$ meet the optimality criteria.  
\finbox
\end{proposition}

Consider the canonical $S$-partition $\{S_{1},\dots ,S_{q}\}$, 
the essential value-sequence 
$\beta_{1}>\beta_{2}>\cdots >\beta_{q}$,
and the matroids $M_{i}$ on $S_{i}$ $(i=1,\dots ,q)$.  
We can use the notions and apply the results of Section~\ref{value-fixed} 
formulated for $M_{1}$ to each $M_{i}$ \ $(i=1,\dots ,q)$.  
To follow the pattern of ${\cal F}_{1}$ introduced in \eqref{(scriptF1)}, let 
\begin{equation}
 {\cal F}_{i}:= \{X\subseteq S_{i}:  \ \beta_{i}\vert X\vert = p_{i}(X)\} ,
\end{equation} 
where $p_{i}$ was defined by $p_{i}(X)=p(C_{i-1}\cup X) - p(C_{i-1})$ 
for $X\subseteq S_{i}$.  
Since $\beta_{i}\vert X\vert \geq p_{i}(X)$ for every $X\subseteq S_{i}$
and $p_{i}$ is supermodular, 
${\cal F}_{i}$ is closed under taking intersection and union.  
Let $F_{i}$ denote the unique largest member of ${\cal F}_{i}$, 
that is, $F_{i}$ is the union of the members of ${\cal F}_{i}$.  
Both $F_{i}=\emptyset $ and $F_{i}=S_{i}$ are possible.

\begin{theorem}  \label{co-loopi} 
For an element $s\in S_{i}$ \ $(i=1,\dots ,q)$, 
the following properties are pairwise equivalent.%
\smallskip

\noindent 
{\rm (A)} \ $s$ is value-fixed.  
\smallskip

\noindent 
{\rm (B)} \ $m(s)=\beta_{i}$ holds for every dec-min element $m$ of $\odotZ{B}$.  
\smallskip

\noindent
 {\rm (C)} \ $s\in F_{i}$.  
\smallskip

\noindent
 {\rm (D)} \ $s$ is a co-loop of $M_{i}$.  
\finbox 
\end{theorem}

Define a digraph $D_{i}=(F_{i}, A_{i})$ on node-set $F_{i}$ in which $st$ is
an arc if $s,t\in F_{i}$ and there is no $t\ol s$-set in ${\cal F}_{i}$.
This implies that no arc of $D_{i}$ enters any member of ${\cal F}_{i}$.

\begin{theorem}  \label{opt.dual.set} 
An integral vector $\pi \in {\bf Z}\sp{S}$ 
is an optimal dual solution to the integral minimum square-sum problem 
(that is,  $\pi \in \Pi$) 
if and only if 
the following three conditions hold for each $i=1,\dots ,q:$
\begin{eqnarray} 
&& \pi(s)=2\beta_{i}-1  \quad \hbox{\rm  for every}\ \  s\in S_{i}-F_{i}, 
\label{(Si-Fi)} 
\\&& 
2\beta_{i}-1\leq \pi(s) \leq 2\beta_{i}+1 \quad \hbox{\rm for every}\ \ s\in F_{i}, 
\label{(Fi)} 
\\&& 
\pi(s)-\pi(t) \geq 0 \quad \hbox{\rm whenever \ $s,t\in F_{i}$ \  and \  $st\in A_{i}$}.  
\label{(pispit)} 
\end{eqnarray} 
\end{theorem}

\Proof
\begin{claim}  
\label{O1'} 
Optimality criterion {\rm (O1)} is equivalent to
\begin{equation} 
{\rm (O1')} \ \ \ \ \ 2m(s)-1 \leq \pi(s) \leq 2m(s)+1 \ \ 
\hbox{\rm for \ $s\in S$}. 
\label{(optcrit.1')} 
\end{equation} 
\end{claim} 

\Proof 
When $\pi(s)$ is even, we have the following equivalences:
\begin{eqnarray*} 
m(s)\in \bigg\{ \llfloor {\pi(s) \over 2}\rrfloor , 
           \llceil {\pi(s) \over 2}\rrceil \bigg\}
 & \Leftrightarrow & \pi(s)=2m(s) 
\\ &
\Leftrightarrow & 2m(s)-1 \leq \pi(s) \leq 2m(s)+1.  
\end{eqnarray*}
When $\pi(s)$ is odd, we have the following equivalences:
\begin{eqnarray*} 
m(s)\in \bigg\{ \llfloor {\pi(s) \over 2}\rrfloor , 
                \llceil {\pi(s) \over 2}\rrceil \bigg\} 
& \Leftrightarrow & \pi(s)-1 \leq 2m(s) \leq \pi(s)+1 \ 
\\ & \Leftrightarrow & 2m(s)-1 \leq \pi(s) \leq 2m(s)+1.  
\end{eqnarray*}
\vspace{-2.8\baselineskip} \\
\finbox
\vspace{1.2\baselineskip}

Suppose first that $\pi \in {\bf Z}\sp{S}$ is an optimal dual solution.
Then the optimality criteria (O1$'$) and (O2) formulated in
\eqref{(optcrit.1')} and \eqref{(optcrit.2)} 
hold for every dec-min element $m$ of $\odotZ{B}$.

Let $s$ be an element of $S_{i}-F_{i}$.  
Since $s$ is not value-fixed, 
there are dec-min elements $m$ and $m'$ of $\odotZ{B}$ for which
$m(s)=\beta_{i}-1$ and $m'(s)=\beta_{i}$.  
By applying \eqref{(optcrit.1')} to $m$ and to $m'$, we obtain that
\[
 2\beta_{i}-1 = 2m'(s)-1\leq \pi(s)\leq 2m(s)+1 
   = 2(\beta_{i}-1) + 1 = 2\beta_{i}-1 ,
\]
from which $\pi(s)=2\beta_{i}-1$ follows, and hence \eqref{(Si-Fi)} holds indeed.

Let $s$ be an element of $F_{i}$.  
As $s$ is value-fixed, $m(s)=\beta_{i}$ holds 
for any dec-min element $m$ of $\odotZ{B}$.  
We obtain from \eqref{(optcrit.1')} that
\[
 2\beta_{i}-1 = 2m(s)-1\leq \pi(s)\leq 2m(s)+1 = 2\beta_{i} + 1 
\]
and hence \eqref{(Fi)} holds.

To derive \eqref{(pispit)}, suppose indirectly that $st$ is an arc in $A_{i}$ 
for which $\pi(t)>\pi(s)\geq 2\beta_{i}-1$.  
Let $Z:=\{v\in S: \pi(v)\geq \pi(t)\}$.  
Then $Z$ is a strict $\pi$-top set and hence
$C_{i-1}\subseteq Z\subseteq C_{i-1}\cup F_{i}-s$.  
By Optimality criterion (O2), $Z$ is $m$-tight with respect to $p$.  
Let $X:=Z\cap S_{i}$.  
Then $X\subseteq F_{i}$ and hence
\[
p(Z) = \widetilde m(Z) = \widetilde m(C_{i-1}) + \widetilde m(X) 
     = p(C_{i-1}) + \beta_{i}\vert X\vert , 
\]
from which 
\[
  \beta_{i}\vert X\vert = p(Z) - p(C_{i-1}) =p_{i}(X),
\]
that is, $X$ is in ${\cal F}_{i}$, in contradiction with the
definition of $A_{i}$ which requires that $st$ enters no member of ${\cal F}_{i}$.

Suppose now that $\pi$ meets the three properties formulated in
Theorem \ref{opt.dual.set}.  Let $m\in \odotZ{B}$ be an arbitrary dec-min element.  
Consider an element $s$ of $S_{i}$.  
If $s\in F_{i}$, that is, if $s$ is value-fixed, then $m(s)=\beta_{i}$.  
By \eqref{(Fi)}, we have $2m(s)-1\leq \pi(s)\leq 2m(s)+1$, 
that is, Optimality criterion (O1$'$) holds.  
If $s\in S_{i}-F_{i}$, then $\pi(s)=2\beta_{i}-1$ by \eqref{(Si-Fi)}, 
from which 
\[
 \llfloor {\pi(s) \over 2}\rrfloor = {\pi(s)-1 \over 2} 
 \ = \  \beta_{i}-1
 \ \leq \  m(s) 
 \ \leq \  \beta_{i} 
 \ = \  {\pi(s) +1 \over 2} = \llceil {\pi(s) \over 2}\rrceil ,
\]
showing that Optimality criterion (O1$'$) holds.

To prove optimality criterion (O2), let $Z$ be a strict $\pi$-top set
and let $\mu := \min\{ \pi(v):  v\in Z\}$.  
 Let $i$ denote the largest subscript 
for which $X:=Z\cap S_{i}\not =\emptyset $. 
Then $\mu \leq 2\beta_{i}+1 \leq 2\beta_{i-1}-1\leq \pi(u)$ 
holds for every $u\in C_{i-1}$,  
from which $C_{i-1}\subseteq Z$ as $Z$ is a strict $\pi$-top set.

If $\mu =2\beta_{i}-1$, then $S_{i}\subseteq Z$ as $Z$ 
is a strict $\pi$-top set, from which 
$Z=C_{i}$, implying that $Z$ is an $m$-tight set in this case.  
Therefore we suppose 
$\mu \geq 2\beta_{i}$, from which $X\subseteq F_{i}$ follows.  
Now $X\in {\cal F}_{i}$, for otherwise there is an arc $st\in A_{i}$ $(s,t\in F_{i})$ 
entering $X$, and then $\pi(t)\leq \pi(s)$ holds by Property \eqref{(pispit)}; 
 this contradicts the assumption that $Z$ is a strict $\pi$-top set.  
By $X\in {\cal F}_{i}$ we have 
$\beta_{i}\vert X\vert =p_{i}(X)$ and hence
\begin{eqnarray*} 
\widetilde m(Z) &=& \widetilde m(X)+ \widetilde m(C_{i-1}) 
 = \beta_{i}\vert X\vert + p(C_{i-1}) 
\\ &=& p_{i}(X) + p(C_{i-1}) 
= p(X\cup C_{i-1}) - p(C_{i-1}) + p(C_{i-1}) = p(Z),
\end{eqnarray*}
that is, $Z$ is indeed $m$-tight.  
\finbox \finboxHere

\medskip \medskip 
In \eqref{(pi*def)}, we defined a special dual optimal solution 
$\pi\sp{*}$ by $\pi\sp{*}(s) =2\beta_{i}-1$ 
whenever $s\in S_{i}$ ($i=1,\dots ,q$).  
Theorem \ref{opt.dual.set} and the definition we use 
for L$\sp {\natural}$-convex sets immediately implies the following.

\begin{corollary}   \label{Pistar} 
The set $\Pi$ of optimal dual integral vectors $\pi$ 
in the min-max formula \eqref{(minmax)} of 
Theorem~{\rm \ref{minmax}}
is an {\rm L}$\sp{\natural}$-convex set.  
The unique smallest element of $\Pi$ (that is, the unique smallest
dual optimum) is $\pi\sp{*}$.  
\finbox 
\end{corollary}



\section{Algorithms} 
\label{SCalgo2}

In this section, we consider algorithmic aspects of the problems
investigated so far, and show how to compute efficiently 
a decreasingly minimal element of an M-convex set along with 
its canonical chain (or $S$-partition).

Let $B$ be a non-empty integral base-polyhedron.  
As mentioned earlier, $B$ can be given in the form $B(b)$ where $b$ 
is a (fully) submodular function or in the form $B'(p)$ 
where $p$ is a (fully) supermodular function.  
Here $b$ and $p$ are complementary functions
(that is, $p(X)=b(S)-b(S-X))$ and hence an algorithm described for one
of them can easily be transformed to work on the other.  
In the present description, we use supermodular functions with the remark
that in applications base-polyhedra are often given with $b$.

There is a one-to-one correspondence between $B$ and $p$ but, 
as mentioned earlier, for an intersecting or crossing supermodular function $p$, 
\ $B'(p)$ is also a (possibly empty) base-polyhedron
which is integral if $p$ is integer-valued.  
As already indicated, for obtaining and proving results for $B$ 
(or for $\odotZ{B}$), 
it is much easier to work with a fully supermodular $p$
while in applications base-polyhedra are often 
arise from---intersecting or crossing (or even weaker)---supermodular functions.
Therefore in describing and analysing algorithms, we must consider
these weaker functions as well.

One of the most fundamental algorithms of discrete optimization is 
for minimizing a submodular function, that is, 
for finding a subset $Z$ of $S$ for which 
$b(Z)=\min \{b(X):X\subseteq S\}$.  
There are strongly
polynomial algorithms for this problem 
(for example, Schrijver \cite{Schrijver2000} and Iwata et al.~\cite{IFF01} 
are the first, while Orlin \cite{Orlin09} is one of the fastest), 
and we shall refer to such an algorithm as a {\bf submod-minimizer} subroutine.  The
complexity of Orlin's algorithm \cite{Orlin09}, for example, is
$O(n\sp 6)$ (where $n=\vert S\vert $) and the algorithm calls 
$O(n\sp{5})$ times a routine which evaluates the submodular function in question.  
(An evaluation routine outputs the value $b(X)$ for any input subset $X\subseteq S$).  
This complexity bound is definitely
attractive from a theoretical point of view but in concrete applications 
it is always a challenge to develop faster algorithms for the special case.  
Naturally, submodular function minimization and
supermodular function maximization are equivalent.

\subsection{The basic algorithm for computing a dec-min element}
\label{SCbasicalg}

Our first goal is to describe a natural 
approach---the {\bf basic algorithm}---for finding a decreasingly minimal element 
of an M-convex set $\odotZ{B}$.  
The basic algorithm is polynomial in
$n+\vert p(S)\vert $, 
and hence it is polynomial in $n$ 
when $\vert p(S)\vert $ is small in the sense that it can be bounded 
by a polynomial of $n$.  
This is the case, for example, in an application
when we are interested in strongly connected decreasingly minimal
(=egalitarian) orientations.  
In the general case, where typical
applications arise by defining $p$ with a \lq large\rq \ capacity function, 
a (more complex) strongly polynomial-time algorithm will be
described in the next section.

In order to find a dec-min element of an M-convex set $\odotZ{B}$, 
we assume that a subroutine is available to 
\begin{equation} 
\hbox{compute an integral element of $B$.}\ 
\label{(member)} 
\end{equation} 
When $B=B'(p)$ and $p$ is
(fully) supermodular, a variant
 of Edmonds' polymatroid greedy algorithm 
finds an integral member of $B$. 
 (Namely, take any ordering $s_{1},\dots ,s_{n}$ of $S$, and define 
$m(s_{1}):=p(s_{1})$ and, for $i=2,\dots ,n$, 
$m(s_{i})=p(Z_{i})-p(Z_{i-1})$ where $Z_{i}=\{s_{1},s_{2},\dots,s_{i}\}$.  
Edmonds \cite{Edmonds70} proved that vector $m$ is indeed in $B$).  
This algorithm needs only a subroutine to evaluate $p(Z_{i})$ for
$i=1,\dots ,n$.  If $p$ is intersecting supermodular, 
then Frank and Tardos \cite{FrankJ17} described an algorithm 
which needs $n$ applications of a submod-minimizer routine.  
For crossing supermodular $p$, a more complex algorithm is given in \cite{FrankJ17} 
which terminates after at most $n\sp{2}$ applications of a submod-minimizer.
Note that the latter problem of finding an integral element of a
base-polyhedron $B'(p)$ defined by a crossing supermodular function
$p$ covers such non-trivial problems as the one of finding a
degree-constrained $k$-edge-connected orientation of an undirected
graph, a problem solved first in \cite{FrankJ4}.

Suppose now that an initial integral member $m$ of $B$ is available.
The algorithm needs a subroutine to 
\begin{equation} 
\hbox{ decide for  $m\in \odotZ{B}$ and for $s,t\in S$ if $m':=m+\chi_{s}-\chi_{t}$ 
belongs to $B$.}\ 
\label{(st.step)} 
\end{equation}
Observe that Subroutine \eqref{(st.step)} is certainly available if we can
\begin{equation} 
\hbox{ decide for any $m'\in {\bf Z}\sp{S}$ whether or not $m'$ belongs to $B$,}\
\label{(st.step2)} 
\end{equation}
though applying this more general subroutine is clearly
slower than a direct algorithm to realize \eqref{(st.step)}.

Note that  $m'=m+\chi_{s}-\chi_{t}$ is in $B$
precisely if there is no $m$-tight $t\ol s$-set (with respect to $p$),
and this is true even if $B$ is defined by a crossing supermodular function $p$.  
Subroutine \eqref{(st.step)} can be carried out by a
single application of a submod-minimizer.

As long as possible, apply the 1-tightening step 
(as described in Section \ref{SCchardecmin}).  
Recall that a 1-tightening step replaces $m$ by 
$m':=m+\chi_{s}-\chi_{t}$ 
where $s$ and $t$ are elements of $S$ for which $m(t)\geq m(s)+2$ and $m'$
belongs to $\odotZ{B}$.

By Theorem \ref{equi.1}, when no more 1-tightening step is available,
the current $m$ is a decreasingly minimal member of $\odotZ{B}$ 
and the algorithm terminates.  
In order to estimate the number of 1-tightening steps, 
observe that a single 1-tightening step decreases the
square-sum of the components.  
Since the largest square-sum of an
arbitrary integral vector $z$ with $\widetilde z(S)=p(S)$ is 
$p(S)\sp{2}$ and 
$\widetilde z(S)=p(S)$ holds for all members $z$ of $\odotZ{B}$,
we conclude that the number of 1-tightening steps is at most  $p(S)\sp{2}$.  
Therefore if $\vert p(S)\vert $ is bounded by a polynomial of $n$, 
then the basic algorithm to compute a dec-min element of $\odotZ{B}$ is strongly polynomial.

At this point, we postpone the description of the algorithm for
arbitrary $p$, and show how the canonical chain can be computed.

\subsection{Computing the essential value-sequence, the canonical chain and partition} 
\label{SCcompcanchain}

Let $B=B'(p)$ be again an integral base-polyhedron whose unique 
fully supermodular bounding function is $p$.  
In the algorithm, we must be able to compute, 
for a given (dec-min) member $m$ of $\odotZ{B}$, 
the smallest $m$-tight set $T_{m}(u)$ containing a given element $u\in S$.
Here $m$-tightness is with respect to $p$, that is, a set $X$ is $m$-tight 
if $\widetilde m(X)=p(X)$.  
It is fundamental, however, to
emphasize that $T_{m}(u)$ can be computed even in the case when $p$ is
not explicitly available and $B$ is defined by a weaker function, for
example, by a crossing supermodular function.  
Namely, recall from Claim \ref{st} that 
$T_{m}(u)=\{s:  m+\chi_{s} -\chi_{u}\in B\}$ 
and hence $T_{m}(u)$ is indeed computable by at most $n$ applications 
of routine \eqref{(st.step)}.

\begin{algorithm}  \rm  \label{algo1}
Given a dec-min element $m$ of $\odotZ{B}$, 
the following procedure computes 
the canonical chain 
${\cal C}\sp{*}=\{C_{1},C_{2},\dots ,C_{q}\}$
and the canonical partition 
${\cal P}\sp{*}=\{S_{1},S_{2}, \dots , S_{q}\}$ of $S$
along with the essential value-sequence 
$\beta_{1}>\beta_{2}>\cdots >\beta_{q}$ belonging to $\odotZ{B}$.

\begin{enumerate}
\item 
 Let $\beta_{1}$ denote the largest value of $m$.
Let $C_{1} := \bigcup \{  T_{m}(u):  m(u) =\beta_{1} \}$,
$S_{1} := C_{1}$, and $i := 2$.

\item 
In the general step $i \geq 2$, 
the pairwise disjoint non-empty sets $S_{1},S_{2},\dots ,S_{i-1}$ 
and a chain $C_{1} \subset C_{2} \subset \dots \subset C_{i-1}$
have already been computed along with the essential values 
$\beta_{1}>\beta_{2}> \dots > \beta_{i-1}$.
If $C_{i-1} = S$, set $q:= i-1$ and stop.
Otherwise, let 
\begin{align*}
\beta_{i} & := \max\{m(s):  s\in S-C_{i-1} \},
\\
C_{i}  &  := \bigcup \{ T_{m}(u):  m(u)\geq \beta_{i}\}, 
\\
S_{i}   & :=C_{i}-C_{i-1},
\end{align*}
and go to the next step for $i:=i+1$.
\finbox
\end{enumerate}
\end{algorithm}  

Corollary \ref{Cim} implies that the sequence $\beta_{1},\beta_{2},\dots ,\beta_{q}$ 
provided by this algorithm is indeed the essential value-sequence belonging to $\odotZ{B}$, 
and similarly the chain $C_{1} \subset C_{2} \subset \cdots \subset C_{q}$ is the
canonical chain while the partition $\{S_{1},S_{2}, \dots ,S_{q} \}$ is the
canonical partition.

We emphasize that the basic algorithm in Section \ref{SCbasicalg}
for computing a dec-min element $m$ of $\odotZ{B}$ 
is polynomial only in $\vert p(S)\vert $, 
meaning that it is polynomial in $n$ only if $\vert p(S)\vert $ is small 
(that is, $\vert p(S)\vert $ is bounded by a power of $n$).  
On the other hand, Algorithm \ref{algo1} to compute the essential
value-sequence and the canonical chain is strongly polynomial for
arbitrary $p$ (independently of the magnitude of $\vert p(S)\vert )$,
provided that a dec-min element $m$ of $\odotZ{B}$ is already available
as well as Oracle (7.2).

\paragraph{Adaptation to the intersection with a box} \ 
Algorithm \ref{algo1} can be adapted to the case when we have specific upper
and lower bounds on the members of $\odotZ{B}= \odotZ{B'}(p)$.  
Let
$f:S\rightarrow {\bf Z} \cup \{-\infty \}$ and 
$g:S\rightarrow {\bf Z} \cup \{+\infty \}$ 
be bounding functions with $f\leq g$ and 
let $T(f,g)$ denote the box defined by $f$ and $g$.  
It is a basic fact on integral base-polyhedra that the intersection 
$B\sq :=B\cap T(f,g)$ 
is also a (possibly empty) integral base-polyhedron.  
Assume that ${B\sq}$
is non-empty (which is, by a known theorem,
equivalent to requiring that $\widetilde f\leq \overline{p}$ and 
$\widetilde g\geq p$ when $p$ is fully supermodular.  
Here $\overline{p}$ denotes the complementary submodular function of $p$).

Let $m$ be an element of $\odotZ{B\sq}$.  
Let $T_{m}(u)$ denote the smallest $m$-tight set containing $u$ with respect to $p$, 
and let $T\sq_{m}(u)$ be the smallest $m$-tight set containing $u$ with respect to $p\sq$.

\begin{claim} \label{CLTm(u)}
\[
T\sq_{m}(u) = 
\begin{cases} \{u\} & 
\ \ \hbox{\rm if}\ \ \ m(u)=f(u),
\\ 
T_{m}(u) - \{v:  m(v) = g(v)\} & \ \ \hbox{\rm if}\ \ \ m(u)>f(u).
\end{cases}
\]
\end{claim}  

\Proof 
We have $T\sq_{m}(u)= \{ s: m-\chi_{u} + \chi_{s} \in B\sq \}$.
Since
$B\sq = B \cap T(f,g)$,
we have
$m-\chi_{u} + \chi_{s} \in B\sq$
if and only if
(i)
$m-\chi_{u} + \chi_{s} \in B$ and (ii) $m-\chi_{u} + \chi_{s} \in T(f,g)$ hold.
For $s \not = u$,  (i) holds 
if and only if
$s \in T_{m}(u)$,
and (ii) holds 
if and only if
$m(u) > f(u)$ and $m(s) < g(s)$.
Hence follows the claim.
\finbox \medskip

The claim implies that Algorithm \ref{algo1}
can be adapted easily to compute the
canonical chain and partition belonging to $\odotZ{B\sq}$
along with its essential value-sequence.

Our next goal is to describe a strongly polynomial algorithm to
compute a dec-min element of $\odotZ{B}$ in the general case when no
restriction is imposed on the magnitude of $\vert p(S)\vert $. 
To this end, we need an algorithm to maximize 
$\llceil {p(X) \over \vert X\vert} \rrceil$. 
We describe an algorithm for a more general case but this
generality will be needed only in a forthcoming paper \cite{Frank-Murota.4} 
dealing with dec-min elements of the intersection of two M-convex sets.

\subsection{Maximizing $\llceil {p(X) \over b(X)}\rrceil$ 
           with the Newton--Dinkelbach (ND) algorithm} 
\label{pbmaximizing}

On a ground-set $S$ with $n\geq 1$ elements, we are given 
an integer-valued set-function $p$ with $p(\emptyset )=0$.
  (Here $p(S)$ is finite but $p(X)$ may otherwise be
$-\infty $. However, $p(X)$ is never $+\infty $.)  
Moreover, we are also given a non-negative, finite integer-valued set-function $b$. 
Both $p$ and $b$ are integer-valued.  
Our present goal is to describe a variation of the
Newton--Dinkelbach algorithm to compute the maximum 
$\llceil {p(X) \over b(X)}\rrceil $ 
over the subsets $X$ of $S$.  
An excellent overview by Radzik \cite{Radzik} analyses this method concerning 
(among others) the special problem of maximizing 
${p(X) \over \vert X\vert }$ and
describes a strongly polynomial algorithm.  
We present a variation of the ND-algorithm whose specific feature is that it works
throughout with integers $\llceil {p(X) \over b(X)}\rrceil$. 
This has the advantage that the proof is simpler than the original one 
working with the fractions ${p(X) \over b(X)}$.

Let $M$ denote the largest value of $b$.  The algorithm works if a
subroutine is available to 
\begin{equation} 
\hbox{ find a subset of $S$ maximizing
$p(X) - \mu b(X)$ \ $(X\subseteq S)$ \ for any fixed integer $\mu \geq 0$. }\ 
\label{(ND.routine)} 
\end{equation}
This routine will actually be needed only for special values
of $\mu $ when 
$\mu =\llceil p(X)/\ell\rrceil$ 
(where $X\subseteq S$ and $1\leq \ell\leq M$). 
Note that we do not have to assume that $p$ is supermodular and $b$ is submodular, 
the only requirement for the ND-algorithm is that Subroutine \eqref{(ND.routine)} be available.  
Via a submod-minimizer, this is certainly the case when $p$ happens to be
supermodular and $b$ submodular, since then $\mu b-p$ is submodular when $\mu \geq 0$.  
But \eqref{(ND.routine)} is also available in the
more general case when the function $p'$ defined by 
$p'(X):=p(X)-\mu b(X)$ 
is only crossing supermodular.  
Indeed, for a given ordered pair of elements $s,t\in S$, 
the restriction of $p'$ on the family of $s\ol t$-sets is fully supermodular, 
and therefore we can apply a submod-minimizer 
to each of the $n(n-1)$ ordered pairs $(s,t)$ to get the requested maximum of $p'$.

We call a value $\mu $ \ {\bf good} \ if $\mu b(X)\geq p(X)$ \ 
[i.e., $p(X)-\mu b(X)\leq 0$] \ for every $X\subseteq S$.  
A value that is not good is called {\bf bad}.  
We assume that there is a good $\mu $,
which is equivalent to requiring that $p(X)\leq 0$ whenever $b(X)=0$.
We also assume that $\mu =0$ is bad.

Our goal is to compute the minimum $\mu_{\rm min}$ of the good integers.  
In other words, we want to maximize \ $\llceil {p(X) \over b(X)}\rrceil $ \ 
over the subsets of $S$ with $b(X)>0$.

The algorithm starts with the bad $\mu_{0}:=0$.  
Let 
\[
  X_{0} \in \arg\max \{ p(X)-\mu_{0}b(X) :  \ X\subseteq S \},
\]
that is, \ $X_{0}$ is a set maximizing the function $p(X)-\mu_{0}b(X)=p(X)$.  
Note that the badness of $\mu_{0}$ implies that
$p(X_{0}) > 0$.  Since, by the assumption, there is a good $\mu $, 
it follows that $\mu b(X_{0}) \geq p(X_{0}),$ and hence $b(X_{0})>0$.

The procedure determines one by one a series of pairs $(\mu_{j},X_{j})$
for subscripts $j=1,2,\dots $ where each integer $\mu_{j}$ 
is an (intermediate) tentative candidate 
for $\mu $ while $X_{j}$ is a non-empty subset of $S$.  
Suppose that the pair $(\mu_{j-1},X_{j-1})$
has already been determined for a subscript $j\geq 1$.  
Let $\mu_{j}$ be the smallest integer 
for which $\mu_{j}b(X_{j-1})\geq p(X_{j-1})$,
that is,
\[
 \mu_{j}:  =\llceil {p(X_{j-1}) \over b(X_{j-1})}\rrceil . 
\]

If $\mu_{j}$ is bad, that is, if there is a set $X\subseteq S$ with
$p(X) -\mu_{j}b(X) > 0$, then let 
\[
  X_{j} \in \arg\max \{ p(X)-\mu_{j}b(X) :  \ X\subseteq S \},
\] 
that is, \ $X_{j}$ is a set
maximizing the function \ $p(X)-\mu_{j}b(X)$.  
(If there are more than one maximizing set, we can take any).  
Since $\mu_{j}$ is bad, $X_{j}\not =\emptyset $ and $p(X_{j}) - \mu_{j}b(X_{j})>0$.

\begin{claim}  \label{betano} 
If $\mu_{j}$ is bad for some subscript $j\geq 0$, then $\mu_{j} < \mu_{j+1}$.  
\end{claim}

\Proof 
The badness of $\mu_{j}$ means that 
$p(X_{j})-\mu_{j}b(X_{j}) > 0$ from which
\[
\mu_{j+1} = \llceil {p(X_{j}) \over b(X_{j})}\rrceil 
  = \llceil {p(X_{j})-\mu_{j}b(X_{j}) \over b(X_{j}) }\rrceil + \mu_{j} 
\ > \  \mu_{j}.  
\]
\vspace{-2.8\baselineskip} \\
\finbox
\vspace{1.5\baselineskip}

Since there is a good $\mu $ and the sequence $\mu_{j}$ is 
strictly monotone increasing by Claim \ref{betano}, 
there will be a first subscript $h\geq 1$ for which $\mu_{h}$ is good.  
The algorithm terminates by outputting this $\mu_{h}$ 
(and in this case $X_{h}$ is not computed or needed anymore).

\begin{theorem}   \label{Msteps} 
If $h$ is the first subscript during the run of the algorithm 
for which $\mu_{h}$ is good, then 
$\mu_{\rm min}=\mu_{h}$ 
(that is, $\mu_{h}$ is the requested smallest good $\mu $-value)
and $h\leq M$, where $M$ denotes the largest value of $b$.  
\end{theorem}

\Proof 
Since $\mu_{h}$ is good and $\mu_{h}$ is the smallest integer for
which $\mu_{h}b(X_{h-1})\geq p(X_{h-1})$, 
the set $X_{h-1}$ certifies that no good integer $\mu $ 
can exist which is smaller than $\mu_{h}$,
that is, $\mu_{\rm min}=\mu_{h}$.

\begin{proposition}  \label{Xno} 
If $\mu_{j}$ is bad for some subscript $j\geq 1$, then $b(X_{j}) < b(X_{j-1})$.  
\end{proposition}

\Proof
As $\mu_{j}$ \ ($= \llceil {p(X_{j-1})\over b(X_{j-1})}\rrceil$) is bad, 
we obtain that
\begin{align*}
&
p(X_{j})-\mu_{j}b(X_{j}) >0 = p(X_{j-1}) - { p(X_{j-1}) \over b(X_{j-1})} b(X_{j-1}) 
\\ &
\geq p(X_{j-1}) - \llceil {p(X_{j-1}) \over b(X_{j-1})}\rrceil b(X_{j-1})
= p(X_{j-1}) - \mu_{j}b(X_{j-1}) ,
\end{align*}
from which we get
\[ 
\hbox{(A)}\quad p(X_{j}) - \mu_{j}b(X_{j}) > p(X_{j-1}) - \mu_{j}b(X_{j-1}).
\]
Since $X_{j-1}$ maximizes $p(X) - \mu_{j-1}b(X)$, it follows that
\[ \quad \quad
\hbox{(B)}\quad p(X_{j-1}) - \mu_{j-1}b(X_{j-1}) 
\geq p(X_{j}) - \mu_{j-1}b(X_{j}).  
\]
By adding up (A) and (B), we obtain
\[
 (\mu_{j} - \mu_{j-1})b(X_{j-1}) > (\mu_{j} - \mu_{j-1})b(X_{j}).
\]
As $\mu_{j}$ is bad, so is $\mu_{j-1}$, and hence, 
by applying Claim \ref{betano} to $j-1$ in place of $j$, we obtain that
$\mu_{j} > \mu_{j-1}$, from which we arrive at $b(X_{j}) < b(X_{j-1})$, 
as required.  
\finbox

\medskip

Proposition \ref{Xno} implies that 
$M \geq b(X_{0}) > b(X_{1})> \cdots >b(X_{h-1})$,
from which $1\leq b(X_{h-1}) \leq M -(h-1)$, and hence $h\leq M$ follows.  
\finbox \finboxHere 
\medskip

Note that the ND-algorithm is definitely polynomial in the special
case when $M$ is bounded by a power of $n$, since the number of
phases is bounded by $M$.  
In the special case when $b(X)=\vert X\vert$, $M = n$.

In several applications, the requested general purpose
submod-minimizer can be superseded by a direct and more efficient
algorithm such as the ones for network flows or for matroid partition.

\subsection{Computing a dec-min element in strongly polynomial time}
\label{strong.pol}

In the present context, we need the ND-algorithm above only in the
special case when $b$ is the cardinality function, that is,
$b(X)=\vert X\vert $ for each $X\subseteq S$.  
Note that in this special case we have 
$M=\vert S\vert $, and hence the sequence of bad $\mu_{i}$
values has at most $\vert S\vert $ members by Theorem \ref{Msteps}.

After at most $\vert S\vert $ applications of Subroutine \eqref{(ND.routine)}, 
the ND-algorithm terminates with the smallest integer $\beta_{1}$ 
for which $\odotZ{B}$ has a $\beta_{1}$-covered member $m$.  
It is well-known that such an $m$ can easily be computed
with a greedy-type algorithm, as follows.  
Since there is a $\beta_{1}$-covered member of $B$, the vector 
$(\beta_{1},\beta_{1},\dots ,\beta_{1})$ 
belongs to the so-called supermodular polyhedron 
$S'(p):=\{x:  \widetilde x(X)\geq p(X)$ for every $X\subseteq S\}$.  
Consider the elements of $S$ in an arbitrary order $\{s_{1},\dots ,s_{n}\}$.  
Let $m(s_{1}):=\min \{z:  (z,\beta_{1},\beta_{1},\dots ,\beta_{1})\in S'(p)\}$. 
In the general step, if the components 
$m(s_{1}),\dots ,m(s_{i-1})$ have already been determined, let
\begin{equation} 
m(s_{i}):= \min \{z: 
 (m(s_{1}),m(s_{2}),\dots ,m(s_{i-1}), z,\beta_{1},\beta_{1},\dots ,\beta_{1})\in S'(p)\}.  
\label{(m.def)} 
\end{equation}

This computation can be carried out by $n$ applications of a subroutine
for a submodular function minimization.  
(Note that the previous algorithm to compute a $\beta_{1}$-covered integral member of $B$ 
is nothing but a special case of the algorithm that finds an integral member 
of a base-polyhedron given by an intersecting submodular function.)

Given a $\beta_{1}$-covered integral element of $B$, our next goal is
to obtain a pre-dec-min element of $\odotZ{B}$.  
To this end, we apply 1-tightening steps.  
That is, as long as possible, we pick two elements $s$ and $t$  of $S$ 
for which $m(t)=\beta_{1}$ and  $m(s)\leq \beta_{1} -2$ 
such that there is no $m$-tight $t\ol s$-set, reduce $m(t)$ by 1 
and increase $m(s)$ by 1. 
In this way, we obtain another integral element of $B$ 
for which the largest component continues to be $\beta_{1}$ 
(as $\beta_{1}$ was chosen to be the smallest upper bound) but the
number of $\beta_{1}$-valued components is strictly smaller.
Therefore, after at most $\vert S\vert -1$ such 1-tightening steps, we
arrive at a vector for which no 1-tightening step 
(with $m(t)=\beta_{1}$ and $m(s)\leq \beta_{1} -2$) 
is possible anymore, and hence this final vector 
is a pre-decreasingly minimal element of $\odotZ{B}$ by Theorem \ref{optkrit}.  
We use the same letter $m$ to denote this pre-dec-min element.

Recall that $T_{m}(t)$ denoted the unique smallest tight set containing $t$ 
when $p$ is fully supermodular.  
But $T_{m}(t)$ can be described without explicitly referring to $p$ 
since an element $s\in S$ belongs to $T_{m}(t)$ precisely if 
$m':=m - \chi_t + \chi_{s}$ is in $B$, 
and this is computable by subroutine \eqref{(st.step)}.  
Therefore we can compute $S_{1}(m)$ (defined in \eqref{(S1m.def)}).  
It was proved in Theorem \ref{smallest} that $S_{1}(m)$ 
is the first member $S_{1}$ of the canonical $S$-partition associated with $\odotZ{B}$.

The restriction $m_{1}:=m\vert S_{1}$ is a near-uniform member of the
restriction of $\odotZ{B}$ to $S_{1}$, and by Theorem \ref{felbont}, 
if $m_{1}'$ is a dec-min element of $\odotZ{B_{1}'}$, 
then $(m_{1},m_{1}')$ is a dec-min element of $\odotZ{B}$, 
where $B_{1}'$ is the base-polyhedron obtained from $B$ by contracting $S_{1}$.  
Such a dec-min element $m_{1}'$ can be computed 
by applying iteratively the computation described above for computing $m_{1}$.



\section{Applications} 
\label{alk}

\subsection{Background}

There are two major sources of applicability of the results in the
preceding sections.  
One of them relies on the fact that the class of
integral base-polyhedra is closed under several operations.  
For example, a face of a base-polyhedron is also a base-polyhedron, and so
is the intersection of an integral box with a base-polyhedron $B$.
Also, the sum of integral base-polyhedra $B_{1},\dots ,B_{k}$ is a
base-polyhedron $B$ which has, in addition, the integer decomposition
property meaning that any integral element of $B$ can be obtained as
the sum of $k$ integral elements by taking one from each $B_{i}$.  
This latter property implies that the sum of M-convex sets is M-convex.  
We also mention the important operation of taking an aggregate of a
base-polyhedron, to be introduced below in Section \ref{matroidalk}.

The other source of applicability is based on the fact that not only
fully super- or submodular functions can define base-polyhedra but
some weaker functions as well.  
For example, if $p$ is an integer-valued crossing 
(in particular, intersecting) supermodular
function with finite $p(S)$, then $B=B'(p)$ is a (possibly empty)
integral base-polyhedron (and $\odotZ{B}$ is an M-convex set).  
This fact will be exploited in solving dec-min orientation problems when
both degree-constraints and edge-connectivity requirements must be fulfilled.  
In some cases even weaker set-functions can define base-polyhedra.  
This is why we can solve dec-min problems concerning
edge- and node-connectivity augmentations of digraphs.

\subsection{Matroids} 
\label{matroidalk}

Levin and Onn \cite{Levin-Onn} solved algorithmically the following problem.  
Find $k$ bases of a matroid $M$ on a ground-set $S$ such
that the sum of their characteristic vectors be decreasingly minimal.
Their approach, however, does not seem to work in the following natural extension.  
Suppose we are given $k$ matroids $M_{1},\dots ,M_{k}$
on a common ground-set $S$, and our goal is to find a basis $B_{i}$ of
each matroid $M_{i}$ in such a way that the vector  
$\sum [\chi_{B_{i}} : i=1,\dots ,k]$ 
is decreasingly minimal.  
Let $B_{\sum }$ denote the sum of the base-polyhedra of the $k$ matroids.  
By a theorem of Edmonds, the integral elements of $B_{\sum }$ 
are exactly the vectors of form 
$\sum [\chi_{B_{i}} :  i=1,\dots ,k]$ 
where $B_{i}$ is a basis of $M_{i}$.  
Therefore the problem is to find a dec-min element of $\odottZ{B_{\sum }}$.  
This can be found by the basic algorithm described in Section \ref{SCbasicalg}.
Let us see how the requested subroutines are available in this special case.  
The algorithm starts with an arbitrary member $m$ of $\odottZ{B_{\sum }}$ 
which is obtained by taking a basis $B_{i}$ from each matroid $M_{i}$, 
and these bases define $m:=\sum_{i} \chi_{B_{i}}$.

To realize Subroutine \eqref{(st.step)}, we mentioned that it suffices
to realize Subroutine \eqref{(st.step2)}, which requires for a given
integral vector $m'$ with $\widetilde m'(S)= \sum_{i}r_{i}(S)$ to decide
whether $m'$ is in $\odottZ{B_{\sum }}$ or not.  
But this can simply be done by Edmonds' matroid intersection algorithm.  
Namely, let
$S_{1},\dots ,S_{k}$ be disjoint copies of $S$ and $M_{i}'$ an isomorphic
copy of $M_{i}$ on $S_{i}$.  
Let $N_{1}$ be the direct sum of matroids
$M_{i}'$ on ground-set $S':=S_{1}\cup \cdots \cup S_{k}$.  
Let $N_{2}$ be a partition matroid on $S'$ in which a subset $Z$ is a basis 
if it contains exactly $m'(s)$ members of the $k$ copies of $s$ 
for each $s\in S$.  
Then $m'$ is in $\odottZ{B_{\sum }}$ precisely if $N_{1}$ and
$N_{2}$ have a common basis.

In conclusion, with the help of Edmonds' matroid intersection algorithm, 
Subroutine \eqref{(st.step)} is available, and hence the
basic algorithm can be applied.  
(Actually, the algorithm can be sped up by looking into the details 
of the matroid intersection algorithm for $N_{1}$ and $N_{2}$.)

\medskip \medskip

Another natural problem concerns a single matroid $M$ on a ground-set $T$.  
Suppose we are given a partition ${\cal P}= \{T_{1},\dots ,T_{n}\}$
of $T$ and we consider the intersection vector 
$(\vert Z\cap T_{1}\vert,\dots ,\vert Z\cap T_{n}\vert )$ 
assigned to a basis $Z$ of $M$.  
The problem is to find a basis for which the intersection vector is decreasingly minimal.

To solve this problem, we recall an important construction of
base-polyhedra, called the aggregate.  
Let $T$ be a ground-set and $B_{T}$ an integral base-polyhedron in ${\bf R}\sp T$.  
Let ${\cal P}=\{T_{1},\dots ,T_{n}\}$ be a partition of $T$ into non-empty subsets
and let $S=\{s_{1},\dots ,s_{n}\}$ be a set 
whose elements correspond to the members of ${\cal P}$.  
The aggregate $B_{S}$ of $B_{T}$ is defined as follows.
\begin{equation} 
\hbox{ $B_{S}:= \{ (y_{1},\dots ,y_{n}):$ there is an $x\in B_{T}$ with
   $y_{i}=\widetilde x(T_{i}) \ (i=1,\dots ,n)\}$.   }\ 
\label{(aggre)} 
\end{equation}
A basic theorem concerning base-polyhedra states that $B_{S}$ is a
base-polyhedron, moreover, for each integral member $(y_{1},\dots ,y_{n})$
of $B_{S}$, the vector $x$ in \eqref{(aggre)} can be chosen integer-valued.  
In other words,
\begin{equation} 
\hbox{ $\odotZ{B_{S}}:= \{ (y_{1},\dots ,y_{n}):$ there is an $x\in
\odotZ{B_{T}}$ with $y_{i}=\widetilde x(T_{i}) \ (i=1,\dots ,n)\}$.   }\
\label{(aggre.odot)} 
\end{equation}
We call $\odotZ{B_{S}}$ the {\bf aggregate} of $\odotZ{B_{T}}$.

Returning to our matroid problem, let $B_{T}$ denote the base-polyhedron
of matroid $M$.  Then the problem is nothing but finding a dec-min
element of $\odottZ{B_{S}}$.

We can apply the basic algorithm (concerning M-convex sets) for this
special case since the requested subroutines are available through
standard matroid algorithms.  Namely, Subroutine \eqref{(member)} is
available since for any basis $Z$ of $M$, the intersection vector
assigned to $Z$ is nothing but an element of $\odotZ{B_{S}}$.

To realize Subroutine \eqref{(st.step)}, we mentioned that it suffices
to realize Subroutine \eqref{(st.step2)}.  
Suppose we are given a vector $y\in {\bf Z}_{+}\sp{S}$ 
(Here $y$ stands for $m'$ in \eqref{(st.step2)}).  
Suppose that $\widetilde y(S)=r(T)$
 (where $r$ is the rank-function of matroid $M$) and that 
$y(s_{i})\leq \vert T_{i}\vert$ for $i=1,\dots ,n$.

Let $G=(S,T;E)$ denote a bipartite graph where $E=\{ts_{i}:  t\in T_{i},
i=1,\dots ,n\}$.  By this definition, the degree of every node in $T$
is 1 and hence the elements of $E$ correspond to the elements of $M$.
Let $M_{1}$ be the matroid on $E$ corresponding to $M$ (on $T$).  Let
$M_{2}$ be a partition matroid on $E$ in which a set $F\subseteq E$ is a
basis if $d_{F}(s_{i})=y(s_{i})$.  
By this construction, the vector $y$ is
in $\odotZ{B_{S}}$ precisely if the two matroids $M_{1}$ and $M_{2}$ have a
common basis.  This problem is again tractable 
by Edmonds' matroid intersection algorithm.

\medskip 
As a special case, we can find a spanning tree of a
(connected) directed graph for which its in-degree-vector is
decreasingly minimal.  Since the family of unions of $k$ disjoint
bases of a matroid forms also a matroid, we can also compute $k$
edge-disjoint spanning trees in a digraph whose union has a
decreasingly minimal in-degree vector.

Another special case is when we want to find a spanning tree of a
connected bipartite graph $G=(S,T;E)$ whose in-degree vector
restricted to $S$ is decreasingly minimal.

\subsection{Flows} 
\label{SCflows}

\subsubsection{A base polyhedron associated with net-in-flows}
\label{SCflowbase}

Let $D=(V,A)$ be a digraph endowed with integer-valued bounding functions 
$f:A\rightarrow {\bf Z}\cup \{-\infty \}$ and
$g:A\rightarrow {\bf Z}\cup \{+\infty \}$ for which $f\leq g$.  
We call a vector (or function) $z$ on $A$ {\bf feasible} 
if $f\leq z\leq g$.  
The {\bf net-in-flow} $\Psi_z$ of $z$ is a vector on $V$ and defined by 
$\Psi_z(v) = \varrho_z(v) -\delta_z(v)$, 
where 
$\varrho_z(v):=\sum [z(uv):  uv\in A]$ and $\delta_z(v):=\sum [z(vu):  uv\in A]$.  
If $m$ is the net-in-flow of a vector $z$, 
then we also say that $z$ is an {\bf $m$-flow}.

A variation of Hoffman's classic theorem on feasible circulations
\cite{Hoffman60} is as follows.

\begin{lemma} 
\label{LMnetinflowbase}
An integral vector $m:V\rightarrow {\bf Z} $ is the net-in-flow
of an integral feasible vector 
(or in other words, there is an integer-valued feasible $m$-flow) 
if and only if 
$\widetilde m(V)=0$ and
\begin{equation} 
\varrho_f(Z) - \delta_g(Z) \leq \widetilde m(Z)
 \ \ \hbox{\rm holds whenever}\ \ \ Z\subseteq V, 
\label{(feasmflow)} 
\end{equation}
where \ 
$\varrho_f(Z):=\sum [f(a):  a\in A \hbox{\  \rm and $a$ enters $Z$} ]$ 
\ and \ 
$\delta_g(Z):=\sum [g(a):  a\in A \hbox{\  \rm and $a$ leaves $Z$}]$.
\finbox 
\end{lemma}

Define a set-function $p_{fg}$ on $V$ by 
\[
  p_{fg}(Z):=\varrho_f(Z)-\delta_g(Z).
\] 
Then $p_{fg}$ is (fully) supermodular 
(see, e.g.  Proposition 1.2.3 in \cite{Frank-book}).  
Consider the base-polyhedron
$B_{fg}:=B'(p_{fg})$ and the M-convex set $\odottZ{B_{fg}}$.  
By Lemma \ref{LMnetinflowbase} 
the M-convex set $\odottZ{B_{fg}}$ consists
exactly of the net-in-flow integral vectors $m$.

By the algorithm described in Section \ref{SCalgo2}, we can compute a
decreasingly minimal element of $\odottZ{B_{fg}}$ in strongly polynomial time.  
By relying on a strongly polynomial push-relabel algorithm, 
we can check whether or not \eqref{(feasmflow)} holds.  
If it does not, then the push-relabel algorithm can compute a set most
violating \eqref{(feasmflow)} 
(that is, a maximizer of $\varrho_f(Z) - \delta_g(Z) - \widetilde m(Z)$) 
while if \eqref{(feasmflow)} does hold, 
then the push-relabel algorithm computes an integral valued feasible $m$-flow.  
Therefore the requested oracles in the general algorithm 
for computing a dec-min element are available through a
network flow algorithm, and we do not have to rely on a
general-purpose submodular function minimizing oracle.

For the sake of an application of this algorithm to capacitated
dec-min orientations in Section \ref{capori}, we remark that the
algorithm can also be used to compute a dec-min element of the
M-convex set obtained from $\odottZ{B_{fg}}$ 
by translating it with a given integral vector.

\subsubsection{Discrete version of Megiddo's flow problem}

Megiddo \cite{Megiddo74}, \cite{Megiddo77} considered the following problem.  
Let $D=(V,A)$ be a digraph endowed with a non-negative
capacity function $g:A\rightarrow {\bf R}+$.  
Let $S$ and $T$ be two disjoint non-empty subsets of $V$.  
Megiddo described an algorithm to compute a feasible flow from $S$ to $T$ 
with maximum flow amount $M$
for which the net-in-flow vector restricted on $S$ 
is (in our terms) increasingly maximal.  
Here a feasible flow is a vector $x$ on $A$ 
for which $\Psi_{x}(v) \leq 0$ for $v\in S$, 
$\Psi_{x}(v)\geq 0$ for $v\in T$,
and $\Psi_{x}(v)=0$ for $v\in V-(S\cup T)$.  
The flow amount $x$ is $\sum [\Psi_{x}(t):t\in T]$.

We emphasize that Megiddo solved the continuous (fractional) case and
did not consider the corresponding discrete (or integer-valued) flow problem.  
To our knowledge, this natural optimization problem has not been investigated so far.

To provide a solution, suppose that $g$ is integer-valued.  
Let $f\equiv 0$ and consider the net-in-flow vectors belonging to feasible vectors.  
These form a base-polyhedron $B_{1}$ in ${\bf R}\sp V$.  
Let $B_{2}$ denote the base polyhedron obtained from $B_{1}$ by intersecting
it with the box defined by 
$z(v)\leq 0$ for $v\in S$, $z(v)\geq 0$ for $v\in T$ 
and $z(v)=0$ for $v\in V-(S\cup T)$.   

The restriction of $B_{2}$ to $S$ is a g-polymatroid $Q$ in ${\bf R}\sp{S}$.  
And finally, we can consider the face of $Q$ defined by $\widetilde z(S)=-M$.  
This is a base-polyhedron $B_3$ in ${\bf R}\sp{S}$, 
and the discrete version of Megiddo's flow problem is equivalent 
to finding an inc-max element of $\odottZ{B_3}$. 
(Recall that an element of an M-convex set is dec-min precisely if it is inc-max.)

It can be shown that in this case again the general submodular
function minimizing subroutine used in the algorithm to find a dec-min
element of an M-convex set can be replaced by a max-flow min-cut algorithm.

In Part III \cite{Frank-Murota.3}
we solve a more general discrete dec-min problem,
in which we are to find an integral feasible flow
that is dec-min on an arbitrarily specified edge set.

\subsection{Further applications}

\subsubsection{Root-vectors of arborescences}

A graph-example comes from packing arborescences.  
Let $D=(V,A)$ be a digraph and $k>0$ an integer.  
We say that a non-negative integral vector 
$m:V\rightarrow {\bf Z}_{+}$
is a {\bf root-vector} if
there are $k$ edge-disjoint spanning arborescences such that each node
$v\in V$ is the root of $m(v)$ arborescences.  
Edmonds \cite{Edmonds73} classic result on disjoint arborescences implies that
$m$ is a root-vector 
if and only if 
$\widetilde m(V)=k$ and
$\widetilde m(X)\geq k-\varrho(X)$ holds for every subset $X$ with
$\emptyset \subset X\subset V$.  
Define set-function $p$ by $p(X):  = k-\varrho(X)$ 
if $\emptyset \subset X\subseteq V$ and $p(\emptyset):=0$.  
Then $p$ is intersecting supermodular, so $B'(p)$ is an integral base-polyhedron.  
The intersection $B$ of $B'(p)$ with the
non-negative orthant is also a base-polyhedron, and the theorem of
Edmonds is equivalent to stating that a vector $m$ is a root-vector 
if and only if 
$m$ is in $\odotZ{B}$.

Therefore the general results on base-polyhedra can be specialized to
obtain $k$ disjoint spanning arborescences whose root-vector is
decreasingly minimal.

\subsubsection{Connectivity augmentations}

Let $D=(V,A)$ be a directed graph and $k>0$ an integer.  
We are interested in finding a so-called augmenting digraph 
$H=(V,F)$ of $\gamma $ arcs for which $D+H$ is $k$-edge-connected or
$k$-node-connected.  
In both cases, the in-degree vectors of the
augmenting digraphs are the integral elements of an integral
base-polyhedron \cite{FrankJ23}, \cite{FrankJ31}.  
Obviously, the in-degree vectors of the augmented digraphs 
are the integral elements of an integral base-polyhedron.

Again, our results on general base-polyhedra can be specialized to
find an augmenting digraph whose in-degree vector is decreasingly minimal.



\section{Orientations of graphs} 
\label{SCori1}

Let $G=(V,E)$ be an undirected graph.  
For $X\subseteq V$, let $i_{G}(X)$
denote the number of edges induced by $X$
while $e_{G}(X)$
is the number of edges with at least one end-node in $X$.  
Then $i_{G}$ is supermodular, $e_{G}$ is submodular, and they are
complementary functions, that is, 
$i_{G}(X) = e_{G}(V)- e_{G}(V-X)$.  
Let $B_{G}:= B(e_{G}) = B'(i_{G})$ 
denote the base-polyhedron defined by $e_{G}$ or $i_{G}$.

We say that a function $m:V\rightarrow {\bf Z} $ is 
the {\bf in-degree vector} of an orientation $D$ of $G$ 
if $\varrho_{D}(v)=m(v)$ for each node $v\in V$.  
An in-degree vector $m$ obviously meets the equality
$\widetilde m(V)=\vert E\vert $. 
The following basic result, sometimes called the Orientation lemma, 
is due to Hakimi \cite{Hakimi}.

\begin{lemma}[Orientation lemma] \label{LMorientationlemma}
Let $G=(V,E)$ be an undirected graph and 
$m: V\rightarrow {\bf Z}$ an integral vector for which 
$\widetilde m(V)=\vert E\vert$. 
Then $G$ has an orientation with in-degree vector $m$ 
if and only if
\begin{equation} 
\widetilde m(X) \leq e_{G}(X)\ \ \hbox{\rm for every subset \ $X\subseteq V$,}\ 
\label{(eG)} 
\end{equation}
which is equivalent to
\begin{equation} 
\widetilde m(X) \geq i_{G}(X) \ \  \hbox{\rm for every subset \ $X\subseteq V$.}\  
\label{(iG)} 
\end{equation}
\vspace{-1\baselineskip}\\
\finbox
\end{lemma}

This immediately implies the following claim.

\begin{claim} \label{oribase} 
The in-degree vectors of orientations of $G$
are precisely the integral elements of base-polyhedron 
$B_{G}$\ $(=B(e_{G}) = B'(i_{G}))$, 
that is, the set of in-degree vectors of orientations of $G$ 
is the M-convex set $\odottZ{B_{G}}$.  
\finbox
\end{claim}

The proof of Lemma \ref{LMorientationlemma} is algorithmic 
(see, e.g.,  Theorem 2.3.2 of \cite{Frank-book})
and the orientation corresponding to a given $m$ can be constructed easily.

\subsection{Decreasingly minimal orientations}
\label{SCdecminori}

Due to Claim \ref{oribase}, we can apply the earlier results on
dec-min elements to the special base-polyhedron $B_{G}$.  
Recall that Borradaile et al.~\cite{BIMOWZ} called an orientation of $G$
egalitarian if its in-degree vector is decreasingly minimal 
but we prefer the term {\bf dec-min orientation} 
since an orientation with an increasingly maximal in-degree vector 
also has an intuitive
egalitarian feeling.  Such an orientation is called {\bf inc-max}.
For example, Theorem \ref{equi.1} immediately implies the following.

\begin{corollary}  \label{dm-im} 
An orientation of $G$ is dec-min 
if and only if 
it is inc-max.  
\finbox 
\end{corollary}

Note that the term dec-min orientation is asymmetric in the sense that
it refers to in-degree vectors.  
One could also aspire for finding an orientation 
whose out-degree vector is decreasingly minimal.  
But this problem is clearly equivalent to the in-degree version and hence in
the present work we do not consider out-degree vectors, 
with a single exception in Section \ref{inout}.

By Theorem \ref{equi.1}, an element $m$ of $\odottZ{B_{G}}$ 
is decreasingly minimal 
if and only if 
there is no 1-tightening step for $m$.  
What is the meaning of a 1-tightening step in terms of orientations?

\begin{claim} \label{reori} 
Let $D$ be an orientation of $G$ with in-degree vector $m$.  
Let $t$ and $s$ be nodes of $G$.  
The vector $m':=m + \chi_s - \chi_t$ is in $B_{G}$ 
if and only if 
$D$ admits a dipath from $s$ to $t$.  
\end{claim}  

\Proof 
$m'\in B_{G}$ holds precisely if there is no $t\ol s$-set $X$
which is tight with respect to $i_{G}$, 
that is, $\widetilde m(X) = i_{G}(X)$.  
Since 
$\varrho(Y) + i_{G}(Y) = \sum [\varrho(v):v\in Y]  = \widetilde m(Y)$
holds for any set $Y\subseteq V$, 
the tightness of $X$ is equivalent to requiring that 
$\varrho(X)=0$.  
Therefore $m'\in B_{G}$ 
if and only if 
$\varrho(Y) > 0$ holds for every $t\ol s$-set $Y$, 
which is equivalent to the existence of a dipath of $D$ from $s$ to $t$.  
\finbox 
\medskip

Recall that a 1-tightening step at a member $m$ of $B_{G}$ consists of
replacing $m$ by $m'$ provided that $m(s)\geq m(t)+2$ and $m'\in B_{G}$.
By Claim \ref{reori}, a 1-tightening step at a given orientation of
$G$ corresponds to reorienting an arbitrary dipath 
from a node $s$ to node $t$  for which 
$\varrho(s)\geq \varrho(t)+2$.  
Therefore Theorem \ref{equi.1} immediately implies 
the following basic theorem of Borradaile et al.~\cite{BIMOWZ}.

\begin{theorem}[Borradaile et al.~\cite{BIMOWZ}] \label{pathrev} 
An orientation $D$ of a graph $G=(V,E)$ is decreasingly minimal 
if and only if 
no dipath exists from a node $s$ to a node $t$ 
for which $\varrho(t)\geq \varrho(s)+2$.  
\finbox 
\end{theorem}  

Note that this theorem also implies Corollary \ref{dm-im}.  
It immediately gives rise to an algorithm for finding a dec-min orientation.  
Namely, we start with an arbitrary orientation of $G$.
We call a dipath {\bf feasible} if $\varrho(t)\geq \varrho(s)+2$
holds for its starting node $s$ and end-node $t$.  
The algorithm consists of reversing feasible dipaths as long as possible.  
Since the sum of the squares of in-degrees always drops 
when a feasible dipath is reversed, 
and originally this sum is at most $\vert E\vert \sp{2}$,
the dipath-reversing procedure terminates after at most 
$\vert E\vert \sp{2}$ reversals.  
By Theorem \ref{pathrev}, when no more feasible dipath exists, 
the current orientation is dec-min.  
The basic algorithm concerning general base-polyhedra in 
Section~\ref{SCbasicalg}
is nothing but an extension of the algorithm of Borradaile et al.

It should be noted that they suggested to choose at every step the
current feasible dipath in such a way that the in-degree of its
end-node $t$ is as high as possible, and they proved that the
algorithm in this case terminates after 
at most $\vert E\vert \vert V\vert $ dipath reversals.  

Note that we obtained Corollary \ref{dm-im} as a special case of a
result on M-convex sets but it is also a direct consequence 
of Theorem \ref{pathrev}.

\medskip

\subsection{Capacitated orientation} 
\label{capori}

Consider the following capacitated version 
of the basic dec-min orientation problem of Borradaile et al.~\cite{BIMOWZ}.  
Suppose that a positive integer $\ell(e)$ 
is assigned to each edge $e$ of $G$.  
Denote by $G\sp +$ the graph arising from $G$ 
by replacing each edge $e$ of $G$ with $\ell(e)$ parallel edges.  
Our goal is to find a dec-min orientation of $G\sp +$.  
In this case, an orientation of $G\sp +$ is described 
by telling that, among the $\ell(e)$ parallel
edges connecting the end-nodes $u$ and $v$ of $e$ 
how many are oriented toward $v$ 
(implying that the rest of the $\ell(e)$ edges are oriented toward $u$).  
In principle, this problem can be solved by applying the algorithm 
described above to $G\sp +$, and this algorithm
is satisfactory when $\ell$ is small in the sense that its largest
value can be bounded by a power of $\vert E\vert $. 
The difficulty in the general case is that the algorithm will be polynomial 
only in the number of edges of $G\sp +$, 
that is, in $\widetilde \ell(E)$, 
and hence this algorithm is not polynomial in $\vert E\vert $.

We show how the algorithm in 
Section~\ref{SCflowbase}
can be used to solve
the decreasingly minimal orientation problem in the capacitated case
in strongly polynomial time.  
To this end, let $D=(V,A)$ be an arbitrary orientation of $G$ 
serving as a reference orientation.
Define a capacity function $g$ on $A$ by $g(\ora e):=\ell(e)$,
where $\ora e$ denotes the arc of $D$ obtained by orienting $e$.

We associate an orientation of $G\sp +$ with an integral vector 
$z:A\rightarrow {\bf Z}_{+}$
with $z\leq g$ as follows.  For an arc $uv$ of $D$, 
orient $z(uv)$ parallel copies of $e=uv\in E$ toward $v$ and
$g(uv)-z(uv)$ parallel copies toward $u$.  
Then the in-degree of a node $v$ is 
$m_z(v):= \varrho_z(v) + \delta_{g-z}(v) 
  = \varrho_z(v) - \delta_z(v) +\delta_g(v)$.
Therefore our goal is to find an integral vector $z$ on $A$ 
for which $0\leq z\leq g$ and the vector $m_z$ on $V$ is dec-min.  
Consider the set of net-in-flow vectors
$\{(\Psi_z(v):  v\in V) :  0\leq z\leq g\}$.  
In Section~\ref{SCflowbase},
we proved that this is a base-polyhedron $B_{1}$.  
Therefore the set of vectors $(m_z(v):v\in V)$ is also a base-polyhedron $B$ arising from
$B_{1}$ by translating $B_{1}$ with the vector $(\delta_g(v):v\in V)$.

As remarked at the end of 
Section~\ref{SCflowbase},
a dec-min element of
$\odottZ{B}$ can be computed in strongly polynomial time by relying on
a push-relabel subroutine for network flows 
(and not using a general-purpose submodular function minimizer).

\subsection{Canonical chain and essential value-sequence for orientations}
\label{SCoricano}

In Section \ref{SCcompcanchain}, we described Algorithm \ref{algo1} for an arbitrary
M-convex set $\odotZ{B}$ that computes, from a given dec-min element
$m$ of $\odotZ{B}$, the canonical chain and essential value-sequence
belonging to $\odotZ{B}$.  
That algorithm needed an oracle for computing 
the smallest $m$-tight set $T_{m}(u)$ containing $u$.  
Here we show how this general algorithm can be turned into 
a pure graph-algorithm in the special case of dec-min orientations.

To this end, consider the special M-convex set, denoted by $\odotZ{B_{G}}$, 
consisting of the in-degree vectors of the orientations 
of an undirected graph $G=(V,E)$.  
By the Orientation lemma, 
$B_{G}=B'(i_{G})$ where $i_{G}(X)$ denotes the number of edges induced by $X$.  
Recall that $i_{G}$ is a fully supermodular function.  
For an orientation $D$ of $G$ with in-degree vector $m$, the smallest
$m$-tight set $T_{m}(t)$ (with respect to $i_{G}$) containing a node $t$
will be denoted by $T_{D}(t)$.

\begin{claim} \label{mtight} 
Let $D$ be an arbitrary orientation of $G$ with in-degree vector $m$.  
{\rm (A)}\  A set $X\subseteq V$ is $m$-tight (with respect to $i_{G}$) 
if and only if $\varrho_{D}(X)=0$. 
{\rm (B)} \ The smallest $m$-tight set $T_{D}(t)$ containing a node $t$ 
is the set of nodes from which $t$ is reachable in $D$.
\end{claim}

\Proof We have
\[
 \varrho_{D}(X) + i_{G}(X) 
 = \sum [\varrho_{D}(v):  v\in X] = \widetilde m(X) \geq i_{G}(X) ,
\]
from which $X$ is $m$-tight 
(that is, $\widetilde m(X) =i_{G}(X)$) precisely if $\varrho_{D}(X)=0$, 
and Part (A) follows.
Therefore the smallest $m$-tight set $T_{D}(t)$ containing $t$ is the smallest set 
containing $t$ with in-degree 0, and hence $T_{D}(t)$ is
indeed the set of nodes from which $t$ is reachable in $D$, as stated
in Part (B).  
\finbox \medskip

By Claim \ref{mtight}, $T_{D}(t)$ is easily computable, and hence
Algorithm \ref{algo1} for general M-convex sets can easily be
specialized to graph orientations.  
By applying Corollary \ref{COchardecmin}
to $p:=i_{G}$ and recalling from Claim \ref{mtight} that $C_{i}$ is
$m$-tight, in the present case, precisely if $\varrho_{D}(C_{i})=0$, 
we obtain the following.

\begin{theorem} \label{decminori} 
An orientation of $D$ of $G$ is dec-min if
and only of $\varrho_{D}(C_{i})=0$ for each member $C_{i}$ of the canonical
chain and $\beta_{i} -1\leq \varrho_{D}(v)\leq \beta_{i}$ holds for every
node $v\in S_{i}$ $(i=1,\dots ,q)$.  
\finbox 
\end{theorem}

We remark that the members of the canonical partition computed by our
algorithm for $\odotZ{B_{G}}$ is exactly the non-empty members of the
so-called density decomposition of $G$ introduced by Borradaile et al.~\cite{BMW}.

\subsection{Cheapest dec-min orientations}
\label{cheapdecmin}

In Section \ref{SCmatrdescdmset}, 
we indicated that, given a (linear) cost-function $c$
on the ground-set, there is an algorithm to compute a cheapest dec-min
element of an M-convex set.  
In the special case of dec-min orientations, 
this means that if $c$ is a cost-function on the node-set of $G=(V,E)$, 
then we have an algorithm to compute a dec-min orientation of $G$ for which 
$\sum [c(v)\varrho(v):v\in V]$ is minimum.

But the question remains:  what happens if, instead of a cost-function
on the node-set, we have a cost-function $c$ on $\vec E_{2}$, where
$\vec E_{2}$ arises from $E$ by replacing each element $e=uv$ ($=vu$) of
$E$ by two oppositely oriented arcs $uv$ and $vu$, and we are
interested in finding a cheapest orientation with specified properties?
  (As an orientation of $e$ consists of replacing $e$ by
one of the two arcs $uv$ and $vu$ and the cost of its orientation is,
accordingly, $c(uv)$ or $c(vu)$.  
Therefore we can actually assume that $\min \{c(uv),c(vu)\}=0.$)

It is important to remark that the minimum cost 
in-degree constrained orientation problem itself 
(where the in-degree vector of the
orientation meets specified upper and lower bound constraints, and
decreasing minimality does not play any role) 
can be reduced with a standard technique to a minimum cost
feasible flow problem in a digraph with small integral capacities.
This latter problem is tractable in strongly polynomial time via the
classic min-cost flow algorithm of Ford and Fulkerson 
(that is, we do not need the more sophisticated min-cost flow algorithm of 
Tardos, which is strongly polynomial for an arbitrary capacity).  
Actually, we
shall need a version of this minimum cost orientation problem when
some of the edges are already oriented, and this slight extension is
also tractable by network flows.

Theorem \ref{decminori} implies that 
the problem of finding a cheapest dec-min orientation is equivalent to 
finding a cheapest 
in-degree constrained orientation
by orienting edges connecting $C_{i}$ and $V-C_{i}$ toward $V-C_{i}$
($i=1,\dots ,q$).  
Here the in-degree constraints are given by 
$\beta_{i}-1\leq \varrho_{D}(v)\leq \beta_{i}$ for $v\in S_{i}$ ($i=1,\dots ,q$).

Note that Harada et al.~\cite{HOSY07} provided a direct algorithm for
 the minimum cost version of the so-called semi-matching problem, 
which problem includes the minimum cost dec-min orientation problem. 
For this link, see Section \ref{SCsemimatching}.

\subsection{Orientation with dec-min in-degree vector and dec-min out-degree vector} 
\label{inout}

We mentioned that dec-min and inc-max orientations always concern in-degree vectors.  
As an example to demonstrate the advantage of the general base-polyhedral view, 
we outline here one exception when
in-degree vectors and out-degree vectors play a symmetric role.  
The problem is to characterize undirected graphs admitting an orientation
which is both dec-min with respect to its in-degree vector and dec-min
with respect to its out-degree vector.

For the present purposes, we let $d_{G}$ denote the degree vector of $G$, 
that is, $d_{G}(v)$ is the number of edges incident to $v\in V$.
(This notation differs from the standard set-function meaning of $d_{G}$.)

Let $B_{\rm in}$ denote the convex hull of the in-degree vectors of
orientations of $G$, and $B_{\rm out}$ the convex hull of out-degree
vectors of orientations of $G$. 
 (Earlier $B_{\rm in}$ was denoted by
$B_{G}$ but now we have to deal with both out-degrees and in-degrees.)
As before, $\odottZ{B_{\rm in}}$ 
is the set of in-degree vectors of
orientations of $G$, and $\odottZ{B_{\rm out}}$ is the set of
out-degree vectors of orientations of $G$. 
Let 
$\odottZ{B_{\rm in}\sp{\bullet}}$ 
denote the set of dec-min in-degree vectors of orientations
of $G$, and 
$\odottZ{B_{\rm out}\sp{\bullet}}$ 
the set 
of dec-min out-degree vectors of orientations of $G$.  
By Theorem \ref{matroid-eltolt} both
$\odottZ{B_{\rm in}\sp{\bullet}}$ 
and
$\odottZ{B_{\rm out}\sp{\bullet}}$ 
are matroidal M-convex sets.

Note that the negative of a (matroidal) M-convex set is also a
(matroidal) M-convex set, and the translation of a (matroidal)
M-convex set by an integral vector is also a (matroidal) M-convex set.
Therefore $d_{G}- \odottZ{B_{\rm out}\sp{\bullet}}$ is a matroidal M-convex set.
Clearly, a vector $m_{\rm in}$ is the in-degree vector of an
orientation $D$ of $G$ 
precisely if $d_{G}-m_{\rm in}$ is the out-degree vector of $D$.

We are interested in finding an orientation whose in-degree vector is
dec-min and whose out-degree vector is dec-min.  
This is equivalent to finding a member 
$m_{\rm in}$ of $\odottZ{B_{\rm in}\sp{\bullet}}$ 
for which the vector $m_{\rm out}:= d_{G}-m_{\rm in}$ 
is in the matroidal M-convex set $\odottZ{B_{\rm out}\sp{\bullet}}$.  
But this latter is equivalent to requiring that $m_{\rm in}$ 
is in the M-convex set $d_{G}- \odottZ{B_{\rm out}\sp{\bullet}}$.  
That is, the problem is equivalent to finding an element
of the intersection of the matroidal M-convex sets 
$\odottZ{B_{\rm in}\sp{\bullet}}$ and $d_{G}- \odottZ{B_{\rm out}\sp{\bullet}}$.  
This latter problem can
be solved by Edmonds matroid intersection algorithm \cite{Edmonds79}.




\section{In-degree constrained orientations of graphs}
\label{SCori2}

In this section we first describe an algorithm to find a dec-min
in-degree constrained orientation.  Second, we develop a complete
description of 
the set of dec-min in-degree constrained orientations,  
which gives rise to an algorithm to compute a
cheapest dec-min in-degree constrained orientation.

\subsection{Computing a dec-min in-degree constrained orientation}
\label{SCdecminfg}

Let $f:V\rightarrow {\bf Z}\cup \{-\infty \}$ be a lower bound function 
and $g:V\rightarrow {\bf Z}\cup \{+\infty \}$ an upper bound function 
for which $f\leq g$.  
We are interested in in-degree constrained orientations $D$ of $G$, 
by which we mean that $f(v)\leq \varrho_{D}(v)\leq g(v)$ for every $v\in V$.  
Such an orientation is called {\bf $(f,g)$-bounded}, 
and we assume that $G$ has such an orientation. 
 (By a well-known orientation theorem, such an orientation exists 
if and only if $i_{G}\leq \widetilde g$ and $\widetilde f\leq e_{G}$.

As before, let $\odotZ{B_{G}}$ denote the M-convex set 
of the in-degree vectors of orientations of $G$, 
and let $\odotZ{B\sq_{G}}$ denote the intersection of $\odotZ{B_{G}}$ 
with the integral box $T(f,g)$.  
That is, $\odotZ{B\sq_{G}}$ is the set 
of in-degree vectors of $(f,g)$-bounded orientations of $G$.  
Let $D$ be an $(f,g)$-bounded orientation of $G$
with in-degree vector $m$.  
We denote the smallest tight set containing a node $t$ by 
$T\sq_{D}(t)$ ($=T\sq_m(t))$.  
By applying
Claim \ref{CLTm(u)}
to $\odotZ{B\sq_{G}}$, we obtain that
\begin{equation} \label{(Tsqdt)} 
 T\sq_{D}(t)
 = \begin{cases}
    \{ t  \} & \ \ \hbox{if}\ \ \ \varrho_{D}(t)=f(t), 
   \\
    T_{D}(t) - \{s:  \varrho_{D}(s) = g(s)\} & \ \ \hbox{if}\ \ \ \varrho_{D}(t)>f(t), 
 \end{cases} 
\end{equation}
implying that, in case $\varrho_{D}(t)>f(t)$, the set $T\sq_{D}(t)$ consists of
those nodes $s$ from which $t$ is reachable and 
for which $\varrho_{D}(s)<g(s)$.

Formula \eqref{(Tsqdt)} implies for distinct nodes $s$ and $t$ that the
vector $m':= m + \chi_s - \chi_t$ belongs to $\odotZ{B\sq_{G}}$ precisely
if there is an $st$-dipath 
(i.e. a dipath from $s$ to $t$) for which 
$\varrho_{D}(s)<g(s)$ and $\varrho_{D}(t)>f(t)$.  
We call such a dipath $P$
of $D$ {\bf reversible}.  Note that the dipath $P'$ of $D'$ obtained
by reorienting $P$ is reversible in $D'$.

If $P$ is a reversible $st$-dipath of $D$ for which 
$\varrho_{D}(t)\geq \varrho_{D}(s)+2$, 
then the orientation $D'$ is decreasingly smaller than $D$.  
We call such a dipath {\bf improving}.  
Therefore, reorienting an improving $st$-dipath 
corresponds to a 1-tightening step.  
Hence Theorem \ref{equi.1}  
implies the following extension of Theorem~\ref{pathrev}.

\begin{theorem} \label{Borfg} 
An $(f,g)$-bounded orientation $D$ of $G$ is
dec-min if and only if there is no improving dipath, that is, a dipath
from a node $s$ to a node $t$ for which 
$\varrho_{D}(t) \geq \varrho_{D}(s)+2$, $\varrho_{D}(s)<g(s)$, 
and $\varrho_{D}(t)>f(t)$.  
\finbox 
\end{theorem}

In Section \ref{SCbasicalg} we have presented an algorithm that computes
a dec-min element of an arbitrary M-convex set. 
By specializing it to $\odotZ{B\sq_{G}}$, 
we conclude that 
in order to construct a dec-min $(f,g)$-bounded orientation of $G$, 
one can start with an arbitrary $(f,g)$-bounded orientation, and then
reorient (currently) improving dipaths one by one, 
as long as such a dipath exists.  
As we pointed out after Theorem \ref{pathrev}, 
after at most $\vert E\vert \sp 2$ improving dipath reorientations, 
the algorithm terminates with a dec-min $(f,g)$-bounded orientation of $G$.

\medskip

\paragraph{Canonical chain and essential value-sequence 
for $(f,g)$-bounded orientations} \ 
In Section \ref{SCcompcanchain}, 
we indicated that Algorithm \ref{algo1} can immediately be applied to 
compute the canonical chain, 
the canonical partition, and the essential value-sequence 
belonging to the intersection $\odotZ{B\sq}$ of an arbitrary M-convex set $\odotZ{B}$ 
with an integral box $T(f,g)$.  

This algorithm needs only the original subroutine to compute $T_m(u)$
since, by Claim \ref{CLTm(u)}, 
$T\sq_m(u)$ is easily computable from $T_m(u)$.  
As we indicated above, in the special case of orientations, 
the corresponding sets $T_{D}(t)$ and $T\sq_{D}(t)$ are
immediately computable from $D$.  
Therefore this extended algorithm
can be used in the special case when we are interested in dec-min
$(f,g)$-bounded orientations of $G=(V,E)$.  
The algorithm starts with a dec-min $(f,g)$-bounded orientation $D$ of $G$
and outputs the canonical chain ${\cal C}\sq = \{C\sq_{1},\dots ,C\sq_{q}\}$, 
the canonical partition ${\cal P}\sq = \{S\sq_{1},\dots ,S\sq_{q}\}$, and 
the essential value-sequence $\beta \sq_{1}>\cdots >\beta \sq_{q}$.
In view of Corollary \ref{COchardecmin} 
we also define bounding functions $f\sp{*}$ and $g\sp{*}$ as
\begin{align*}
f\sp{*}(v) & := \beta \sq_{i} -1 \ \hbox{if}\ \ v\in S_{i} 
\qquad (i=1,\dots,q), 
\\
g\sp{*}(v) & := \beta \sq_{i} \phantom{{} -1}  \ \hbox{if}\ \ v\in S_{i} 
\qquad  (i=1,\dots ,q). 
\end{align*}
We say that the small box 
\begin{equation} \label{(fgbounded)} 
T\sp{*}:=T(f\sp{*},g\sp{*}) 
\end{equation}
belongs to $\odotZ{B\sq_{G}}$.  
Clearly, $f\leq f\sp{*}$ and $g\sp{*} \leq g$, 
and hence $T(f\sp{*},g\sp{*})\subseteq T(f,g)$.
In Section \ref{SCoricheapest} 
below we assume that these data are available.

\begin{remark} \rm \label{Tspec} 
A special case of in-degree constrained
orientations is when we have a prescribed subset $T$ of $V$ and a
non-negative function $m_T:T\rightarrow {\bf Z}_+$ serving as an
in-degree specification on $T$, and we are interested in orientations
of $G$ for which $\varrho(v)=m_T(v)$ holds for every $t\in T$.  We
call such an orientation {\bf $T$-specified}. 
This notion will have applications in Section \ref{SCsemimatching}.
\finbox
\end{remark}

\subsection{Cheapest dec-min in-degree constrained orientations}
\label{SCoricheapest}

We are given a cost-function $c$ on the possible orientations of the edges of $G$ 
and our goal is to find a cheapest dec-min $(f,g)$-bounded orientation of $G$.  
This will be done with the help of a purely graphical description 
of the set of all dec-min $(f,g)$-bounded orientations, 
which is given in Theorem \ref{THfgbounded}.

As a preparation, we derive the following claim as an immediate
consequence of our structural result (Theorem \ref{matroid-eltolt}).  
Let $m$ be a dec-min element of an M-convex set $\odotZ{B}$.  Suppose that
$m':=m + \chi_s - \chi_t$ is in $\odotZ{B}$ 
(that is, $s\in T_m(t)$).
Since $m$ is dec-min, $m(t)\leq m(s)+1$.  
If $m(t)=m(s)+1$, then $m'$ and $m$ are value-equivalent and 
hence $m'$ is also a dec-min element of $\odotZ{B}$.  
We say that $m'$ is obtained from $m$ by an {\bf elementary step}.

\begin{claim} \label{smallstep} 
Any dec-min element of $\odotZ{B}$ can be obtained from 
a given dec-min element $m$ by a sequence of at most $\vert S\vert $ elementary steps.  
\end{claim}

\Proof 
By Theorem \ref{matroid-eltolt}, 
the set of dec-min elements of
$\odotZ{B}$ is a matroidal M-convex set in the sense that it
can be obtained from a matroid $M\sp{*}$ 
by translating the incidence vectors of the bases of $M\sp{*}$ 
by the same integral vector $\Delta\sp{*}$.  
A simple property of matroids is that any basis can be obtained from a given basis 
through a sequence of at most $\vert S\vert $ bases 
such that each member of the series can be obtained from 
the preceding one by taking out one element and adding a new one.
The corresponding change in the translated vector is exactly an elementary step.
\finbox 
\medskip

For a subset $E_{0}\subseteq E$ and for an orientation $A_{0}$ of $E_{0}$,
we say that an orientation $D$ of $G$ is {\bf $A_{0}$-extending} if
every element $e$ of $E_{0}$ is oriented in $D$ in the same direction as in $A_{0}$.

\begin{theorem} \label{THfgbounded} 
Let $G=(V,E)$ be an undirected graph admitting an $(f,g)$-bounded orientation.  
Let $\odotZ{B\sq_{G}}$ denote the M-convex set consisting of the in-degree vectors 
of $(f,g)$-bounded orientations of $G$, 
and let $T\sp{*}$ be the small box, belonging to $\odotZ{B\sq_{G}}$, 
as defined in \eqref{(fgbounded)}.
There are a subset $E_{0}$ of $E$ and an orientation $A_{0}$ of $E_{0}$ such
that an $(f,g)$-bounded orientation $D$ of $G$ is a dec-min  $(f,g)$-bounded orientation
if and only if $D$ is an orientation of $G$ extending
$A_{0}$ and the in-degree vector of $D$ belongs to $T\sp{*}$.  
\end{theorem}

\Proof Let $D$ be a dec-min $(f,g)$-bounded orientation of $G$, and
let $m$ denote its in-degree vector.  Consider the canonical chain
${\cal C}\sq = \{C\sq_{1},\dots ,C\sq_{q}\}$, the canonical partition
${\cal P}\sq = \{S\sq_{1},\dots ,S\sq_{q}\}$, and the essential
value-sequence $\beta \sq_{1}>\cdots >\beta \sq_{q}$ belonging to
$\odotZ{B\sq_{G}}$.

For $i\in \{1,\dots ,q\}$, define $$F_{i}:=\{v:  v\in S\sq_{i}, f(v) =
\beta \sq_{i}\}.$$ Since $f(v)\leq m(v)\leq \beta \sq_{i}$ holds for every
element $v$ of $S\sq_{i}$, we obtain that $f(v)=m(v)=\beta \sq_{i}$ for $v
\in F_{i}$.  Note that $F_{i}$ does not depend on $D$.

\begin{claim} \label{nostpath} 
For every $h=1,\dots ,i$, there is no dipath $P$ 
from a node $s\in V-C\sq_{i}$ with $m(s)<g(s)$ 
to a node $t\in S\sq_{h}$ with $\beta \sq_{h}>f(t)$.  
\end{claim}

\Proof Suppose indirectly that there is such a dipath $P$.  
If $m(t) = \beta \sq_{h}$, then $P$ would be an improving dipath which is
impossible since $D$ is dec-min $(f,g)$-bounded.
Therefore $m(t)=\beta \sq_{h}-1$.  
But a property of the canonical partition is
that there is an element $t'$ of $S\sq_{h}-F_{h}$
for which $m(t')=\beta\sq_{h}$ 
and $t\in T\sq_{D}(t')$.  
This means that $t'$ is reachable from
$t$ in $D$, and therefore there is a dipath from $s$ to $t'$ in $D$
which is improving, a contradiction again.  
\finbox

\medskip

We are going to define a chain ${\cal Z}$ of subsets 
$Z_{1} \supseteq Z_{2} \supseteq \cdots \supseteq Z_{q} \ (=\emptyset )$ 
of $V$ with the help
of $D$, and will show that 
this chain actually does
not depend on $D$.
Let 
\begin{equation} \label{(Zidef)} 
 Z_{i}:= \{ t: 
  \mbox{$t$ is reachable in $D$ from a node $s\in V-C\sq_{i}$ with $\varrho_{D}(s)<g(s)$}   \}. 
\end{equation}
Note that $Z_{i-1}\supseteq Z_{i}$ follows from the definition, where 
equality holds precisely if $\varrho_{D}(s) =g(s)$ for each $s\in S_{i}$.

\begin{lemma} \label{Zisame} 
Every dec-min $(f,g)$-bounded orientation defines the same family ${\cal Z}$.  
\end{lemma}

\Proof 
By Claim \ref{smallstep}, it suffices to prove that a single elementary step 
does not change ${\cal Z}$.  
An elementary step in $\odotZ{B\sq_{G}}$ corresponds 
to the reorientation of an $st$-dipath $P$ in $D$ 
where $s,t\in S\sq_{h}-F_{h}$, $m(t)=\beta \sq_{h}$ and $m(s)=\beta\sq_{h}-1$ hold 
for some $h\in \{1,\dots ,q\}$.  
We will show for $i\in \{1,\dots ,q\}$ 
that the reorientation of $P$ does not change $Z_{i}$.

If $h\leq i$, then Claim \ref{nostpath} implies that 
$Z_{i}\cap S\sq_{h} \subseteq F_{h}$.  
Since $\delta_{D}(Z_{i})=0$, the dipath $P$ is disjoint from $Z_{i}$, 
implying that reorienting $P$ does not affect $Z_{i}$.

Suppose now that $h\geq i+1$.  
Since reorienting $P$ results in a dec-min $(f,g)$-bounded orientation $D'$, 
we get that $m(s)+1\leq g(s)$ and hence $s\in Z_{i}-C\sq_{i}$.
Since $\delta_{D}(Z_{i})=0$, we obtain that $t\in Z_{i}-C\sq_{i}$.  
Since 
$\varrho_{D'}(t) = \varrho_{D}(t)-1 < g(t)$
and the set of nodes reachable from $s$ in $D$ is equal to the set of
nodes reachable from $t$ in $D'$, 
it follows that the reorientation of $P$ does not change $Z_{i}$.  
\finbox 
\medskip
\medskip

Let $E_{0}$ consist of those edges of $G$ which connect $Z_{i}$ with
$V-Z_{i}$ for some $i=1,\dots ,q$.  
Let $A_{0}$ denote the orientation of
$E_{0}$ obtained by orienting each edge connecting $Z_{i}$ and $V-Z_{i}$ toward $Z_{i}$.

\begin{lemma} \label{LME0A0}
The subset $E_{0}\subseteq E$ and its orientation $A_{0}$ meet 
the requirements in the theorem.  
\end{lemma}

\Proof 
Consider first an arbitrary dec-min $(f,g)$-bounded orientation $D$ of $G$.  
Then $\delta_{D}(Z_{i})=0$ and hence $D$ extends $A_{0}$.
Moreover, by a basic property of the canonical chain, the in-degree vector 
of $D$ belongs to $T\sp{*}$.

Conversely, let $D$ be an orientation of $G$ extending $A_{0}$ 
whose in-degree vector belongs to $T\sp{*}$, 
that is, 
\begin{equation*}
 \hbox{ $f\sp{*}(v) \leq \varrho_{D}(v) \leq g\sp{*}(v)$ for every $v\in V$.  }\
\end{equation*}
Then $D$ is clearly $(f,g)$-bounded.

\begin{claim} \label{CLnoimproving}
There is no improving dipath in $D$.  
\end{claim}

\Proof 
Suppose, indirectly, that $P$ is an improving $st$-dipath, that is, 
a dipath from $s$ to $t$ such that 
$\varrho_{D}(t) \geq \varrho_{D}(s)+2$, 
$\varrho_{D}(t)>f(t)$, and $\varrho_{D}(s)<g(s)$.  
Suppose that $t$ is in $S\sq_{i}$ for some $i\in \{1,\dots ,q\}$.  
If $s$ is in $S_k\sq$ for some $k\in \{1,\dots ,q\}$, then
\[
 \beta \sq_k-1 \leq \varrho_{D}(s) \leq \varrho_{D}(t) -2\leq \beta
\sq_{i}-2,
\]
that is, $\beta \sq_k < \beta \sq_{i}$, and hence $k>i$,
implying that $s$ is in $V-C\sq_{i}$.  This and $\varrho_{D}(s)<g(s)$
imply that $s$ is in $Z_{i}$.  
Furthermore, $\beta \sq_{i}\geq \varrho_{D}(t)>f(t)$ implies that $t$ is not in $F_{i}$, 
and since 
$S_{i}\sq \cap Z_{i}\subseteq F_{i}$, 
we obtain that $t$ is not in $Z_{i}$.  
On the other hand, 
we must have $t \in Z_{i}$,
since there is a dipath 
from $s \in V-C\sq_{i}$ to $t$ and $\varrho_{D}(s)<g(s)$.
This is a contradiction.
\finbox
\medskip

By proving Claim \ref{CLnoimproving}, we have shown Lemma \ref{LME0A0}.
Thus the proof of the theorem is completed.
\finbox \finboxHere 
\medskip

\paragraph{Algorithm for computing a cheapest dec-min $(f,g)$-bounded orientation} \ 
First we compute a dec-min $(f,g)$-bounded orientation $D$ 
of $G$ with the help of the algorithm outlined in Section \ref{SCdecminfg}.  
Second, by applying the algorithm described in the same section, 
we compute the canonical chain and partition belonging to $\odotZ{B\sq_{G}}$ 
along with the essential value-sequence.  
Once these data are available, the sets $Z_{i}$
($i=1,\dots ,q$) defined in \eqref{(Zidef)} are easily computable.
Lemma \ref{Zisame} ensures that these sets $Z_{i}$ do not depend on the
starting dec-min $(f,g)$-bounded orientation $D$.  
Let $E_{0}$ be the
union of the set of edges connecting some $Z_{i}$ with $V-Z_{i}$, and
define the orientation $A_{0}$ of $E_{0}$ 
by orienting each edge between $Z_{i}$ and $V-Z_{i}$ toward $Z_{i}$.

Theorem \ref{THfgbounded} implies that, once $E_{0}$ and its orientation
$A_{0}$ are available, the problem of computing a cheapest dec-min
$(f,g)$-bounded orientation of $G$ reduces to finding cheapest
in-degree constrained 
(namely, $(f\sp{*},g\sp{*})$-bounded) 
orientation of a mixed graph.  
We indicated already in Section \ref{cheapdecmin} 
that such a problem is easily sovable by the strongly polynomial
min-cost flow algorithm of Ford and Fulkerson in a digraph with identically 1 capacities.

\begin{remark} \rm \label{RMorivar}
In Section \ref{capori} we have considered the capacitated dec-min orientation problem
in the basic case where no in-degree constraints are imposed.
With the technique presented there, 
we can cope with the capacitated, min-cost, in-degree constrained variants as well.
Furthermore, 
the algorithms above can easily be extended, with a slight modification, 
to the case when one is interested in orientations of mixed graphs. 
\finbox
\end{remark}

\subsection{Dec-min $(f,g)$-bounded orientations minimizing the in-degree of $T$} 
\label{SCminm(T)}

One may consider $(f,g)$-bounded orientations of $G$ when the
additional requirement is imposed that the in-degree of a specified subset $T$ 
of nodes be as small as possible.  
We shall show that these orientations of $G$ can be described 
as $(f',g')$-bounded orientations
of a mixed graph arising from $G$ by orienting the edges 
between a certain subset $X_{T}$ 
of nodes and its complement $V-X_{T}$ toward $V-X_{T}$.

It is more comfortable, however, to show the analogous statement for a
general M-convex set $\odotZ{B'}(p)\subseteq {\bf Z}\sp V$ defined by a
(fully) supermodular function $p$ for which $\odotZ{B\sq}:=
\odotZ{B'}(p)\cap T(f,g)$ is non-empty.  (Here, instead of the usual
$S$, we use $V$ to denote the ground-set of the general M-convex set.
We are back at the special case of graph orientations when $p=i_G$.)
We assume that each of $p$, $f$, and $g$ is finite-valued.

Let $p\sq$ denote the unique (fully) supermodular function defining $B\sq$.  
This function can be expressed with the help of $p$, $f$, and $g$, as follows
(see, for example, Theorem 14.3.9 in \cite{Frank-book}):
\begin{equation}
p\sq(Y) = \max \{p(X) + \widetilde f(Y-X) - \widetilde g(X-Y):  \ X\subseteq V \} 
\qquad(Y \subseteq V).
\label{(pfgb)} 
\end{equation}
As $B\sq$ is defined by the supermodular function $p\sq$ (that is, $B\sq=B'(p\sq)$), we have
\begin{equation}
\min \{\widetilde m(T):  m\in \odotZ{B\sq} \}=p\sq(T) .
\label{(mminim)} 
\end{equation}
This implies that the set of elements of $\odotZ{B\sq}$ minimizing $\widetilde m(T)$ 
is itself an M-convex set.  
Namely, it is the set of integral elements of the base-polyhedron arising from
$B\sq$ by taking its face defined by 
$\{m\in B\sq:  \widetilde m(T)=p\sq(T)\}$.  
The next theorem shows how this M-convex set can be described 
in terms of $f$, $g$, and $p$, without referring to $p\sq$.

\begin{theorem}  \label{minmT} 
There is a box $T(f',g')\subseteq T(f,g)$ and a subset 
$X_{T}\subseteq V$ such that an element $m\in \odotZ{B\sq}$
minimizes $\widetilde m(T)$ 
if and only if 
$\widetilde m(X_{T})=p(X_{T})$ and $m\in \odotZ{B}\cap T(f',g')$.  
\end{theorem}

\Proof Let $X_{T}$ be a set maximizing the right-hand side of
\eqref{(pfgb)}.

\begin{claim}  \label{opcri} 
An element $m\in \odotZ{B\sq}$ is a minimizer of
the left-hand side of \eqref{(mminim)} 
if and only if 
the following three optimality criteria hold:
\begin{align*}
& \widetilde m(X_{T})=p(X_{T}),
\\
& v\in T-X_T \quad \mbox{\rm implies} \quad m(v)=f(v) ,
\\ 
& v\in X_T-T \quad \mbox{\rm implies} \quad m(v)=g(v).
\end{align*}
\end{claim}

\Proof 
For any $m\in \odotZ{B\sq}$ and $X\subseteq V$, we have
$\widetilde m(T) = \widetilde m(X) + \widetilde m(T-X) - \widetilde m(X-T) 
 \geq p(X) + \widetilde f(T-X) - \widetilde g(X-T)$.  
Here we have equality if and only if 
$\widetilde m(X) =p(X)$, $\widetilde m(T-X)= \widetilde f(T-X)$, 
and $\widetilde m(X-T)= \widetilde g(X-T)$, 
implying the claim.  
\finbox

\medskip 
Define $f'$ and $g'$ as follows:
\begin{align}
 f'(v) &:= 
 \begin{cases} 
   g(v) & \ \ \hbox{if}\ \ \ v\in X_{T}-T, \cr
   f(v) & \ \ \hbox{if}\ \ \ v\in V- (X_{T}-T), 
 \end{cases} 
\label{(f'def)}
\\
g'(v) &:= 
 \begin{cases} 
  f(v) & \ \ \hbox{if}\ \ \ v\in T-X_{T}, \cr
  g(v) & \ \ \hbox{if}\ \ \ v\in V- (T-X_{T}).  
 \end{cases}
\label{(g'def)} 
\end{align}
The claim implies that $T(f',g')$ and $X_{T}$ meet the requirement of the theorem.  
\finbox \finboxHere \medskip

As the set of elements of $\odotZ{B\sq}$ minimizing $\widetilde m(T)$
is itself an M-convex set, all the algorithms developed earlier can be
applied once we are able to compute set $X_{T}$ occurring in Theorem~\ref{minmT}.  
(By definitions \eqref{(f'def)} and \eqref{(g'def)}, 
$X_{T}$ immediately determines $f'$ and $g'$).

The following straightforward algorithm computes an element 
$m\in \odotZ{B\sq}$ minimizing the left-hand side of \eqref{(mminim)}
and a subset $X_{T}$ maximizing the right-hand side of \eqref{(pfgb)}.  
Start with an arbitrary element $m\in \odotZ{B\sq}$.  
By an improving step we mean the change of 
$m$ to $m':=m + \chi_{s} - \chi_{t}$ for some elements
$s\in V-T, t\in T$ for which $m(s)<g(s)$, $m(t)>f(t)$,  and $s\in T_{m}(t)$, 
where $T_{m}(t)$ is the smallest $m$-tight set 
(with respect to $p$) containing $t$.  
Clearly, $m'\in \odotZ{B\sq}$, and $\widetilde m'(T)=\widetilde m(T)-1$.  
The algorithm applies improving steps as long as possible.  
When no more improving step exists, the set 
$X_{T}:= \cup (T_{m}(t):  t\in T, m(t)>f(t))$ 
meets the three optimality criteria.  
The algorithm is polynomial if $\vert p(X)\vert $ is
bounded by a polynomial of $\vert V\vert $.

By applying Theorem \ref{minmT} to the special case of $p=i_G$, we obtain the following.

\begin{corollary}  \label{grafminmt} 
Let $G=(V,E)$ be a graph admitting an
$(f,g)$-bounded orientation.  There is a box $T(f',g')\subseteq
T(f,g)$ and a subset $X_{T}\subseteq V$ such that an $(f,g)$-bounded
orientation of $G$ minimizes the in-degree of $T$ if and only if $D$
is an $(f',g')$-bounded orientation for which $\varrho _D(X_{T})=0$.
\finbox
\end{corollary}

In this case, the algorithm above to compute $X_{T}$ starts with an
$(f,g)$-bounded orientation $D$ of $G$, whose in-degree vector is
denoted by $m$.  As long as there is an $st$-dipath $P$ with $s\in
V-T, t\in T, m(s)<g(s)$, and $m(t)>f(t)$, reorient $P$.  When no such
a dipath exists anymore, the set $X_{T}$ of nodes from which a node
$t\in T$ with $m(t)>f(t)$ is reachable in $D$, along with the bounding
functions $f'$ and $g'$ defined in \eqref{(f'def)} and in
\eqref{(g'def)}, meet the requirement in the corollary.  \medskip

\paragraph{Minimum cost version} \ 
It follows that, in order to
compute a minimum cost dec-min $(f,g)$-bounded orientation for which
the in-degree of $T$ is minimum, we can apply the algorithm described
in Section~\ref{SCoricheapest} for the mixed graph
obtained from $G$ by orienting each edge between $X_{T}$ and $V-X_{T}$
toward $V-X_{T}$.

\begin{remark} \rm 
Instead of a single subset $T$ of $V$, we may consider a
chain ${\cal T}$ of subsets $T_{1}\subset T_{2}\subset \cdots \subset T_{h}$ of $V$.  
Then $\cal T$ defines a face $B\sq_{\rm face}$ of the base-polyhedron $B\sq$.  
Namely, an element $m$ of $B\sq$ belongs to
$B\sq_{\rm face}$ precisely if 
$\widetilde m(T_{i})=p\sq(T_{i})$
for each $i\in \{1,\dots ,h\}$.  
This implies that the integral elements of
$B\sq_{\rm face}$ simultaneously minimize $\widetilde m(T_{i})$ 
for each $i\in \{1,\dots ,h\}$ (over the elements of $\odotZ{B\sq}$).
Therefore, we can consider $(f,g)$-bounded orientations of $G$ with
the additional requirement that each of the in-degrees of 
$T_{1},T_{2}, \dots ,T_{h}$ 
is (simultaneously) minimum.  Corollary \ref{grafminmt}
can be extended to this case, implying that we have an algorithm to
compute a minimum cost dec-min $(f,g)$-bounded orientation of $G$ that
simultaneously minimizes the in-degree of each member of the chain
$\{T_{1}, T_{2}, \dots , T_{h}$\}. 
\finbox
\end{remark}

\subsection{Application in resource allocation: semi-matchings}
\label{SCsemimatching}

We proved for a general M-convex set $\odotZ{B}$ 
(Corollary \ref{COdecmin=squaresum} and
Theorem \ref{THdecmin=diffsum}) 
that for an element $m$ of $\odotZ{B}$ the following
properties are equivalent:  
\ (A) \ $m$ is dec-min, 
\ (B) \ the square-sum of the components is minimum, 
\ (C) \ the difference-sum of the components of $m$ is minimum.  
Therefore the corresponding equivalences hold 
in the special case of in-degree constrained 
(in particular, $T$-specified)
 orientations of undirected graphs.

As an application of this equivalence, 
we show first how a result of Harvey et al.~\cite{HLLT} 
concerning a resource allocation problem
 (mentioned already in Section \ref{SCresource}) follows immediately.  
They introduced the notion of a semi-matching of 
a simple bipartite graph $G=(S,T;E)$ as a subset $F$ of edges 
for which $d_{F}(t)=1$ holds for every node $t\in T$, 
and solved the problem of finding a semi-matching $F$ for which
$\sum [d_{F}(s) (d_{F}(s)+1):  s\in S]$ is minimum.  
The problem was motivated by practical applications 
in the area of resource allocation in computer science.  
Note that
\begin{align*}
& \sum [d_{F}(s) (d_{F}(s)+1):  s\in S] 
 = \sum [d_{F}(s)\sp{2} :s\in S] + \sum [d_{F}(s):  s \in S] 
\\ &
= \sum [d_{F}(s)\sp{2} :s\in S] + \vert F\vert 
= \sum [d_{F}(s)\sp{2}:s\in S] + \vert T\vert ,
\end{align*}
and therefore the problem of Harvey et al.~is equivalent to
finding a semi-matching $F$ of $G$ that minimizes the square-sum of
degrees in $S$.

By orienting each edge in $F$ toward $S$ and each edge in $E-F$ toward $T$, 
a semi-matching can be identified with the set of arcs directed toward $S$ 
in an orientation of $G=(S,T;E)$ in which the out-degree of every node $t\in T$ is 1 
(that is, $\varrho(t)=d_{G}(t)-1$), and
$d_{F}(s)=\varrho(s)$ for each $s\in S$.  
Since $\varrho(t)$ for $t\in T$ is the same in these orientations, 
it follows that the total sum of $\varrho(v)\sp{2}$ over $S\cup T$ 
is minimized precisely if 
$\sum [\varrho(s)\sp{2}: s\in S] = \sum [d_{F}(s)\sp{2}:s\in S]$ is minimized.
Therefore the semi-matching problem of Harvey et al. is nothing but a
special dec-min $T$-specified orientation problem.
Note that not only semi-matching problems can be managed
with graph orientations, but conversely, 
an orientation of a graph $G=(V,E)$ can also be interpreted 
as a semi-matching of the bipartite
graph obtained from $G$ by subdividing each edge by a new node.  
This implies, for example, that the algorithm of Harvey et al.~to 
compute a semi-matching minimizing 
$\sum [d_F(v)\sp{2}: v\in S]$ 
is able to compute an orientation of a graph $G$ 
for which $\sum [\varrho(v)\sp{2}:v\in S]$ is minimum.  
Furthermore, an orientation of a hypergraph means that 
we assign an element of each hyper-edge $Z$ to $Z$ as its head.  
In this sense, semi-matchings of bipartite graphs and
orientations of hypergraphs are exactly the same.  
Several graph orientation results have been extended to hypergraph orientation, 
for an overview, see, e.g.  \cite{Frank-book}.

Bokal et al.~\cite{BBJ} extended the results to subgraphs of $G$
meeting a more general degree-specification on $T$ when, rather than
the identically 1 function, one imposes an arbitrary
degree-specification $m_{T}$ on $T$ satisfying 
$0\leq m_{T}(t)\leq d_{G}(t)$ $(t\in T)$.  
The same orientation approach applies in this more general setting.  
We may call a subset $F$ of edges an {\bf
$m_{T}$-semi-matching} if $d_{F}(t)= m_{T}(t)$ for each $t\in T$.  
The extended resource allocation problem is to find an $m_{T}$-semi-matching
$F$ that minimizes $\sum [ d_{F}(s)\sp{2}:  s\in S]$.  
This is equivalent to finding a $T$-specified orientation of $G$ 
for which the square-sum
of the in-degrees is minimum and the in-degree specification 
in $t\in T$ is 
$m_{T}'(t):  = d_{G}(t) - m_{T}(t)$.  
Therefore this extended resource allocation problem is equivalent 
to finding a dec-min $T$-specified orientation of $G$.

The same orientation approach, when applied to in-degree constrained
orientations, allows us to extend the $m_{T}$-semi-matching problem when
we have upper and lower bounds imposed on the nodes in $S$.  
This may be a natural requirement in practical applications
where the elements of $S$ correspond to available resources (e.g. computers), 
the elements of $T$ correspond to users, and we are interested in a fair
($=$ dec-min $=$ square-sum minimizer) distribution 
(=$m_{T}$-semi-matchings) of the resources when the load (or burden) of
each resource is requested to meet a specified upper and/or lower bound. 
Note that in the resource allocation framework, the degree $d_{F}(s)$ of
node $s\in S$ may be interpreted as the burden of $s$, and hence a
difference-sum minimizer semi-matching minimizes the total sum of
burden-differences.

Katreni{\v c} and Semani{\v s}in \cite{Katrenic13}
investigated the problem of finding a dec-min \lq maximum
$(f,g)$-semi-matching\rq \ problem where there is a lower-bound
function $f_{T}$ on $T$ and an upper bound function $g_{S}$ on $S$ 
(in the present notation) and one is interested in maximum cardinality
subgraphs of $G$ meeting these bounds.  
They describe an algorithm to compute a dec-min subgraph of this type.  
With the help of the orientation model discussed in Section \ref{SCminm(T)}
(where, besides the in-degree bounds on the nodes, 
the in-degree of a specified subset $T$ was requested to be minimum), 
we have a strongly polynomial algorithm for an extension of the model of
\cite{Katrenic13} when there may be upper and lower bounds on both $S$ and $T$.  
Actually, even the minimum cost version of this problem 
was solved in Section \ref{SCminm(T)}.

In another variation, we also have degree bounds $(f_{S},g_{S})$ on $S$
and $(f_{T},g_{T})$ on $T$, 
but we impose an arbitrary positive integer $\gamma $ 
for the cardinality of $F$.  
We consider degree-constrained subgraphs $(S,T;F)$ of $G$ 
for which $\vert F\vert =\gamma $, and want
to find such a subgraph for which 
$\sum [d_F(s)\sp{2} :  s\in S]$
is minimum.  
(Notice the asymmetric role of $S$ and $T$.)  
This is equivalent to finding an in-degree constrained orientation $D$ of $G$
for which 
$\varrho_D(S)=\gamma $ and 
$\sum [\varrho_D(s)\sp{2}: s\in S]$ 
is minimum.  
Here the corresponding in-degree bound $(f,g)$
on $S$ is the given $(f_{S},g_{S})$ 
while $(f,g)$ on $T$ is defined for $t\in T$ by
\[
 f(t):= d_{G}(t) -g_{T}(t) \quad \mbox{and} \quad g(t):= d_{G}(t) -f_{T}(t).
\]

Let $B$ denote the base-polyhedron spanned by the in-degree vectors of
the degree-constrained orientations of $G$.  Then the restriction of
$B$ to $S$ is a g-polymatroid $Q$.  By intersecting $Q$ with the hyperplane 
$\{x:  \widetilde x(S)=\gamma \}$, 
we obtain an integral base-polyhedron $B_{S}$ in ${\bf R}\sp{S}$, 
and then the elements of $\odotZ{B_{S}}$ are exactly 
the in-degree vectors of the requested orientations restricted to $S$.  
That is, the elements of $\odotZ{B_{S}}$ are the restriction of 
the degree-vectors of the requested subgraphs of $G$ to $S$.  
Since $B_{S}$ is a base-polyhedron, a dec-min element of $\odotZ{B_{S}}$ 
will be a solution to our minimum degree-square sum problem.

\medskip

We briefly indicate that a capacitated version of the semi-matching
problem can also be formulated as 
a dec-min in-degree constrained and capacitated orientation problem 
(cf., Section~\ref{capori} and Remark~\ref{RMorivar}).
Let $G=(S,T;E)$ be again a bipartite graph, $\gamma $ a positive integer, 
and $f_{V}$ and $g_{V}$ integer-valued bounding functions on
$V:=S\cup T$ for which $f_{V}\leq g_{V}$.  
In addition, an integer-valued capacity function $g_E$ 
is also given on the edge-set $E$, and we are interested 
in finding a non-negative integral vector $z:E\rightarrow {\bf Z}_{+}$ 
for which 
$\widetilde z(E)=\gamma $, $z\leq g_E$ and
$f_{V}(v)\leq d_z(v)\leq g_{V}(v)$
for every $v\in V$.  
(Here $d_z(v):=\sum [z(uv):  uv\in E]$.) 
We call such a vector {\bf feasible}.  
The problem is to find a feasible vector $z$ whose degree vector 
restricted to $S$ (that is, the vector $(d_z(s):  s\in S)$ is decreasingly minimal.

By replacing each edge $e$ with $g_E(e)$ parallel edges, it follows
from the uncapacitated case above that the vectors 
$\{(d_z(s):  s\in S) :  \mbox{$z$ is a feasible integral vector} \}$ 
form an M-convex set.  
In this case, however, the basic algorithm is not necessarily polynomial
since the values of $g_E$ may be large.  
Therefore we need the general
strongly polynomial algorithm described in Section \ref{strong.pol}.
 In this case the general Subroutine
\eqref{(ND.routine)} can be realized via max-flow min-cut computations.

\paragraph{Minimum cost dec-min semi-matchings} \ 
Harada et al.~\cite{HOSY07} developed an algorithm to solve the minimum
cost version of the original semi-matching problem of Harvey et al.~\cite{HLLT}.  
As the dec-min in-degree bounded orientation problem covers all the
extensions of semi-matching problems mentioned above, the minimum cost
version of these extensions can also be solved with the strongly
polynomial algorithms developed in Section \ref{SCoricheapest}
for minimum cost dec-min in-degree bounded orientations.

\medskip

We close this section with some historical remarks.  
The problem of Harvey et al.~is closely related to earlier investigations in the
context of minimizing a separable convex function over 
(integral elements of) a base-polyhedron.  
For example, Federgruen and Groenevelt \cite{Fed-Gro86} provided 
a polynomial time algorithm in 1986.  
Hochbaum and Hong \cite{Hochbaum-Hong} in 1995 developed 
a strongly polynomial algorithm;
their proof, however, included a technical gap,
which was fixed by Moriguchi, Shioura, and Tsuchimura \cite{MST} in 2011.  
For an early book on resource allocation, see the
one by Ibaraki and Katoh \cite{Ibaraki-Katoh.88} 
while three more recent surveys are due to 
Katoh and Ibaraki \cite{Katoh-Ibaraki.98} from 1998, 
to Hochbaum \cite{Hoc07} from 2007, and 
to Katoh, Shioura, and Ibaraki \cite{KSI} from 2013.
Algorithmic aspects of minimum degree square-sum problems for general
graphs were discussed by Apollonio and Seb{\H o} \cite{Apollonio-Sebo}.



\section{Orientations of graphs with edge-connectivity requirements}
\label{SCori3}

In this section, we investigate various edge-connectivity requirements
for the orientations of $G$.  
The main motivation behind these
investigations is a conjecture of Borradaile et al.~\cite{BIMOWZ} 
on decreasingly minimal strongly connected orientations.  
Our goal is to prove their conjecture in a more general form.

\subsection{Strongly connected orientations}
\label{SCstrori}

Suppose that $G$ is 2-edge-connected, implying that it has a strong
orientation by a theorem of Robbins \cite{Robbins}.  
We are interested in dec-min strong orientations, meaning 
that the in-degree vector is decreasingly minimal over the strong orientations of $G$.  
This problem of Borradaile et al.~was motivated by a practical application
concerning optimal interval routing schemes.

Analogously to Theorem \ref{pathrev}, they described a natural way to
improve a strong orientation $D$ to another one whose in-degree
vector is decreasingly smaller.  
Suppose that there are two nodes $s$ and $t$ for which 
$\varrho(t)\geq \varrho(s)+2$ 
and there are two edge-disjoint dipaths from $s$ to $t$ in $D$.  
Then reorienting an arbitrary $st$-dipath of $D$ results in another strongly connected
orientation of $D$ which is clearly decreasingly smaller than $D$.
Borradaile et al.~conjectured the truth of the converse (and this
conjecture was the starting point of our investigations).  The next
theorem states that the conjecture is true.

\begin{theorem}   \label{proved.conj} 
A strongly connected orientation $D$ of $G=(V,E)$ is decreasingly minimal 
if and only if 
there are no two arc-disjoint $st$-dipaths in $D$ 
for nodes $s$ and $t$ with $\varrho(t)\geq \varrho(s)+2$.
\end{theorem}

\Proof 
Suppose first that there are nodes $s$ and $t$ with 
$\varrho(t)\geq \varrho(s)+2$ 
such that there are two arc-disjoint $st$-dipaths of $D$.  
Let $P$ be any $st$-dipath in $D$ and let $D'$
denote the digraph arising from $D$ by reorienting $P$.  
Then $D'$ is strongly connected, 
since if it had a node-set $Z$ ($\emptyset \subset Z\subset V$) 
with no entering arcs, then $Z$ must be a $t\ol s$-set
and $P$ enters $Z$ exactly once.  
But then 
$0 = \varrho_{D'}(Z) =\varrho_{D}(Z)-1 \geq 2-1=1$, 
a contradiction.  
Therefore $D'$ is indeed strongly connected
and its in-degree vector is decreasingly 
smaller than that of $D$.

To see the non-trivial part, define a set-function $p_{1}$ as follows:
\begin{equation} 
p_{1}(X):= 
\begin{cases}
 0 & \ \ \hbox{if}\ \ \ X=\emptyset,
\cr
\vert E\vert & \ \ \hbox{if}\ \ \ X=V ,
\cr 
i_{G}(X)+1 & \ \ \hbox{if}\ \ \ \emptyset \subset X\subset V. 
\end{cases} 
\end{equation}
Then $p_{1}$ is crossing supermodular and hence $B_{1}:=B'(p_{1})$ is a base-polyhedron.

\begin{claim}   \label{mstrongori} 
An integral vector $m$ is the in-degree vector of a strong orientation of $G$
if and only if 
$m$ is in $\odottZ{B_{1}}$.
\end{claim}  

\Proof 
If $m$ is the in-degree vector of a strong orientation of $G$,
then $\widetilde m(V)= \vert E\vert = p_{1}(V)$, \ 
$\widetilde m(\emptyset) = 0 = p_{1}(\emptyset )$, and
\[
\widetilde m(Z) = \sum [\varrho(v):v\in Z] 
= \varrho(Z) + i_{G}(Z) \geq 1 + i_{G}(Z)= p_{1}(Z)
\]
for $\emptyset \subset Z\subset V$, that is, $m\in \odottZ{B_{1}}$.

Conversely, let $m\in \odottZ{B_{1}}$.  
Then $m\in B_{G}$ and hence by Claim \ref{oribase}, 
$G$ has an orientation $D$ with in-degree vector $m$.  
We claim that $D$ is strongly connected.  
Indeed,
\[
\varrho(Z)
 = \sum [\varrho(v):v\in Z] - i_{G}(Z) 
 = \widetilde m(Z) - i_{G}(Z) 
 \geq p_{1}(Z)-i_{G}(Z) =1
\]
whenever $\emptyset \subset Z\subset V$.  
\finbox

\begin{claim} \label{strongreori} 
Let $D$ be a strong orientation of $G$ with in-degree vector $m$.  
Let $t$ and $s$ be nodes of $G$.  
The vector $m':=m + \chi_s - \chi_t$ is in $B_{1}$ 
if and only if 
$D$ admits two arc-disjoint dipaths from $s$ to $t$.  
\end{claim}

\Proof 
$m'\in B_{1}$ holds precisely if there is no $t\ol s$-set $X$
which is  $m$-tight
 with respect to $p_{1}$, 
that is, 
$\widetilde m(X) = i_{G}(X)+1$.  
Since $\varrho(Y) + i_{G}(Y) = \sum [\varrho(v):v\in Y]= \widetilde m(Y)$ 
holds for any set $Y\subset V$, the tightness of $X$
(that is, $\widetilde m(X)= i_{G}(X)+1$) 
is equivalent to requiring that $\varrho(X)=1$.  
Therefore $m'\in B_{1}$ 
if and only if 
$\varrho(Y) > 1$
holds for every $t\ol s$-set $Y$, 
which is, by Menger's theorem, equivalent to 
the existence of two arc-disjoint $st$-dipaths of $D$.
\finbox 
\medskip

By Theorem \ref{equi.1}, $m$ is a dec-min element of $\odottZ{B_{1}}$ 
if and only if 
there is no 1-tightening step for $m$.  
By Claim \ref{strongreori} 
this is just equivalent to the condition in the theorem 
that there are no two arc-disjoint $st$-dipaths in $D$ 
for nodes $s$ and $t$ for which $\varrho(t)\geq \varrho(s)+2$.  
\finbox \finboxHere

\medskip

An immediate consequence of Claim \ref{mstrongori} and Theorem \ref{equi.1} 
is the following.

\begin{corollary} 
A strong orientation of $G$ is dec-min 
if and only 
if it is inc-max.  
\finbox 
\end{corollary}

We indicated in Section \ref{SCdecminfg}
how in-degree constrained dec-min orientations can
be managed due to the fact that the intersection of an integral
base-polyhedron $B$ with an integral box $T$ is an integral base-polyhedron.  
The same approach works for degree-constrained strong orientations.  
For example, in this case dec-min and inc-max
again coincide and one can formulate the in-degree constrained
version of Theorem \ref{proved.conj}.
In the next section, we overview more general cases.

\subsection{A counter-example for mixed graphs}
\label{SCstroricntex}

Although Robbins' theorem on strong orientability of undirected graphs
easily extends to mixed graphs, 
as was pointed out by Boesch and Tindell \cite{BT}, 
it is not true anymore that a decreasingly minimal strong orientation 
of a mixed graph is always increasingly maximal.
Actually, one may consider two natural variants.

In the first one, decreasing minimality and increasing maximality
concern the total in-degree of the directed graph obtained from the
initial mixed graph after orienting its undirected edges.  
Let $V=\{a,b,c,d\}$.  
Let $E=\{ab,cd\}$ denote the set of undirected edges
and let 
$A=\{ad,ad,ad,da,da,bc,bc,cb\}$ 
denote the set of directed edges of a mixed graph $M=(V,A+E)$.  
There are two strong orientations of $M$.  
In the first one, the orientations of the elements of $E$ are $ba$ and $dc$, 
in which case the total in-degree vector is $(3,1,3,3)$.  
In the second one, the orientations of the elements of
$E$ are $ab$ and $cd$, in which case the total in-degree vector is $(2,2,2,4)$.  
Now $(3,1,3,3)$ is dec-min while $(2,2,2,4)$ is inc-max.

\medskip

In the second variant, we are interested in the in-degree vector of
the digraph obtained by orienting the originally undirected part $E$.
For this version the counterexample is as follows.  
Let $V=\{a,b,c,d,x,y,u,v\}$.  
Let $E=\{ab,cd, au, au, av, av, dy, dy, bx, bx\}$ 
denote the set of undirected edges and let 
$A=\{ad, da, bc, cb\}$ 
denote the set of directed edges of a mixed graph $M=(V,A+E)$.
The undirected part of $M$ is denoted by $G=(V,E)$.

In any strong orientation of $M=(V,A+E)$, the orientations of the
undirected parallel edge-pairs 
$\{au,au\}, \ \{av,av\}, \ \{dy, dy\}, \ \{bx,bx\}$ 
are oriented oppositely, and hence their contribution to
the in-degrees (in the order of $a, b, c, d, u, v, x, y$) 
is $(2, 1, 0, 1, 1, 1, 1, 1)$.

Therefore there are essentially two distinct strong orientations of $M$.  
In the first one, the undirected edges $ab, cd$ are oriented as $ba, dc$, 
while in the second one the undirected edges $ab, cd$ are oriented as $ab, cd$.  
Hence the in-degree vector of the first strong
orientation corresponding to the orientation of $G$ 
(in the order of $a, b, c, d, u, v, x, y$)
 is \ $(3, 1, 1, 1, 1, 1, 1, 1)$.  
The in-degree vector of second strong orientation corresponding to the
orientation of $G$ is $(2,2,0,2,1,1,1,1)$.  
The first vector is inc-max while the second vector is dec-min.

\medskip

These examples give rise to the question:  
what is behind the phenomenon that 
while dec-min and inc-max coincide for strong orientations of undirected graphs, 
they differ for strong orientations of mixed graph?  
The explanation is, as we pointed out earlier, that
for an M-convex set the two notions coincide and the set of 
in-degree vectors of strong orientations of an 
undirected graph is an M-convex set, 
while the corresponding set for a mixed graph is, in general, 
not an M-convex set.  
It is actually the intersection of two M-convex sets.  
In \cite{Frank-Murota.4},
we shall describe an algorithm for computing 
a dec-min element of the intersection of two M-convex sets.  
\medskip

\subsection{Higher edge-connectivity}
\label{SChighori}

An analogous approach works in a much more general setting.  
We say that a digraph covers a set-function $h$ if $\varrho(X)\geq h(X)$
holds for every set $X\subseteq V$.  
The following result was proved in \cite{FrankJ4}.

\begin{theorem}  \label{hori} 
Let $h$ be a finite-valued, non-negative
crossing supermodular function with $h(\emptyset )=h(V)=0$.  A graph
$G=(V,E)$ has an orientation covering $h$ 
if and only if
\[
 e_{\cal P} \geq \sum_{i=1}\sp q h(V_{i}) \quad \hbox{\rm and}\quad 
 e_{\cal P} \geq \sum_{i=1}\sp q h(V-V_{i}) 
\]
hold for every partition ${\cal P}=\{V_{1},\dots ,V_{q}\}$ of $V$, 
where $e_{\cal P}$ denotes the number of edges connecting
distinct parts of $\cal P$.  
\finbox 
\end{theorem}

This theorem easily implies the classic orientation result of
Nash-Williams \cite{NWir} stating that a graph $G$ has a
$k$-edge-connected orientation precisely if $G$ is $2k$-edge-connected.  
Even more, call a digraph $(k,\ell)$-edge-connected ($\ell\leq k$) 
(with respect to a root-node $r_{0}$) 
if $\varrho(X)\geq k$ whenever 
$\emptyset \subset X\subseteq V-r_{0}$ and 
$\varrho(X)\geq \ell$ whenever $r_{0}\in X\subset V$.  
(By Menger's theorem, 
$(k,\ell)$-edge-connectedness is equivalent to requiring that 
 there are $k$ arc-disjoint dipaths from $r_{0}$ to every node and
there are $\ell$ arc-disjoint dipaths from every node to $r_{0}$.)
Then Theorem \ref{hori} implies:

\begin{theorem}  \label{kellori} 
A graph $G=(V,E)$ has a $(k,\ell)$-edge-connected orientation 
if and only if
\[
 e_{\cal P} \geq k(q-1) + \ell 
\]
holds for every $q$-partite partition $\cal P$ of $V$.  
\finbox 
\end{theorem}  

Note that an even more general special case of Theorem \ref{hori} 
can be formulated to characterize graphs
admitting in-degree constrained and $(k,\ell)$-edge-connected orientations.

It is important to emphasize that however general Theorem \ref{hori} is, 
it does not say anything about strong orientations of mixed graphs.  
In particular, it does not imply the pretty but easily provable theorem 
of Boesch and Tindell \cite{BT}.  
The problem of finding decreasingly minimal in-degree constrained 
$k$-edge-connected orientation of mixed graphs 
will be solved in \cite{Frank-Murota.4}.

The next lemma shows why the set of in-degree vectors of orientations
of $G$ covering the set-function $h$ appearing in Theorem \ref{hori}
is an M-convex set, ensuring in this way the possibility 
of applying
the results on decreasing minimization over M-convex sets to general
graph orientation problems.

\begin{lemma}  \label{hbase} 
An orientation $D$ of $G$ covers $h$ 
if and only if 
its in-degree vector $m$ is in the base-polyhedron $B=B'(p)$,
where $p:= h+i_{G}$ is a crossing supermodular function.  
\end{lemma}

\Proof 
Suppose first that $m$ is the in-degree vector of a digraph covering $h$.  
Then $h(X)\leq \varrho(X)= \widetilde m(X) - i_{G}(X)$
for $X\subset V$ and 
$h(V) = 0 = \varrho(V)= \widetilde m(V) - i_{G}(V)$, 
that is, $m$ is indeed in $B$.

Conversely, suppose that $m\in B$.  
Since $h$ is finite-valued and non-negative, we have 
$\widetilde m(X) \geq p(X)\geq i_{G}(X)$ 
for
$X\subset V$ and $\widetilde m(V)=i_{G}(V)$ 
and hence, by the Orientation lemma, 
there is an orientation $D$ of $G$ with in-degree vector $m$.  
Moreover, this digraph $D$ covers $h$ since 
$\varrho_{D}(X) = \widetilde m(X) - i_{G}(X) \geq p(X) - i_{G}(X) = h(X)$ 
holds for $X\subset V$.  
\finbox 
\medskip

By Lemma \ref{hbase}, Theorem \ref{equi.1} can be applied again to the
general orientation problem covering a non-negative and crossing
supermodular set-function $h$ in the same way as it was applied in the
special case of strong orientation above, but we formulate the result
only for the special case of in-degree constrained and
$k$-edge-connected orientations.

\begin{theorem} 
Let $G=(V,E)$ be an undirected graph endowed with a lower bound function 
$f:V\rightarrow {\bf Z} $ and an upper bound function
$g:V\rightarrow {\bf Z} $ with $f\leq g$.  
A $k$-edge-connected and in-degree constrained orientation $D$ of $G$ 
is decreasingly minimal 
if and only if 
there are no two nodes $s$ and $t$ for which 
$\varrho(t) \geq \varrho(s)+2$,  \
$\varrho(t)>f(t)$, $\varrho(s)<g(s)$, \
and there are $k+1$ arc-disjoint $st$-dipaths.  
\finbox 
\end{theorem}

The theorem can be extended even further to in-degree constrained and
$(k,\ell)$-edge-connected orientations ($\ell\leq k$).

\paragraph{An extension}
We say that a digraph $D=(V,A)$ is $k$-edge-connected 
in a specified subset $S$ of nodes 
if there are $k$-arc-disjoint dipaths in $D$ from
any node of $S$ to any other node of $S$.

By relying on Lemma \ref{hbase}, one can derive the following.

\begin{theorem}  \label{hmbase} 
Let $G=(V,E)$ be an undirected graph with a specified subset $S$ of $V$.  
Let $m_{0}$ be an in-degree specification on $V-S$.  
The set of in-degree vectors of those orientations of $G$
which are $k$-edge-connected in $S$ and in-degree specified in $V-S$
is an M-convex set.  
\finbox 
\end{theorem}

By this theorem, we can determine a decreasingly minimal orientation
among those which are $k$-edge-connected in $S$ and in-degree
specified in $V-S$.  
Even additional in-degree constraints can be imposed on the elements of $S$.

\paragraph{Hypergraph orientation}

Let $H=(V,{\cal E})$ be a hypergraph for which we assume that each
hyperedge has at least 2 nodes.  
Orienting a hyperedge $Z$ means that
we designate an element $z$ of $Z$ as its head-node.  
A hyperedge $Z$ with a designated head-node $z\in Z$ 
is a directed hyperedge denoted by $(Z,z)$.  
Orienting a hypergraph means the operation of orienting
each of its hyperedges.  
We say that a directed hyperedge $(Z,z)$
enters a subset $X$ of nodes if $z\in X$ and $Z-X\not =\emptyset $. 
A directed hypergraph is called $k$-edge-connected 
if the in-degree of every non-empty proper subset of nodes is at least $k$.

The following result was proved in \cite{FrankJ50} 
(see, also Theorem 2.22 in the survey paper \cite{FrankS12}).

\begin{theorem}  
The set of in-degree vectors of $k$-edge-connected and
in-degree constrained orientations of a hypergraph forms an M-convex set.  
\finbox 
\end{theorem}  

Therefore we can apply the general results obtained for decreasing
minimization over base-polyhedra.







\begin{thebibliography}{999}



\bibitem{AHO04}
Ahuja, R. K.,   Hochbaum, D. S.,  Orlin, J. B.:
A cut-based algorithm for the nonlinear dual of the minimum cost network flow problem.
Algorithmica {\bf 39}, 189--208 (2004).



\bibitem{Apollonio-Sebo}
Apollonio, N., Seb{\H o}, A.:
Minconvex factors of prescribed size in graphs. 
SIAM Journal on Discrete Mathematics {\bf 23}, 1297--1310 (2009)




\bibitem{AS18} 
Arnold, B.C., Sarabia, J.M.:
Majorization and the Lorenz Order with Applications 
in Applied Mathematics and Economics,
Springer International Publishing, Cham (2018),
(1st edn., 1987)





\bibitem{BT}
Boesch, F.,  Tindell, R.: 
Robbins's theorem for mixed multigraphs.
American  Mathematical  Monthly {\bf 87},  716--719 (1980)


\bibitem{BBJ}
Bokal, D.,  Bre{\v s}ar, B.,   Jerebic, J.:
A generalization of Hungarian method and Hall's theorem with
applications in wireless sensor networks.
Discrete Applied Mathematics {\bf 160}, 460--470  (2012)



\bibitem{BIMOWZ}
Borradaile, G.,   Iglesias, J.,  Migler, T.,  Ochoa, A., Wilfong, G.,  Zhang, L.:
Egalitarian graph orientations.
Journal of Graph Algorithms and Applications {\bf 21}, 687--708 (2017)




\bibitem{BMW}
Borradaile, G.,   Migler, T.,   Wilfong, G.:
Density decompositions of networks.
In: Cornelius, S., Coronges, K., Gon{\c c}alves, B., Sinatra, R., Vespignani, A. (eds.) 
Complex Networks IX. CompleNet 2018,
pp.~15--26. Springer, Cham (2018)







\bibitem{Edmonds70}
Edmonds, J.:
Submodular functions, matroids and certain polyhedra.
In: Guy, R., Hanani, H., Sauer, N., Sch\"onheim, J. (eds.)
Combinatorial Structures and Their Applications,
pp.~69--87. Gordon and Breach, New York (1970)



\bibitem{Edmonds73}
Edmonds, J.:
Edge-disjoint branchings.
In: Rustin, B. (ed.)
Combinatorial Algorithms, 
pp.~91--96.
Academic  Press, New York (1973)


\bibitem{Edmonds79}
Edmonds, J.: 
Matroid intersection.
Annals of Discrete Mathematics {\bf 14}, 39--49 (1979)





\bibitem{Fed-Gro86}
Federgruen, A., Groenevelt, H.:
The greedy procedure for resource allocation problems:
necessary and sufficient conditions for optimality.
Operations Research {\bf 34}, 909--918 (1986) 



\bibitem{Ford-Fulkerson}
Ford, L. R., Jr., Fulkerson, D.R.: 
Flows in Networks.
Princeton University Press, Princeton (1962)




\bibitem{FrankJ4} 
Frank, A.:
On the orientation of graphs.
Journal of Combinatorial Theory, Ser.~B {\bf 28}, 251--261 (1980) 





\bibitem{FrankJ8} 
Frank, A.:
An algorithm for submodular functions on graphs.
Annals of Discrete Mathematics {\bf 16}, 97--120 (1982)



\bibitem{FrankP6}
Frank, A.:
Generalized polymatroids.
In: Hajnal, A., Lov\'{a}sz, L., S{\'o}s, V.T. (eds.) 
Finite and Infinite Sets
(Colloquia  Mathematica Societatis  J{\'a}nos Bolyai {\bf 37}),
pp.~285--294.
North-Holland, Amsterdam (1984)





\bibitem{FrankJ23}
Frank, A.:
Augmenting graphs to meet edge-connectivity requirements.
SIAM Journal on Discrete Mathematics {\bf  5}, 22--53 (1992) 




\bibitem{Frank-book}
Frank, A.:
Connections in Combinatorial Optimization.
Oxford University Press, Oxford (2011)




\bibitem{FrankP2}
Frank,  A.,  Gy\'arf\'as, A.:
How to orient the edges of a graph.
In: Combinatorics
(Colloquia  Mathematica Societatis  J{\'a}nos Bolyai {\bf 18}),
pp.~353--364. 
North-Holland, Amsterdam (1978)






\bibitem{FrankJ31}
Frank, A.,    Jord{\'a}n, T.:
Minimal edge-coverings of pairs of sets.
Journal of Combinatorial Theory, Ser.~B {\bf 65}, 73--110 (1995) 




\bibitem{FrankS12}
Frank, A.,    Kir\'aly, T.:
A survey on covering supermodular functions.
In: Cook, W., Lov{\'a}sz, L., Vygen, J. (eds.)
Research Trends in Combinatorial Optimization,  pp.~87--126. 
Springer, Berlin (2009)



\bibitem{FrankJ50}
Frank, A.,    Kir\'aly, T.,  Kir\'aly, Z.:
On the orientation of graphs and hypergraphs.
Discrete Applied Mathematics {\bf 131}, 385--400 (2003)



\bibitem{Frank-Murota.2}
Frank, A.,  Murota, K.:
Discrete decreasing minimization, Part II:
Views from discrete convex analysis,
arXiv: 1808.08477 (August 2018)


\bibitem{Frank-Murota.3}
Frank, A.,  Murota, K.:
Discrete decreasing minimization, Part III: Network flows,
arXiv: 1907.02673 (July 2019)




\bibitem{Frank-Murota.4}
Frank, A.,  Murota, K.:
Discrete decreasing minimization, Part IV:
Submodular flows and the intersection of two base-polyhedra,
in preparation



\bibitem{Frank-Murota.5}
Frank, A.,  Murota, K.:
Discrete decreasing minimization, Part V:
Weighted cases,
in preparation%



\bibitem{FrankJ17}
Frank, A., Tardos, \'{E}.:
Generalized polymatroids and submodular flows.
Mathematical Programming {\bf 42}, 489--563 (1988)


\bibitem{Fujishige80}
Fujishige, S.:
Lexicographically optimal base of a polymatroid 
with respect to a weight vector.
Mathematics of Operations Research {\bf 5}, 186--196 (1980)


\bibitem{Fujishige84e}
Fujishige, S.:
Structure of polyhedra determined by submodular functions 
on crossing families.
Mathematical Programming {\bf 29}, 125--141 (1984)



\bibitem{Fujishigebook}
Fujishige, S.:
Submodular Functions and Optimization,
2nd edn.
Annals of Discrete Mathematics {\bf 58},
Elsevier, Amsterdam  (2005)





\bibitem{GZSS}
Ghodsi, A.,   Zaharia, M.,    Shenker, S.,   Stoica, I.: 
Choosy: Max-min fair sharing for datacenter jobs with constraints.
In: EuroSys '13 Proceedings of the 8th ACM European Conference on Computer Systems,
pp.~365--378, ACM New York, NY (2013)




\bibitem{Groenevelt91}
Groenevelt, H.:
Two algorithms for maximizing a separable concave function
over a polymatroid feasible region.
European Journal of Operational Research {\bf 54}, 227--236 (1991)





\bibitem{Hakimi}
Hakimi,  S.L.:
On the degrees of the vertices of a directed graph.
Journal of The Franklin Institute {\bf 279}, 290--308  (1965)


\bibitem{HOSY07}
Harada, Y.,   Ono, H.,   Sadakane, K.,  Yamashita, M.:
Optimal balanced semi-matchings for weighted bipartite graphs.
IPSJ Digital Courier {\bf 3},  693--702  (2007)



\bibitem{HLLT}
Harvey, N.J.A.,  Ladner, R.E.,   Lov{\'a}sz, L.,   Tamir, T.:
Semi-matchings for bipartite graphs and load balancing.
Journal of Algorithms {\bf 59}, 53--78 (2006)


\bibitem{Hoc02}
Hochbaum, D.S.:
Solving integer programs over monotone inequalities in three variables: 
A framework for half integrality and good approximations.
European Journal of  Operational Research {\bf 140}, 291--321 (2002)


\bibitem{Hoc07}
Hochbaum, D.S.:
Complexity and algorithms for nonlinear optimization problems.
Annals of Operations Research {\bf 153}, 257--296 (2007)





\bibitem{Hochbaum-Hong}
Hochbaum, D.S., Hong, S.-P.:
About strongly polynomial time algorithms for
quadratic optimization over submodular constraints.
Mathematical Programming {\bf 69}, 269--309 (1995)



\bibitem{Hoffman60}
Hoffman, A.J.:
Some recent applications of the theory of linear inequalities to extremal combinatorial analysis.
In:  Bellman, R.,    Hall, M., Jr. (eds.)
Combinatorial Analysis
(Proceedings of the Symposia of Applied Mathematics {\bf 10}) 
pp.~113--127.
American Mathematical Society, Providence, Rhode Island (1960) 



\bibitem{Ibaraki-Katoh.88}
Ibaraki, T., Katoh, N.:
Resource Allocation Problems: Algorithmic Approaches.
MIT Press, Boston (1988)




\bibitem{IFF01}
Iwata, S., Fleischer, L., Fujishige, S.:
A combinatorial, strongly polynomial-time algorithm for minimizing submodular functions. 
Journal of the ACM {\bf 48}, 761--777 (2001)









\bibitem{Katoh-Ibaraki.98}
Katoh, N., Ibaraki, T.:
Resource allocation problems.
In: Du, D.-Z., Pardalos, P.M. (eds.) 
Handbook of Combinatorial Optimization, Vol.2,
pp.~159--260. Kluwer Academic Publishers, Boston (1998)




\bibitem{KSI}
Katoh, N.,  Shioura, A.,  Ibaraki, T.:
Resource allocation problems.
In: Pardalos, P.M., Du, D.-Z.,   Graham, R.L. (eds.) 
Handbook of Combinatorial Optimization, 2nd ed., 
Vol. 5, 
pp.~2897-2988, Springer, Berlin (2013) 



\bibitem{Katrenic13}
Katreni{\v c}, J., Semani{\v s}in, G.:
Maximum semi-matching problem in bipartite graphs.
Discussiones Mathematicae, Graph Theory {\bf 33}, 559--569 (2013)




\bibitem{Levin-Onn}
Levin, A.,  Onn, S.:
Shifted matroid optimization. 
Operations Research Letters {\bf 44}, 535--539 (2016)



\bibitem{Lovasz83} 
Lov\'{a}sz, L.:
Submodular functions and convexity.
In: ~Bachem, A., ~Gr\"otschel, M.,~Korte, B. (eds.)
Mathematical Programming---The State of the Art,
pp.~235--257. 
Springer, Berlin (1983) 




\bibitem{MOA11}
Marshall, A.W.,   Olkin, I.,  Arnold, B.C.:
Inequalities: Theory of Majorization and Its Applications,
2nd edn. 
Springer, New York (2011),
(1st edn., 1979)


\bibitem{Mar78} 
Maruyama, F.: 
A unified study on problems in information theory via polymatroids.
Graduation Thesis, University of Tokyo, Japan, 1978. (In Japanese.)


\bibitem{Megiddo74}
Megiddo, N.:
Optimal flows in networks with multiple sources and sinks.
Mathematical Programming {\bf 7}, 97--107 (1974)



\bibitem{Megiddo77}
Megiddo, N.:
A good algorithm for lexicographically optimal flows in multi-terminal networks.
Bulletin of the American Mathematical Society {\bf 83}, 407--409 (1977)




\bibitem{MST}
Moriguchi, S., Shioura, A., Tsuchimura, N.:
M-convex function minimization by continuous relaxation 
approach---Proximity theorem and algorithm. 
SIAM Journal on Optimization {\bf 21}, 633--668 (2011)




\bibitem{Murota98a}
Murota, K.:
Discrete convex analysis. 
Mathematical Programming  {\bf 83}, 313--371 (1998)



\bibitem{Murota03}
Murota, K.:
Discrete Convex Analysis.
Society for Industrial and Applied Mathematics, Philadelphia (2003)





\bibitem{Nag07} 
Nagano, K.:
On convex minimization over base polytopes. 
In: Fischetti, M.,   Williamson, D.P. (eds.): 
Integer Programming and Combinatorial Optimization.
Lecture Notes in Computer Science, vol.~4513, pp.~252--266 (2007)


\bibitem{NWir}
Nash-Williams,  C.St.J.A.:
On orientations, connectivity and odd vertex pairings in finite graphs.
Canadian Journal of Mathematics {\bf 12},  555--567  (1960)



\bibitem{Orlin09}
Orlin, J.B.:
A faster strongly polynomial time algorithm for 
submodular function minimization.
Mathematical Programming, Series A {\bf 118}, 237--251 (2009)






\bibitem{Radzik}
Radzik, T.:
Fractional combinatorial optimization.
In: Pardalos, P.M., Du, D.-Z.,    Graham, R.L.
(eds.)
Handbook of Combinatorial Optimization, 2nd edn.,
pp.~1311--1355.
Springer Science+Business Media, New York (2013)







\bibitem{Robbins}
Robbins, H.E.:
A theorem on graphs with an application to a problem of traffic control.
American Mathematical Monthly {\bf 46},  281--283 (1939)




\bibitem{Schrijver2000}
Schrijver, A.:
A combinatorial algorithm minimizing submodular functions
in strongly polynomial time.
Journal of Combinatorial Theory, 
Series B {\bf 80}, 346--355 (2000) 




\bibitem{Tamir95}
Tamir, A.: 
Least majorized elements and generalized polymatroids.
Mathematics of Operations Research {\bf 20}, 583--589 (1995)



\end{thebibliography}
\end{document}